\documentclass[12pt]{article}

\usepackage[margin=1.0in]{geometry}

\usepackage[utf8]{inputenc} 
\usepackage[T1]{fontenc}    
\usepackage{url}            
\usepackage{booktabs}       
\usepackage{amsfonts}       
\usepackage{nicefrac}       
\usepackage{microtype}      
\usepackage{lipsum}

\RequirePackage{amsthm,amsmath}

\RequirePackage[authoryear]{natbib}  

\usepackage{hyperref}
\hypersetup{
  colorlinks   = true,    
  urlcolor     = blue,    
  linkcolor    = blue,    
  citecolor    = blue     
}

\usepackage{dsfont}
\usepackage{amsfonts}
\usepackage{amssymb}
\usepackage{graphics,latexsym,epsfig,bm,textcomp,url,xcolor}

\usepackage{booktabs} 

\usepackage{algorithm}
\usepackage{algorithmic}

\usepackage{graphicx}

\newtheorem{theorem}{Theorem}
\newtheorem{lemma}{Lemma}
\newtheorem{corollary}{Corollary}
\newtheorem{pro}{Proposition}
\newtheorem{defi}{Definition}
\newtheorem{remark}{Remark}

\newcommand{\ind}[1]{\mathds{1}{\Big\{ {#1} \Big\} }}   
\newcommand{\indi}[1]{\mathds{1}{\big\{ {#1} \big\} }}   
\newcommand{\bw}{\boldsymbol{w}}    

\usepackage{authblk}
\title{Distributed Nearest Neighbor Classification}
\author[1]{Jiexin Duan\thanks{PhD candidate, Department of Statistics, Purdue University, West Lafayette, IN 47906. (Email: duan32@purdue.edu).}}
\author[2]{Xingye Qiao\thanks{Associate professor, Binghamton University, State University of New York, Binghamton, NY 13902. (Email: qiao@math.binghamton.edu).}}
\author[1]{Guang Cheng\thanks{Corresponding Author. Professor, Department of Statistics, Purdue University, West Lafayette, IN 47906. (Email: chengg@purdue.edu). Research Sponsored by NSF DMS-1712907, DMS-1811812, DMS1821183, and Office of Naval Research, (ONR N00014-18-2759). }}
\affil[1]{Department of Statistics, Purdue University}
\affil[2]{Department of Mathematical Sciences, Binghamton University}

\begin{document}
\maketitle

\begin{abstract}
Nearest neighbor is a popular nonparametric method for classification and regression with many appealing properties. In the big data era, the sheer volume and spatial/temporal disparity of big data may prohibit centrally processing and storing the data. This has imposed considerable hurdle for nearest neighbor predictions since the entire training data must be memorized. One effective way to overcome this issue is the distributed learning framework. Through majority voting, the distributed nearest neighbor classifier achieves the same rate of convergence as its oracle version in terms of both the regret and instability, up to a multiplicative constant that depends solely on the data dimension. The multiplicative difference can be eliminated by replacing majority voting with the weighted voting scheme. In addition, we provide sharp theoretical upper bounds of the number of subsamples in order for the distributed nearest neighbor classifier to reach the optimal convergence rate. It is interesting to note that the weighted voting scheme allows a larger number of subsamples than the majority voting one. Our findings are supported by numerical studies using both simulated and real data sets.
\end{abstract}



\newpage

\section{Introduction}
\label{sec:intro}
Classification is one of the pillars of statistical learning. The nearest neighbor classifier is among the conceptually simplest and most popular of all classification methods. It is a memory-intensive method in that the entire training data must be memorized to make a prediction. Instead of spending long time to learn a simple rule from the training data, the nearest neighbor classifier defers the computational burden to the prediction stage. The asymptotic properties of the nearest neighbor classification have been studied in \cite{FH51,CH67,DGKL94,S12,CD14,gottlieb2014near,gkm16}, among others.

In the era of big data, due to the unprecedented growth of the sample size  and dimension of the data, denoted as $N$ and $d$, the time and space complexities of nearest neighbor methods are huge. A naive algorithm for $k$-nearest neighbor ($k$NN) classification would compute distances from each query point to all the $N$ training data points, sort the distances, and identify the $k$ smallest distances. Using this naive approach, a single search query has running time between $O(N)$ to $O(N\log(N))$ depending on how efficient the sort method one uses \citep{hoare1961algorithm}. With a large number of training data, having each search query take $O(N)$ time can be prohibitively expensive. The space complexity for storing the training data is $Nd$. For very large data, $k$NN cannot even be conducted on a single machine if the sample size exceeds the memory of the machine. 

A few proposals in the computer science community \citep{anchalia2014k, maillo2015mapreduce} suggested to split the data into multiple local machines and leverage a distributed computing environment, for example, Apache Hadoop that uses the MapReduce paradigm, to process high volume data. However, except for organizations that are known for their ability to exploit large data assets, such as high-tech corporations or research institutions, these distributed computing environments are often not very user friendly or accessible to many average users. The irony is that $k$NN is meant to be a simple yet powerful approach that even a layman can comprehend. There are a group of approximate nearest neighbor search algorithms, such as the locality-sensitive hashing methods \citep{indyk1998approximate}, designed for processing large data sets. However, their implications to the learning performance is less known (with few exceptions such as \citet{gottlieb2014efficient}). Additionally, there are a few techniques that have been empirically used very well, such as the random projection or partition trees \citep{kleinberg1997two,liu2005investigation,dasgupta2008random} and boundary trees \citep{mathy2015boundary}. All  these trees are used for approximate nearest neighbor search. Some theory has shown that \citep{dasgupta2013randomized} a simplification of the random projection tree may be used for exact nearest neighbor search. Currently they lack more theoretical guarantees. Likewise, \cite{ML14} proposed scalable nearest neighbor algorithms for high dimensional data, with little statistical guarantee available for the classification performance of the proposed algorithms.


In this article, we study the nearest neighbor classification in a distributed learning framework designed to alleviate the space and time complexity issues in the big data setting. A few recent works for learning tasks like regression and principal component analysis have fallen under the distributed learning umbrella, e.g., \cite{zhang2013divide, chen2014split, battey2015distributed, zhao2016partially, fan2017distributed, lee2017communication, SC17}. However, distributed {\em classification} is much less understood.


We first propose a Distributed Nearest Neighbor classifier via Majority voting (M-DNN), that is, data are distributed to $s$ subsamples (with subsample size $n\ll N$), nearest neighbor predictions are made at the subsamples, and they are aggregated to cast a single predict using majority voting. This framework can substantially reduce the time and space complexity and is easy to generalize to other classifiers. When computing is done in parallel among all subsamples, the space complexity of M-DNN is $dn$ at each subsample, much lower than $dN$, and the time complexity is reduced to $n\log(n)$ from $N\log(N)$.

M-DNN can actually work without a parallel computing environment. For example, in a multi-cohort medical study, it is fairly common for multiple institutes to collect sensitive patient data separately. Regulations and privacy issues make it impossible to gather all the patient data at a centralized location. Using a distributed learning idea, to predict the class label for a new instant, one can make a prediction at each institute locally, then combine the results to reach a final prediction, without revealing the information of the training data. The theoretical study we conduct in this work can help to understand the learning performance in this scenario.

Our foremost contribution is a proof of an asymptotic expansion form of the regret of the M-DNN method. This proof is not a straightforward extension from \cite{S12}. Specifically, we need an extra normal approximation by the uniform Berry-Esseen theorem \citep{lehmann2004elements} as $s$ diverges, which leads to residual terms that are bounded in a nontrivial way. With carefully chosen weights, the regret of M-DNN can achieve the optimal convergence rate of regret. The only loss is a multiplicative constant which depends on the data dimension only, caused by a Taylor expansion of the normal cumulative distribution function at 0; see Remark \ref{rem:thm1_1}. It was revealed that such a loss is due to the use of majority voting, and hence is dubbed as the majority voting constant.

To eliminate the majority voting constant loss in the regret, we consider an alternative weighting scheme called weighted voting. \cite{dietterich2000ensemble} introduced weighted voting to construct good ensembles of classifiers, which performed better than the base classifiers. \cite{kuncheva2014weighted} proposed a probabilistic framework for classifier ensemble by weighted voting, and conducted some simulations to show that weighted voting outperformed majority voting under certain conditions. We prove that Distributed Nearest Neighbor classifier with Weighted voting (W-DNN) achieve the exact same asymptotic regret as the oracle optimal weighted nearest neighbor (OWNN) \citep{S12}, even with the same multiplicative constant. (We define an ``oracle'' classifier as the one trained on a single machine that has infinite storage and computing power which has access to the entire training data.) Moreover, the time and space complexity of W-DNN are similar to M-DNN.


Our third contribution is to find {\em sharp} upper bounds for the number of subsamples in the proposed M-DNN and W-DNN classifiers, namely, $N^{2/(d+4)}$ and $N^{4/(d+4)}$ respectively, that we can afford in order for them to achieve the optimal convergence rate.  In practice, these upper bounds can provide some theoretical guidance on choosing the number of machines deployed for DNN. 

Much of our findings in this paper is motivated by the study of optimal weighted nearest neighbor (OWNN) by \citet{S12}. The proposed DNN method seems to resemble the bagged nearest neighbor (BNN) method which was closely studied by \citet{hall2005properties,biau2008consistency, biau2010rate}, except that DNN uses data divisions and BNN is based on bootstrap (sub)sampling. In addition, the two methods are fundamentally different in terms of the goals. The DNN method is proposed to deal with big data that cannot be processed by a single machine, while bagging's goal is to improve the stability \citep{buhlmann2002analyzing,Y13} and classification accuracy. 

Lastly, we compare the DNN method with bagging in terms of classification instability (CIS) \citep{SQC16}. For simplicity, we focus on the BNN classifier that applies $1$-NN classifier to each bootstrapped subsample and returns the final classification by majority voting. (Note that \cite{S12} previously cast the ``infinite simulations'' version of BNN as a special case of WNN classifier.) Specifically, it is found that W-DNN is more stable than BNN, while M-DNN is less stable, despite that all three share the same convergence rates of asymptotic regret and CIS.  



The rest of this article is organized as follows. We derive the asymptotic expansion form for the regret of M-DNN in Section \ref{sec:M-DNN}, followed by some asymptotic comparisons between M-DNN and the oracle WNN. In Section \ref{sec:W-DNN},  we shift the focus to W-DNN. We compare the two DNN methods with bagging in Section \ref{sec:cis}. In Section \ref{sec:exp}, we conduct some numerical studies to illustrate the theoretical results. Some concluding remarks are given in Section \ref{sec:conclusion}.

\section{Preliminaries}\label{sec:pre}
Let $(X, Y)\in \mathcal{R}\times\{0,1\}$ be a random couple with a joint distribution $P$ where $\mathcal{R}\subset\mathbb{R}^d$. We regard $X$ as a $d$-dimensional vector of features for an object and $Y$ as a label indicating that the object belongs to one of two classes. Denote the prior probability as $\pi_j:={\mathbb P}(Y=j)$ and the conditional distribution of $X$ given $Y=j$ as $P_j$ for $j=0,1$. Hence, the marginal distribution of $X$ is $\bar{P} =\pi_1 P_1 + (1-\pi_1) P_0$. For a classifier $\phi$: ${\mathbb R}^d \rightarrow  \lbrace 0, 1 \rbrace$, its risk is defined as $$
R(\phi)={\mathbb P}(\phi(X) \neq Y),$$ and is minimized by the Bayes classifier $\phi^{\ast}(x)=\indi{ \eta(x) \geq 1/2}$, where $$\eta(x)={\mathbb P}(Y=1|X=x)$$ is called the regression function. The corresponding risk $R(\phi^*)$ is thus called the Bayes risk. In practice, a classification procedure $\Psi$ is applied to a training data set ${\cal D}:= \lbrace X_i,Y_i\rbrace_{i=1}^n$ to produce a classifier $\widehat{\phi}_n = \Psi({\cal D})$, with the corresponding risk ${\mathbb E}_{\cal D}[ R(\widehat{\phi}_{n}) ]$. Here, ${\mathbb E}_{\cal D}$ denotes the expectation with respect to the distribution
of $\cal D$. The regret of $\Psi$ is defined as: 
$$
{\rm Regret}(\Psi)={\mathbb E}_{\cal D} [ R(\widehat{\phi}_{n})] - R(\phi^{\ast}).
$$

We next introduce a general class of weighted nearest neighbor  (WNN) classifiers. For any query point $x$, let $(X_{(1)},Y_{(1)})$, $(X_{(2)},Y_{(2)})$, $\ldots$ $(X_{(n)},Y_{(n)})$ be the sequence of observations with ascending distance to $x$, and denote $w_{ni}$ as the (non-negative) weight assigned to the $i$-th neighbor of $x$. 
The WNN classifier is defined as
\begin{equation*}
\widehat{\phi}_{n, \bw_n}(x)=\ind{ \widehat{S}_{n,\bw_n}(x) \geq 1/2 }, \;\; {\rm s.t.} \;\;  \sum^n_{i=1} w_{ni}=1,
\end{equation*} 
where $\bw_n$ denotes the weight vector and $\widehat{S}_{n,\bw_n}(x)=\sum^n_{i=1} w_{ni} Y_{(i)}$. When $w_{ni}=k^{-1}$ for $1 \leq i \leq k$, WNN reduces to the $k$-nearest neighbor ($k$NN) classifier, denoted as $\widehat{\phi}_{n,k}(x)$. 

As an important starting point, Proposition \ref{thm:WNN_re} below \citep{S12} provides the asymptotic regret of WNN. Since it is an existing result, we postpone the descriptions of the assumptions and definitions of the constants therein to Appendix.  
\begin{pro}
\label{thm:WNN_re} (Asymptotic Regret for WNN, \citep{S12}) 
Assuming (A1)--(A4) stated in Appendix \ref{sec:assumptions}, we have for each $\beta\in (0,1/2)$,
\begin{equation}
{\rm Regret}(\widehat{\phi}_{n,\bw_n}(x)) \to \Big\{B_1 \sum_{i=1}^n w_{ni}^2 + B_2 \Big (\sum_{i=1}^n \frac{\alpha_i w_{ni}}{n^{2/d}}\Big)^2 \Big\},\;\;{\rm as}\;n\to\infty,\label{eq:WNN_re}
\end{equation}
uniformly for $\bw_n\in W_{n,\beta}$, where $\alpha_i=i^{1+\frac{2}{d}}-(i-1)^{1+\frac{2}{d}}$, constants $B_1, B_2$ are defined in Appendix \ref{sec:defwnb} and $W_{n,\beta}$\footnote{In the case of $k$NN, it means $k$ satisfies $\max(n^{\beta}, (\log n)^2) \le k \le \min(n^{(1-\beta d /4)}, n^{1-\beta})$.} is defined in Appendix \ref{sec:defwnb}.
\end{pro}

We remark that the first term in \eqref{eq:WNN_re} can be viewed as the variance component of regret, while the second term as the squared bias component. By minimizing the asymptotic regret (\ref{eq:WNN_re}) over weights, \cite{S12} has obtained the so-called optimal weighted nearest neighbor (OWNN) classifier.

\section{DNN Classifier via Majority Voting}\label{sec:M-DNN}
In this section, we introduce the first type of distributed WNN based on majority voting, denoted as M-DNN, and then derive its asymptotic regret. A simple comparison reveals that the difference between the regrets of M-DNN and its oracle counterpart is only at the multiplicative constant level, given the weights in local classifiers are carefully chosen.

The main idea of M-DNN is straightforward: 
\begin{itemize}
    \item randomly partition a massive data set $\mathcal D$ with size $N$ into $s$ subsamples;
    \item a local WNN classifier is obtained based on each subsample;
    \item the final classifier is an outcome of majority voting over $s$ classifiers.
\end{itemize}
For simplicity, we assume equal subsample size, say $n:=N/s$ (denote $s=N^{\gamma}$ and $n=N^{1-\gamma}$). The same local weights will be applied to the subsamples to form local WNN classifiers, which are aggregated as in (\ref{eq:M-DNN}) below.

This is summarized in Algorithm \ref{algo:DC_MV}.
\begin{algorithm}
	\caption{DNN via Majority Voting (M-DNN)}
	\label{algo:DC_MV}
	\begin{algorithmic}[1]
	\renewcommand{\algorithmicrequire}{\textbf{Input:}}
	\renewcommand{\algorithmicensure}{\textbf{Output:}}
	\REQUIRE Data set ${\cal D}$, number of partitions $s$, local weight vector $\bw_n$ and query point $x$.
	\ENSURE M-DNN.
	\STATE Randomly split ${\cal D}$ into $s$ subsamples with equal size $n$.
	\FOR {$j = 1$ to $s$}
		\STATE Obtain the WNN classifier $\widehat{\phi}_{n,\bw_n}^{(j)}(x)$ based on the $j$-th subsample.  
	\ENDFOR
	\STATE Majority voting of all classification outcomes $\widehat{\phi}_{n,\bw_n}^{(j)}(x)$: 
\begin{equation}
\label{eq:M-DNN}
\widehat{\phi}_{n,s,\bw_n}^{M}(x)=\ind{\frac{1}{s}\sum_{j=1}^s\widehat{\phi}_{n,\bw_n}^{(j)}(x)\ge 1/2}. \end{equation}
	\RETURN $\widehat{\phi}_{n,s,\bw_n}^{M}(x)$. 
\end{algorithmic} 
\end{algorithm}

Next, we present the first main result of this paper: an asymptotic expansion for the regret of the M-DNN classifier with general weights ($\bw_n$).
\begin{theorem} 
\label{thm:M-DNN_re}
(Asymptotic Regret for M-DNN) Suppose the same conditions as in Proposition~\ref{thm:WNN_re}, and 
\begin{align}
\sum_{i=1}^n w_{ni}^3/(\sum_{i=1}^n w_{ni}^2)^{3/2}= o(s^{-1/2}(\log(s))^{-2}). \label{eq:extra_condition}
\end{align}
We have as $n,s \rightarrow \infty$,
	\begin{equation}
{\rm Regret}(\widehat{\phi}_{n,s,\bw_n}^{M}) = \Big[ B_1\frac{\pi}{2s}  \sum_{i=1}^n w_{ni}^2+ B_2 \Big( \sum_{i=1}^n \frac{\alpha_i w_{ni}}{n^{2/d}}\Big)^2 \Big]\{1+o(1)\}, \label{eq:M-DNN_re}
	\end{equation}
uniformly for $\bw_n\in W_{n,\beta}$.
\end{theorem}
In contrast with Proposition \ref{thm:WNN_re}, the variance term in the asymptotic regret of M-DNN in Theorem \ref{thm:M-DNN_re} is reduced by a factor of $\pi/(2s)$, while the squared bias term remains the same. This variance reduction effect is not surprising given the study of bagging \citep{buhlmann2002analyzing}, and has also been observed in the nonparametric regression setup, e.g., \cite{zhang2013divide}. Rather, the appearance of the constant $\pi/2$ is new and will motivate a new version of distributed classification in Section~\ref{sec:W-DNN}.

\begin{remark}\label{rem:thm1_1}
Theorem \ref{thm:M-DNN_re} is not a straightforward extension from Proposition \ref{thm:WNN_re} as $s$ diverges. Specifically, we need an extra normal approximation of ${\mathbb P}\big(\widehat{\phi}_{n,s,\bw_n}^{M}(x)=0\big)$ by the uniform Berry-Esseen theorem \citep{lehmann2004elements} as $s$ diverges. In fact, the factor $\pi/2$ in \eqref{eq:M-DNN_re}
comes from a simple Taylor expansion of the normal cumulative distribution function at 0; see Lemma \ref{lemma:Phi}. 
\end{remark}
\begin{remark}\label{rem:thm1_2}
Condition \eqref{eq:extra_condition} in Theorem \ref{thm:M-DNN_re} is used to bound the residual term in normal approximation by the nonuniform Berry-Esseen Theorem \citep{GP12}. Since the minimal of the left hand side is $n^{-1/2}$ (corresponding to a WNN classifier where every data point has an equal vote of $1/n$,) this condition suggests that $n^{-1/2}=o(s^{-1/2}(\log(s))^{-2})$, i.e., roughly speaking, $s/n=o(1)$, or $\gamma<1/2$. When $k$NN is trained on each subsample, the condition reduces to $k^{-1/2}=o(s^{-1/2}(\log(s))^{-2})$ which means $s$ has a smaller order than the number of effective nearest neighbors $k$ on each subsample.
\end{remark}

From \cite{S12}, we know that the minimal asymptotic regret of the oracle $K$NN (obtained based on the entire data set) is achieved when $$K=K^{*}:=\Big(\frac{dB_1}{4B_2} \Big)^{d/(d+4)}N^{4/(d+4)}.$$ Intuitively, we may want to choose $k$ in the distributed $k$NN via majority voting, i.e., M-DNN(k), as $\lceil K^*/s \rceil$. However, a direct application of Theorem~\ref{thm:M-DNN_re} reveals that an optimal choice of $k$ turns out to be 
\begin{equation}
k^*= \lceil (\pi/2)^{d/(d+4)}(K^*/s) \rceil\label{eq:opt_M_k}
\end{equation}
\noindent which minimizes the asymptotic regret. Note that (\ref{eq:opt_M_k}) only holds as $s$ diverges. The factor $\pi/2$ in \eqref{eq:M-DNN_re} has led to the additional re-scaling constant  $(\pi/2)^{d/(d+4)}$, which is always greater than one and depends on the dimension only.  

A comparative result at the constant level for general weights is given in Theorem \ref{thm:M-DNN_rr}. It presents an asymptotic regret comparison between the M-DNN and the oracle WNN, as implied by Proposition~\ref{thm:WNN_re} and Theorem~\ref{thm:M-DNN_re}. 

\begin{theorem}
\label{thm:M-DNN_rr}
(Asymptotic Regret Comparison between M-DNN and Oracle WNN) Suppose the same conditions as Theorem~\ref{thm:M-DNN_re} hold. Given an oracle WNN classifier with weights $\bw_N$, denoted as $\widehat{\phi}_{N,\bw_N}(x)$, we have as $n,s \rightarrow \infty$,
\begin{eqnarray*}
\frac{{\rm Regret}(\widehat{\phi}_{n,s,\bw_n}^{M})}{{\rm Regret}(\widehat{\phi}_{N,\bw_N})} &\longrightarrow& Q :=\left(\frac{\pi}{2}\right)^{\frac{4}{d+4}},
\end{eqnarray*}
uniformly for $\bw_n\in W_{n,\beta}$ and $\bw_N\in W_{N,\beta}$, if the weights satisfy 
\begin{eqnarray} 
\frac{1}{s} \sum_{i=1}^n w_{ni}^2/\sum_{i=1}^N w_{Ni}^2 &\longrightarrow& \left(\frac{\pi}{2}\right)^{-\frac{d}{d+4}} \;\;  {\rm and}  \label{eq:M-DNN_rr_weight1}\\
\sum_{i=1}^n \frac{\alpha_i w_{ni}}{n^{2/d}}/\sum_{i=1}^N \frac{\alpha_i w_{Ni}}{N^{2/d}} &\longrightarrow& \left(\frac{\pi}{2}\right)^{\frac{2}{d+4}}. \label{eq:M-DNN_rr_weight2}
\end{eqnarray}
\end{theorem} 

\begin{figure}[ht!]
\begin{center}
\includegraphics[scale=0.35]{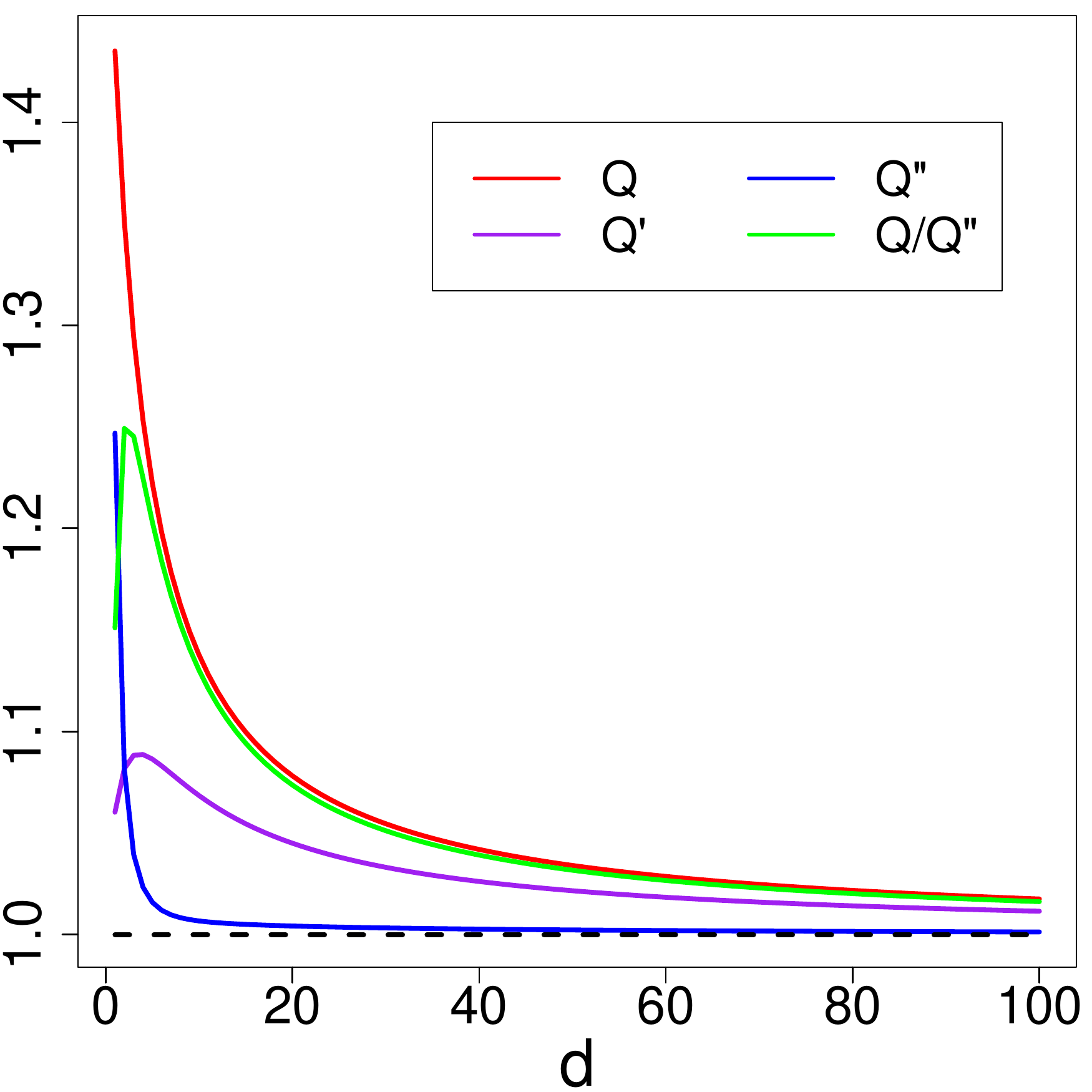}\vspace{-1em}
\caption{\label{fig:Slides_Q_prime_VS_d} $Q$, $Q'$, $Q^{''}$ and $Q/Q^{''}$for different $d$.}
\end{center}
\end{figure}

Theorem \ref{thm:M-DNN_rr} says that the M-DNN can achieve the same regret rate as the oracle version if the local weights are chosen to align with the oracle global weights, according to (\ref{eq:M-DNN_rr_weight1}) and (\ref{eq:M-DNN_rr_weight2}). The ratio of the regrets is a constant $Q$, which only depends on data dimension. In Figure \ref{fig:Slides_Q_prime_VS_d}, we can find that $Q$ is always smaller than $1.5$, and monotonically decreases to $1$ as $d$ grows. This may be viewed as a kind of ``blessing of dimensionality.'' As will be shown, $Q$ results from the majority voting step in Algorithm \ref{algo:DC_MV}, and thus we name it as the majority voting (MV) constant from now on. 

As an illustration, we show how to set local weights by applying Theorem~\ref{thm:M-DNN_rr} to the OWNN method, whose global weights are defined as
\begin{equation}\label{ownn_weights}
w_{i}^*(N, m^*)=\left\{
                \begin{array}{ll}
                \frac{1}{m^*}\Big[1+\frac{d}{2}-\frac{d\alpha_i}{2(m^*)^{2/d}} \Big], \;{\rm if}\;\;i=1,\ldots,m^*,
                  \\
                  0, \;{\rm if}\;\;i=m^*+1,\ldots,N,
                \end{array}
              \right.
\end{equation}
where
\begin{eqnarray*}
m^* &=& \lceil \Big\{\frac{d(d+4)}{2(d+2)}\Big\}^{\frac{d}{d+4}}\Big(\frac{B_1}{B_2} \Big)^{\frac{d}{d+4}}N^{\frac{4}{d+4}} \rceil. \qedhere
\end{eqnarray*} 
According to (\ref{eq:M-DNN_rr_weight1}) and (\ref{eq:M-DNN_rr_weight2}), the local weights in the optimal M-DNN should be assigned as $w_{ni}^*:=w_i^*(n, l^*)$, where 
\begin{equation}
l^*=\lceil(\pi/2)^{d/(d+4)}(m^*/s)\rceil.\label{eq:l_star}
\end{equation}
Interestingly, the above scaling factor is the same as that in \eqref{eq:opt_M_k} for the distributed $k$NN and oracle $K$NN case.


Corollary \ref{thm:opt_M-DNN} summarizes the above findings, and further discovers that $$s^*\asymp N^{2/(d+4)}$$ is a {\em sharp} upper bound for the number of partitions in order for the M-DNN method with optimal weights to achieve the same regret rate of the oracle OWNN. The ratio between the regrets is the same multiplicative constant $Q$ as stipulated in Theorem~\ref{thm:M-DNN_rr}.

\begin{corollary}
\label{thm:opt_M-DNN}
(Optimal M-DNN) Suppose the same conditions as Theorem \ref{thm:M-DNN_re} hold.

\noindent(i)  If $\gamma<2/(d+4)$, the minimum regret of M-DNN is achieved by setting $w_{ni}^*=w_i^*(n, l^*)$ with $l^*$ defined in (\ref{eq:l_star}) and $w_i^*(\cdot,\cdot)$ defined in (\ref{ownn_weights}). In addition, we have as $n,s\rightarrow\infty$,
\begin{eqnarray}
\frac{{\rm Regret}(\widehat{\phi}_{n,s,\bw_n^*}^{M})}{{\rm Regret}(\widehat{\phi}_{N,\bw_N^*})} &\longrightarrow& Q.  \label{eq:opt_M-DNN}
\end{eqnarray}


\noindent(ii) If $\gamma\ge 2/(d+4)$, we have uniformly for $\bw_n\in W_{n,\beta}$,
\begin{eqnarray*}
\liminf_{n,s\rightarrow\infty} \frac{{\rm Regret}(\widehat{\phi}_{n,s,\bw_n}^{M})}{{\rm Regret}(\widehat{\phi}_{N,\bw_N^*})} &\longrightarrow& \infty.
\end{eqnarray*}
\end{corollary} 

Furthermore, after some simple derivations, we find that the M-DNN(k) method, with an optimal choice of local $k$ neighbors, also achieves the same regret rate as the oracle OWNN, with a slightly larger ratio $QQ'$, where $$Q'=2^{-4/(d+4)}\Big((d+4)/(d+2)\Big)^{(2d+4)/(d+4)}>1.$$ Please see Figure \ref{fig:Slides_Q_prime_VS_d} for the unimodal pattern of $Q'$ versus d: $Q'$ increases to its maximal value $1.089$ as $d=4$, and then decreases to $1$ as $d$ grows. 
%

\section{DNN Classifier via Weighted Voting}\label{sec:W-DNN}
In this section, we propose another type of distributed WNN based on weighted voting, denoted as W-DNN, which helps to eliminate the multiplicative loss of regret $Q$ in M-DNN. Specifically, the local classifier $\widehat{\phi}_{n,\bw_n}^{(j)}(x)$ in Algorithm \ref{algo:DC_MV}, which outputs 0 or 1, is replaced by the regression estimator $\widehat{S}_{n,\bw_n}^{(j)}(x):=\sum^n_{i=1} w_{ni} Y_{(i)}^{(j)}$, which outputs a number $\in[0,1]$. The resulting classifier is defined as $$\widehat{\phi}_{n,s,\bw_n}^{W}(x):=\ind{\frac{1}{s}\sum_{j=1}^s\widehat{S}_{n,\bw_n}^{(j)}(x)\ge 1/2}.$$ The superscript $M$ in all notations used in Section \ref{sec:M-DNN} will be replaced by $W$ in this section.

The above simple change leads to a different asymptotic expansion of regret from M-DNN, as stated in Theorem~\ref{thm:W-DNN_re}. Specifically, the variance term is reduced by $1/s$ in contrast to $\pi/(2s)$ for the case of M-DNN. Additionally, Condition \eqref{eq:extra_condition} in Theorem \ref{thm:M-DNN_re} is not required in Theorem \ref{thm:W-DNN_re}.
\begin{theorem} 
\label{thm:W-DNN_re}
(Asymptotic Regret for W-DNN) Assuming the same conditions as in Proposition~\ref{thm:WNN_re}, we have for each $\beta\in (0,1/2)$, as $n,s \rightarrow \infty$,
		\begin{equation}
{\rm Regret}(\widehat{\phi}_{n,s,\bw_n}^{W}) = \Big [ B_1\frac{1}{s}  \sum_{i=1}^n w_{ni}^2+ B_2 \Big( \sum_{i=1}^n \frac{\alpha_i w_{ni}}{n^{2/d}}\Big)^2 \Big]\{1+o(1)\}, \label{eq:W-DNN_re}
		\end{equation}
uniformly for $\bw_n\in W_{n,\beta}$.
\end{theorem}

A consequence of this new asymptotic expansion is that the optimal local choice of $k$ in W-DNN(k) leading to the same regret as the optimal oracle $K$NN is the intuitive choice $k^{\dag}= \lceil K^*/s \rceil$, different from \eqref{eq:opt_M_k} in M-DNN(k).

Unsurprisingly, a similar result can be obtained for W-DNN with general weights. Specifically, Theorem \ref{thm:W-DNN_rr} says that W-DNN is able to achieve the same asymptotic regret as its oracle counterpart {\em without} any regret loss. 
\begin{theorem}
\label{thm:W-DNN_rr}
(Asymptotic Regret Comparison between W-DNN and Oracle WNN) Suppose the same conditions as Theorem \ref{thm:W-DNN_re} hold. We have as $n,s \rightarrow \infty$,
\begin{eqnarray*}
\frac{{\rm Regret}(\widehat{\phi}_{n,s,\bw_n}^{W}  )}{{\rm Regret}(\widehat{\phi}_{N,\bw_N} )} &\longrightarrow& 1,
\end{eqnarray*}
uniformly for $\bw_n\in W_{n,\beta}$ and $\bw_N\in W_{N,\beta}$, if the weights satisfy
\begin{eqnarray}
\frac{1}{s} \sum_{i=1}^n w_{ni}^2/\sum_{i=1}^N w_{Ni}^2 &\longrightarrow& 1\;\; {\rm and} \label{eq:W-DNN_rr_weight1}\\
\sum_{i=1}^n \frac{\alpha_i w_{ni}}{n^{2/d}}/\sum_{i=1}^N \frac{\alpha_i w_{Ni}}{N^{2/d}} &\longrightarrow& 1. \label{eq:W-DNN_rr_weight2}
\end{eqnarray}
\end{theorem} 

Theorem \ref{thm:W-DNN_rr} can be applied to OWNN to identify the local weights for W-DNN that can achieve the same minimal regret, namely $w_{ni}^\dag:=w_i^*(n, l^\dag),$ where 
\begin{equation}
    l^{\dag}= \lceil m^*/s\rceil 
\end{equation} and $w_i^*(\cdot,\cdot)$ is defined in (\ref{ownn_weights}).
Interestingly, due to fewer assumptions made for W-DNN, Corollary \ref{thm:opt_W-DNN} obtains a {\em larger} sharp upper bound $$s^{\dag}\asymp N^{4/(d+4)}$$ than that of M-DNN, suggesting that more machines can be employed in the W-DNN framework.

\begin{corollary}
\label{thm:opt_W-DNN}
(Optimal W-DNN) Suppose the same conditions as Theorem \ref{thm:W-DNN_re} hold.

\noindent(i)  If $\gamma<4/(d+4)$, the minimum regret of W-DNN is achieved by setting $w_{ni}^\dag=w_i^*(n,l^\dag)$ and $w_i^*(\cdot,\cdot)$ defined in (\ref{ownn_weights}). Additionally, we have as $n,s\rightarrow\infty$,
\begin{eqnarray}
\frac{{\rm Regret}(\widehat{\phi}_{n,s,\bw_n^{\dag}}^{W})}{{\rm Regret}(\widehat{\phi}_{N,\bw_N^*})} &\longrightarrow& 1.  \label{eq:opt_W-DNN}
\end{eqnarray}


\noindent(ii) If $\gamma \ge 4/(d+4)$, we have uniformly for $\bw_n\in W_{n,\beta}$,
\begin{eqnarray*}
\liminf_{n,s\rightarrow\infty} \frac{{\rm Regret}(\widehat{\phi}_{n,s,\bw_n}^{W})}{{\rm Regret}(\widehat{\phi}_{N,\bw_N^*})} &\longrightarrow& \infty.
\end{eqnarray*}
\end{corollary}




\section{Asymptotic Comparison with Bagging}\label{sec:cis}
In this section, we compare the DNN method with a similar classification method, i.e., bagging. The purpose of bagging is to improve unstable estimators or classifiers, especially  for high-dimensional data \citep{buhlmann2002analyzing}. In particular, we compare with the Bagged $1$-Nearest Neighbor (BNN) classifier which applies $1$-NN classifier to each bootstrapped subsample and returns the final classification by majority voting. In terms of stability, we find W-DNN is more stable than BNN, while M-DNN is less stable.

To facilitate this comparison, we first introduce the notion of classification instability (CIS) introduced in \cite{SQC16}. For a classification procedure, it is desired that, with high probability, classifiers trained from different samples yield the same prediction for the same object.  Intuitively, CIS is an average probability that the same object is classified to two different classes in two separate runs of a learning algorithm on the data set with the same underlying distribution. 
\begin{defi} 
\label{def:CIS}
(CIS, \cite{SQC16}) Define the classification instability of a classification procedure $\Psi$ as
\begin{equation*}
{\rm CIS}(\Psi) = {\mathbb E}_{{\cal D}_1,{\cal D}_2} \Big[ {\mathbb P}_X\Big(\widehat{\phi}_{n1} (X) \neq \widehat{\phi}_{n2} (X)\Big) \Big], 
\end{equation*}
where $\widehat\phi_{n1}=\Psi(\mathcal D_1)$ and $\widehat\phi_{n2}=\Psi(\mathcal D_2)$ are the classifiers obtained by applying the classification procedure $\Psi$ to samples ${\cal D}_1$ and ${\cal D}_2$ that are two i.i.d. copies of the ${\cal D}$.
\end{defi}
 
Theorem \ref{thm:DNN_ci} provides the asymptotic CIS for M-DNN and W-DNN.
\begin{theorem} 
\label{thm:DNN_ci}
(Asymptotic CIS for M-DNN and W-DNN) Suppose the same conditions as Proposition \ref{thm:WNN_re} hold; additionally, Condition \eqref{eq:extra_condition} is assumed for the case of M-DNN. We have as $n,s \rightarrow \infty$, 
	\begin{eqnarray}
{\rm CIS}(\widehat{\phi}_{n,s,\bw_n}^{M}) &=& B_3 \sqrt{\frac{\pi}{2s}} \Big(\sum_{i=1}^n w_{ni}^2\Big)^{1/2} \{1+o(1)\}, \label{eq:M-DNN_ci}\\
{\rm CIS}(\widehat{\phi}_{n,s,\bw_n}^{W}) &=& B_3\frac{1}{\sqrt{s}}  \Big(\sum_{i=1}^n w_{ni}^2\Big)^{1/2} \{1+o(1)\}, \label{eq:W-DNN_ci}
	\end{eqnarray}
uniformly for $\bw_n\in W_{n,\beta}$, where the constant $B_3 = 4B_1/\sqrt{\pi}>0$.
\end{theorem}

\citet{hall2005properties} showed that, for large $N$, the ``infinite simulation" version of BNN classifier (with or without replacement) is approximately equivalent to a WNN classifier with the weight $$w_{Ni} = q(1-q)^{i-1}/[1-(1-q)^N]$$ for $i=1,\ldots,N$, where $q$ is the resampling ratio $m/N$. Hence, the CIS of BNN can be derived from the general CIS formula given in \cite{SQC16}. We denote BNN as $\widehat{\phi}_{N,q}$ and the optimal BNN as $\widehat{\phi}_{N,q*}$ where the optimal $q*$ is defined in (3.5) of \citet{S12}. 

\begin{corollary}
\label{thm:DNN_BNN}
(Asymptotic CIS Comparison for M-DNN, W-DNN and BNN) Suppose the same conditions as Theorem \ref{thm:DNN_ci} hold. If the weights of M-DNN satisfy \eqref{eq:M-DNN_rr_weight1} and \eqref{eq:M-DNN_rr_weight2}, and the weights of W-DNN satisfy \eqref{eq:W-DNN_rr_weight1} and \eqref{eq:W-DNN_rr_weight2}, we have as $n,s \rightarrow \infty$,
\begin{eqnarray*}
\frac{{\rm Regret}(\widehat{\phi}_{n,s,\bw_n^*}^{M})}{{\rm Regret}(\widehat{\phi}_{N,q*})} \rightarrow \frac{Q}{Q^{''}}>1  \;\;\;&{\rm and}&\;\;\;
\frac{{\rm Regret}(\widehat{\phi}_{n,s,\bw_n^{\dag}}^{W})}{{\rm Regret}(\widehat{\phi}_{N,q*})} \rightarrow \frac{1}{Q^{''}}<1, \\
\frac{{\rm CIS}(\widehat{\phi}_{n,s,\bw_n^*}^{M})}{{\rm CIS}(\widehat{\phi}_{N,q*})} \rightarrow \sqrt{\frac{Q}{Q^{''}}}>1  \;\;\;&{\rm and}&\;\;\;
\frac{{\rm CIS}(\widehat{\phi}_{n,s,\bw_n^{\dag}}^{W})}{{\rm CIS}(\widehat{\phi}_{N,q*})} \rightarrow \sqrt{\frac{1}{Q^{''}}}<1,
\end{eqnarray*}
where $Q^{''}= 2^{-8/(d+4)}\Gamma(2+2/d)^{2d/(d+4)}\Big(\frac{d+4}{d+2}\Big)^{(2d+4)/(d+4)}$.
\end{corollary}

Corollary \ref{thm:DNN_BNN} implies that both the optimal M-DNN and W-DNN have the same regret and CIS rates as the optimal BNN, and their differences are in terms of multiplicative constants $Q/Q^{''}$ and $1/Q^{''}$. In Figure \ref{fig:Slides_Q_prime_VS_d}, both $Q/Q^{''}$ and $Q''$ are larger than $1$. Therefore, in terms of stability, M-DNN is less stable than BNN while W-DNN is more stable. In addition, it is interesting to note that the CIS ratios are square roots of the regret ratios.


\section{Numerical Studies}\label{sec:exp}
In this section, we illustrate the effectiveness of the DNN methods using simulations and real examples. 
\subsection{Simulations}\label{sec:simulation}
In the simulated studies, we compare DNN methods with the oracle $K$NN and the oracle OWNN methods respectively, with slightly different emphases. In comparing DNN(k) ($k$NN is trained at each subsample) with the oracle $K$NN, we aim to verify the main results in Theorem \ref{thm:M-DNN_rr} and Theorem \ref{thm:W-DNN_rr}, namely, the M-DNN and W-DNN can approximate or attain the same performance as the oracle method. In comparing the DNN methods with optimal local weights and the oracle OWNN method, we aim to verify the sharpness of upper bound on $\gamma$ in Corollary \ref{thm:opt_M-DNN} and \ref{thm:opt_W-DNN}. This is by showing that the difference in performance between the DNN methods and the oracle OWNN deviates when $\gamma$ is greater than the theoretical upper bound.

Three general simulation settings are considered. Simulation 1 allows a relatively easy classification task, Simulation 2 examines the bimodal effect, and Simulation 3 combines bimodality with dependence between variables. 

In Simulation 1, $N =27000$ and $d=4,6,8$. The two classes are generated as $P_1\sim N(0_d,\mathbb{I}_d)$ and $P_0\sim N(\frac{2}{\sqrt{d}}1_d,\mathbb{I}_d)$ with the prior class probability $\pi_1={\mathbb P}(Y=1)=1/3$. Simulation 2 has the same setting as Simulation 1, except that both classes are bimodal with $P_1\sim 0.5N(0_d,\mathbb{I}_d)+0.5 N(3_d,2\mathbb{I}_d)$ and $P_0\sim 0.5N(1.5_d,\mathbb{I}_d)+0.5 N(4.5_d,2\mathbb{I}_d)$. Simulation 3 has the same setting as Simulation 2, except that $P_1\sim 0.5N(0_d,\Sigma)+0.5 N(3_d,2\Sigma)$ and $P_0\sim 0.5N(1.5_d,\Sigma)+0.5 N(4.5_d,2\Sigma)$ with $\pi_1=1/2$, where $\Sigma$ is the Toeplitz matrix whose $j$-th entry of the first row is $0.6^{j-1}$.

Let $s=N^{\gamma}$, and we choose the exponent $\gamma=0.0, 0.1 \dots 0.8$. When comparing the $k$NN methods, the number of neighbors $K$ in the oracle $K$NN is chosen as $K=N^{0.7}$. The number of local neighbors in M-DNN(k) and W-DNN(k) are chosen as $k=\lceil (\pi/2)^{d/(d+4)}K/s \rceil$ and $k=\lceil K/s \rceil$ as suggested by Theorem \ref{thm:M-DNN_rr} and Theorem \ref{thm:W-DNN_rr} respectively. These $k$ values are truncated at 1, since we cannot have a fraction of an observation. In comparing with the oracle OWNN method, the $m^{*}$ parameter in OWNN is tuned using cross-validation. The parameter $l$ in M-DNN and W-DNN for each subsample are chosen as $l^*=\lceil(\pi/2)^{d/(d+4)}(m^*/s)\rceil$ and $l^{\dag}= \lceil m^*/s\rceil$ as stated in Corollary \ref{thm:opt_M-DNN} and Corollary \ref{thm:opt_W-DNN} respectively. For both comparisons, the test set is independently generated with $1000$ observations. We repeat the simulation for $1000$ times for each $\gamma$ and $d$. Here the empirical risk (test error) and the computation time are calculated for each of the methods.

\begin{figure}[!ht]
	\centering\vspace{-0.5em}
\includegraphics[width=0.32\textwidth,height=0.28\textwidth]{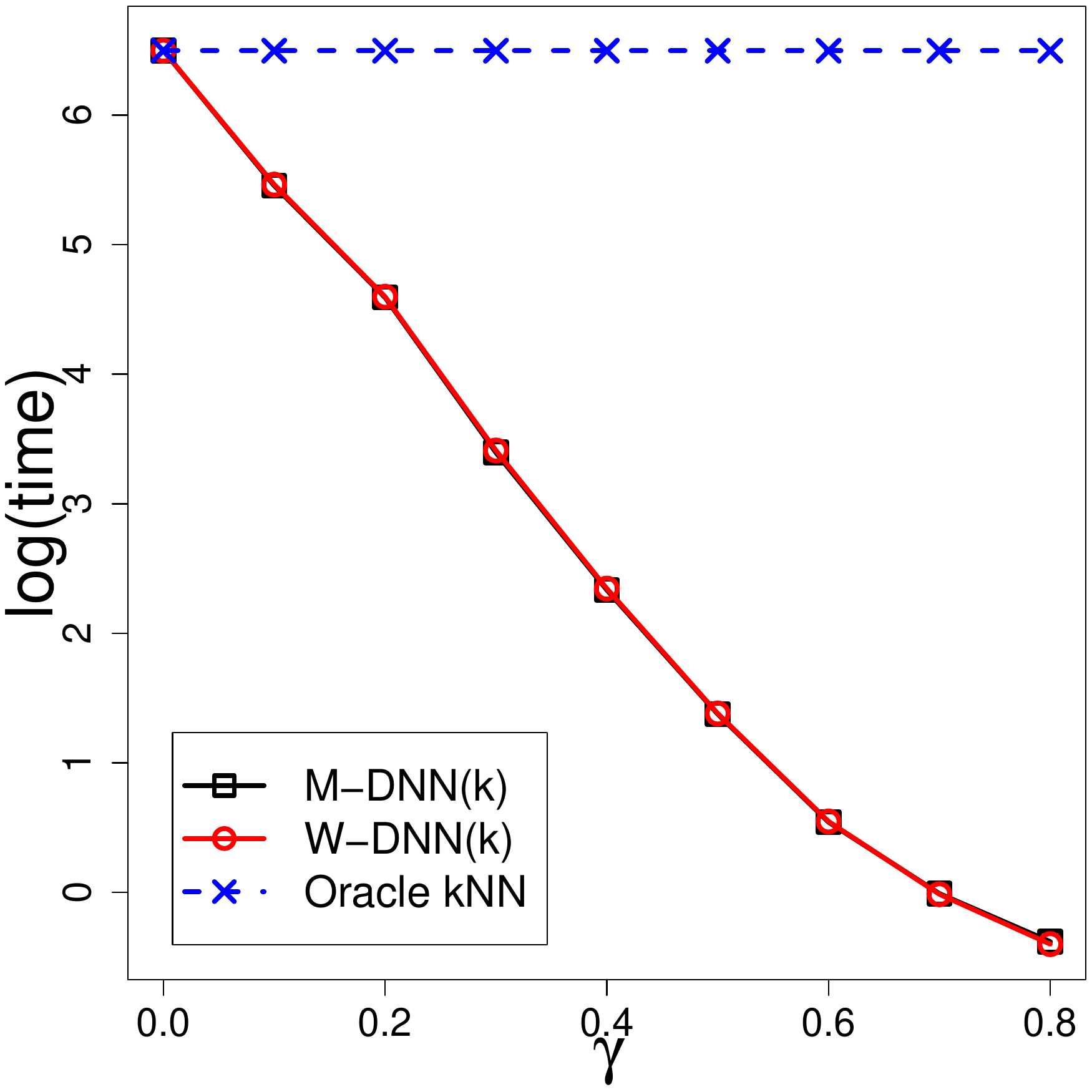}
\includegraphics[width=0.32\textwidth,height=0.28\textwidth]{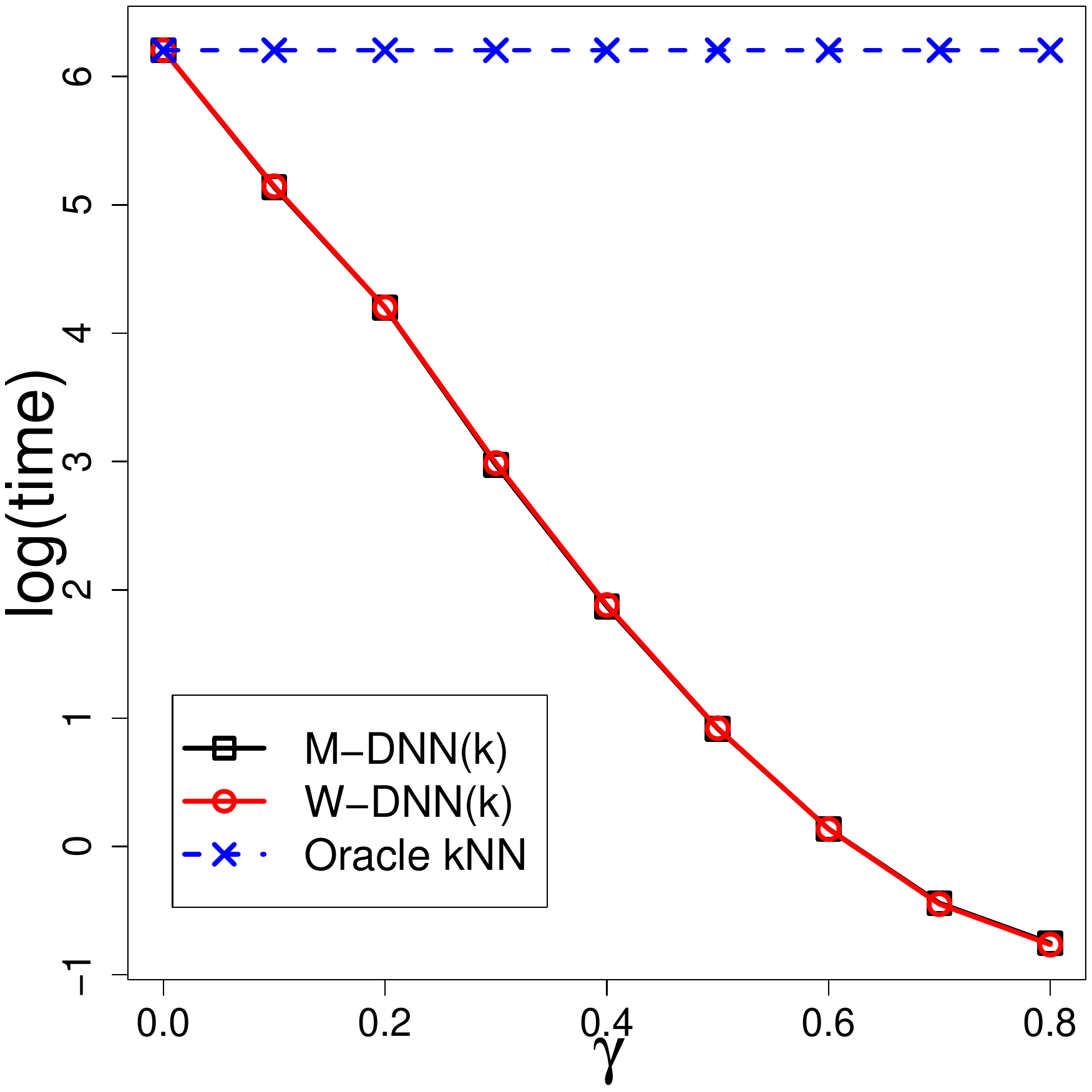}
\includegraphics[width=0.32\textwidth,height=0.28\textwidth]{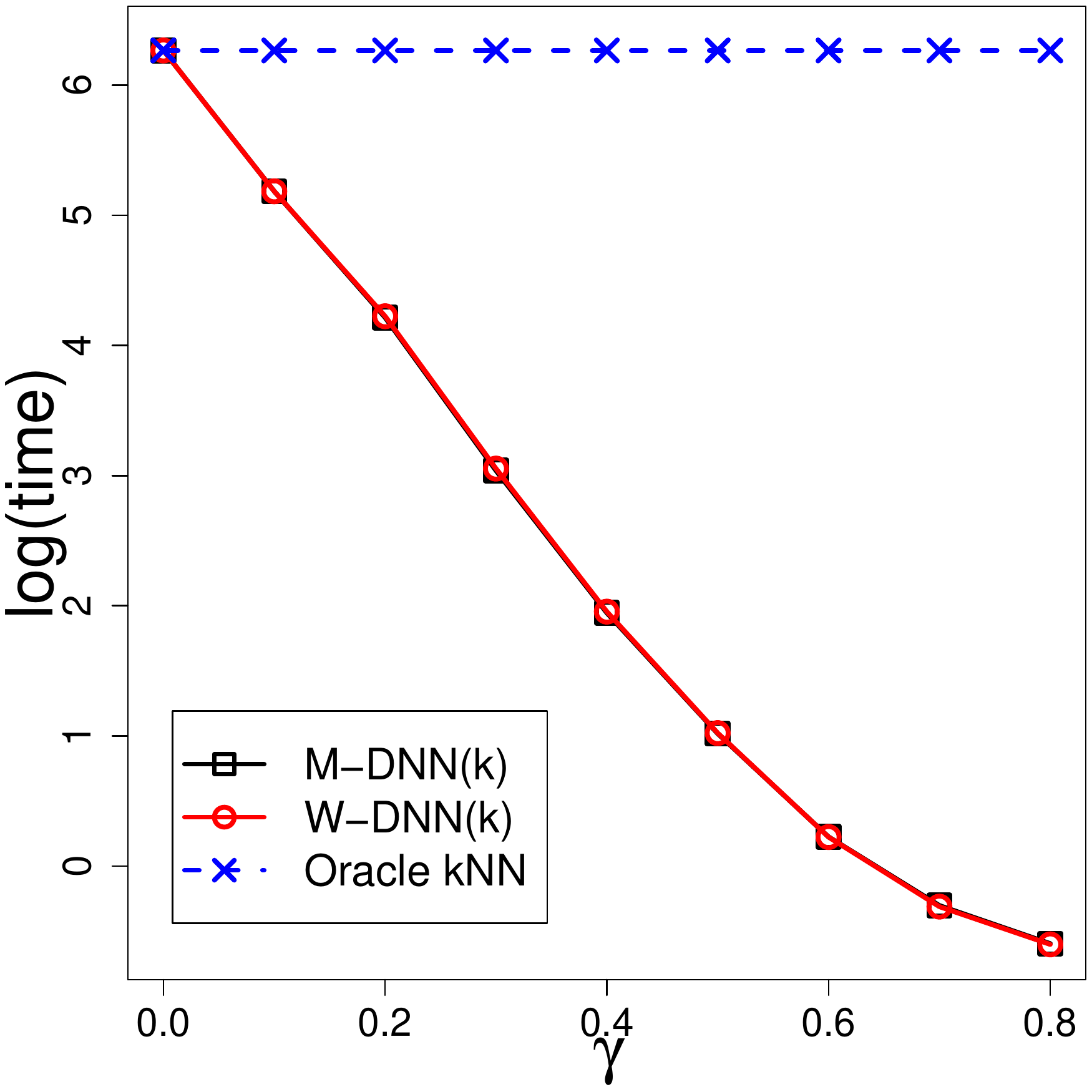}
\vspace{-1em}	
\caption{Computation time (seconds) of M-DNN(k), W-DNN(k), and oracle $K$NN for different $\gamma$. Left/middle/right: Simulation $1/2/3$, $d=4/6/8$.} \label{fig:sim_time_gamma}
\vspace{-0.5em}
\end{figure}

\begin{figure}[!ht]
	\centering\vspace{-1em}
	\includegraphics[width=0.32\textwidth,height=0.28\textwidth]{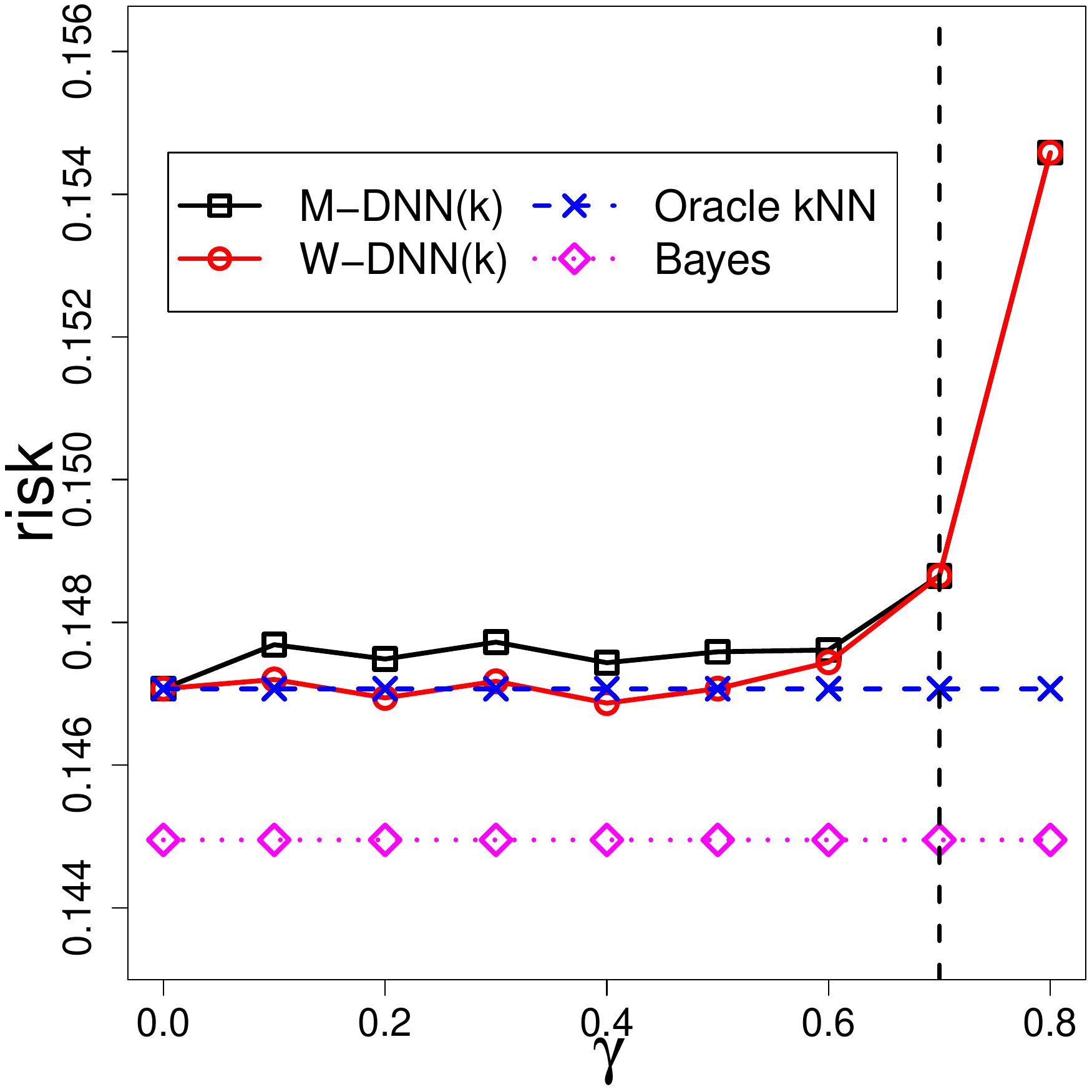}
    \includegraphics[width=0.32\textwidth,height=0.28\textwidth]{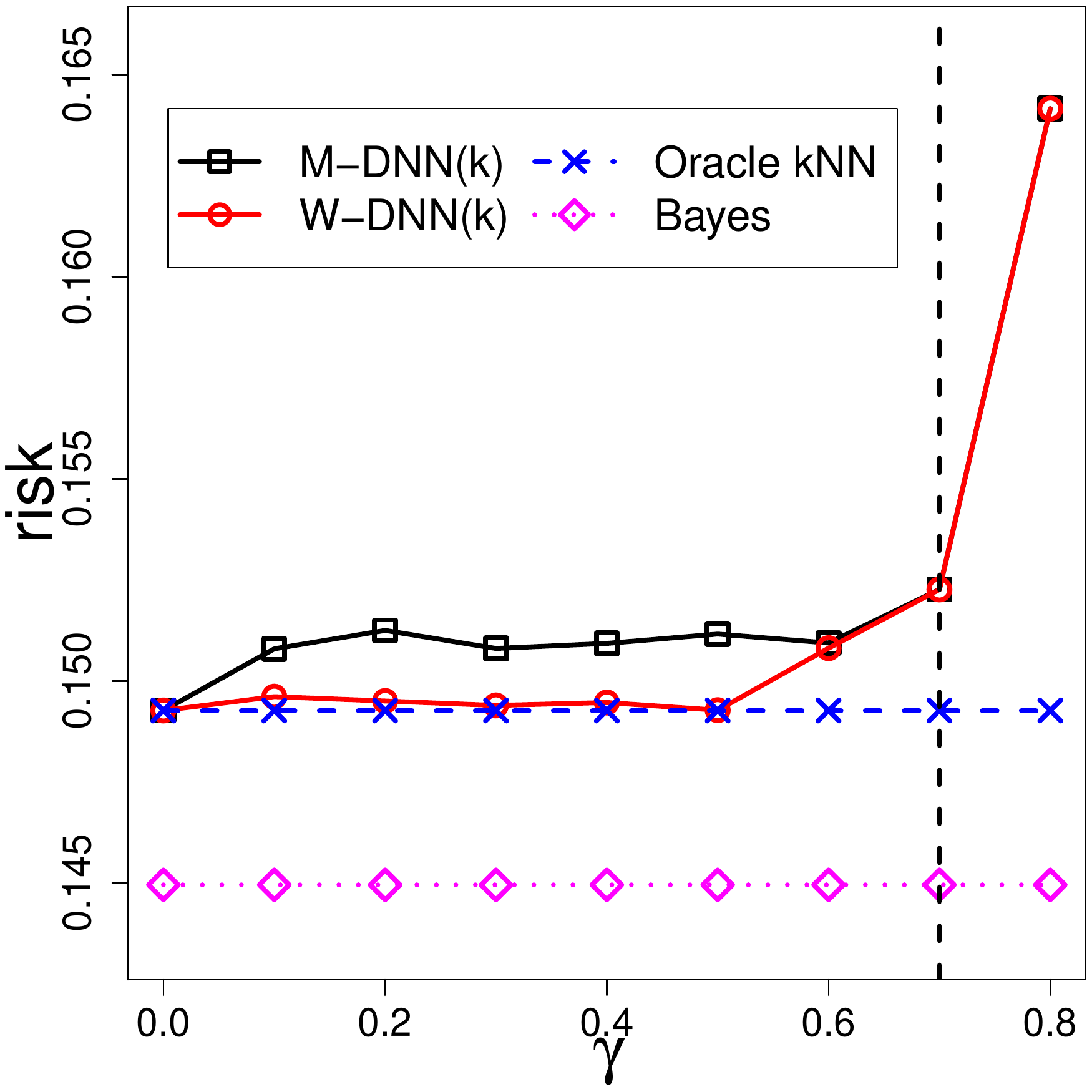}
    \includegraphics[width=0.32\textwidth,height=0.28\textwidth]{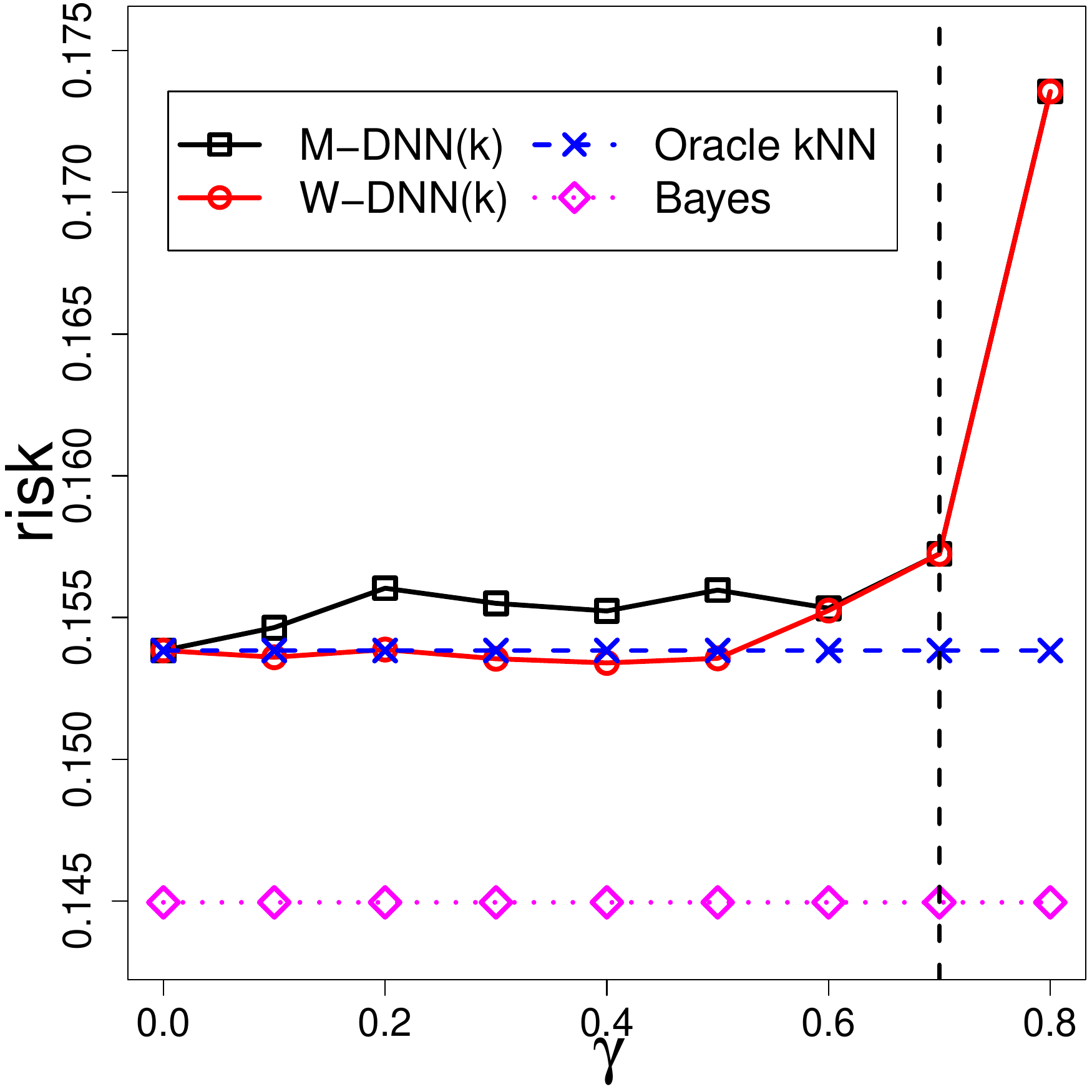}
    \includegraphics[width=0.32\textwidth,height=0.28\textwidth]{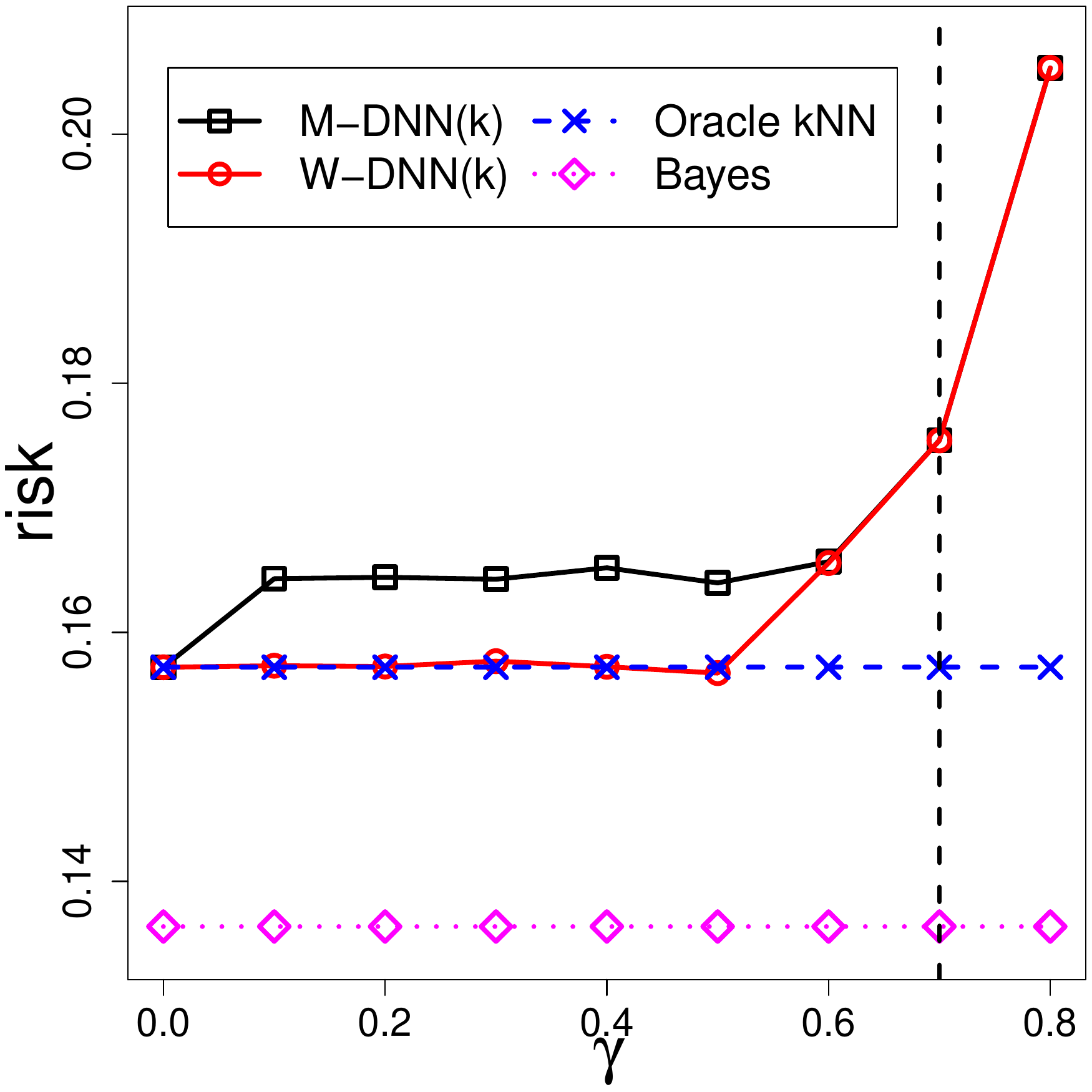}
    \includegraphics[width=0.32\textwidth,height=0.28\textwidth]{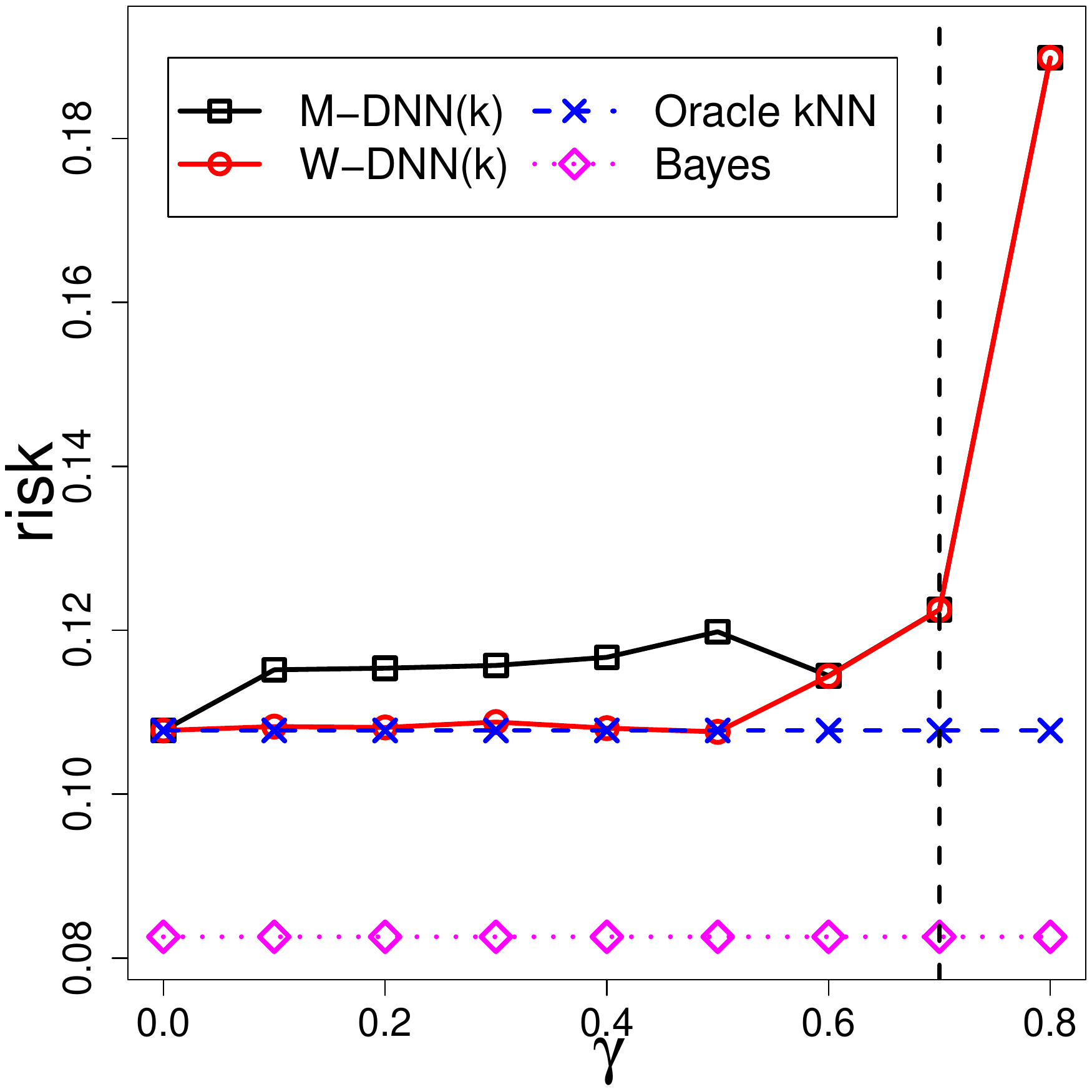}
    \includegraphics[width=0.32\textwidth,height=0.28\textwidth]{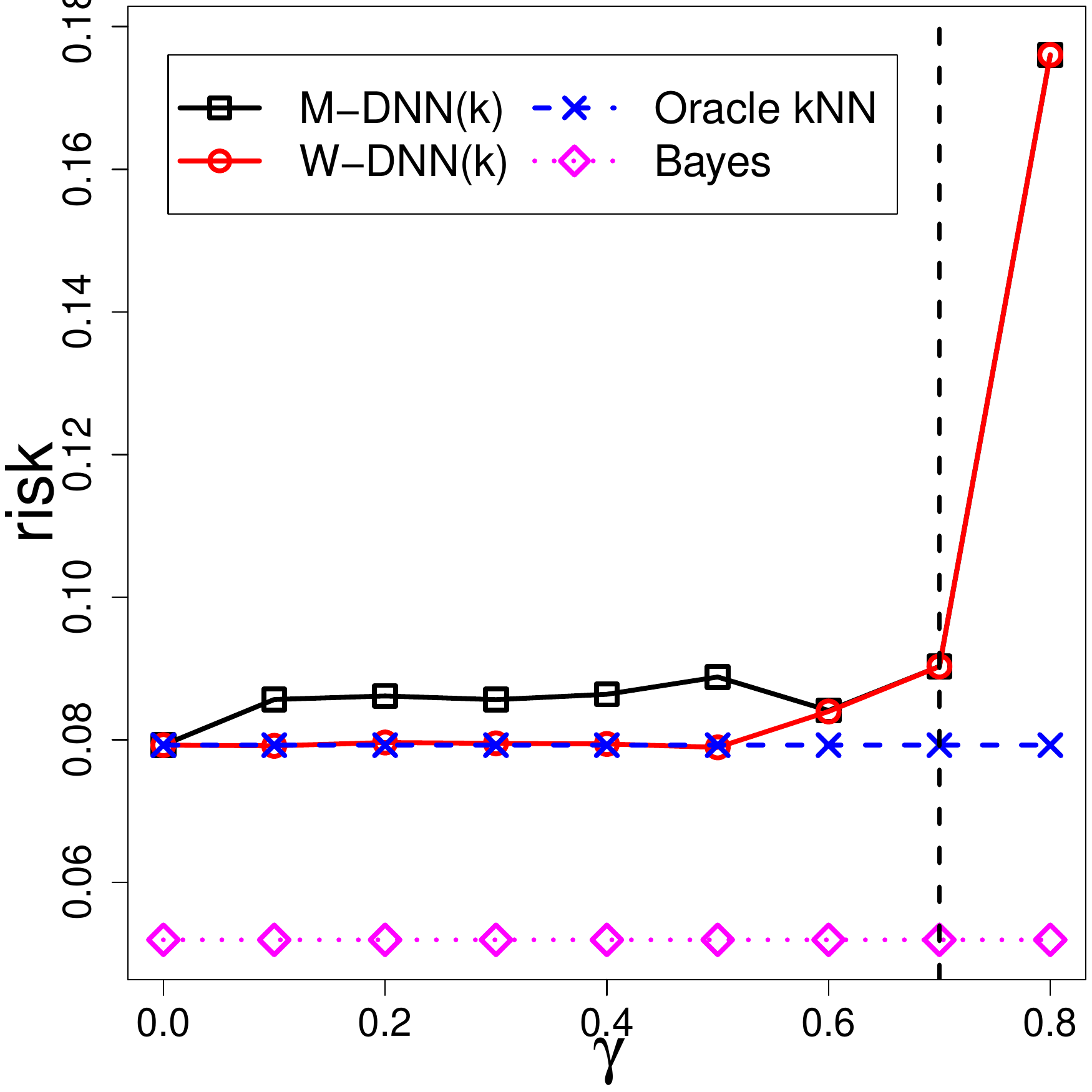}
    \includegraphics[width=0.32\textwidth,height=0.28\textwidth]{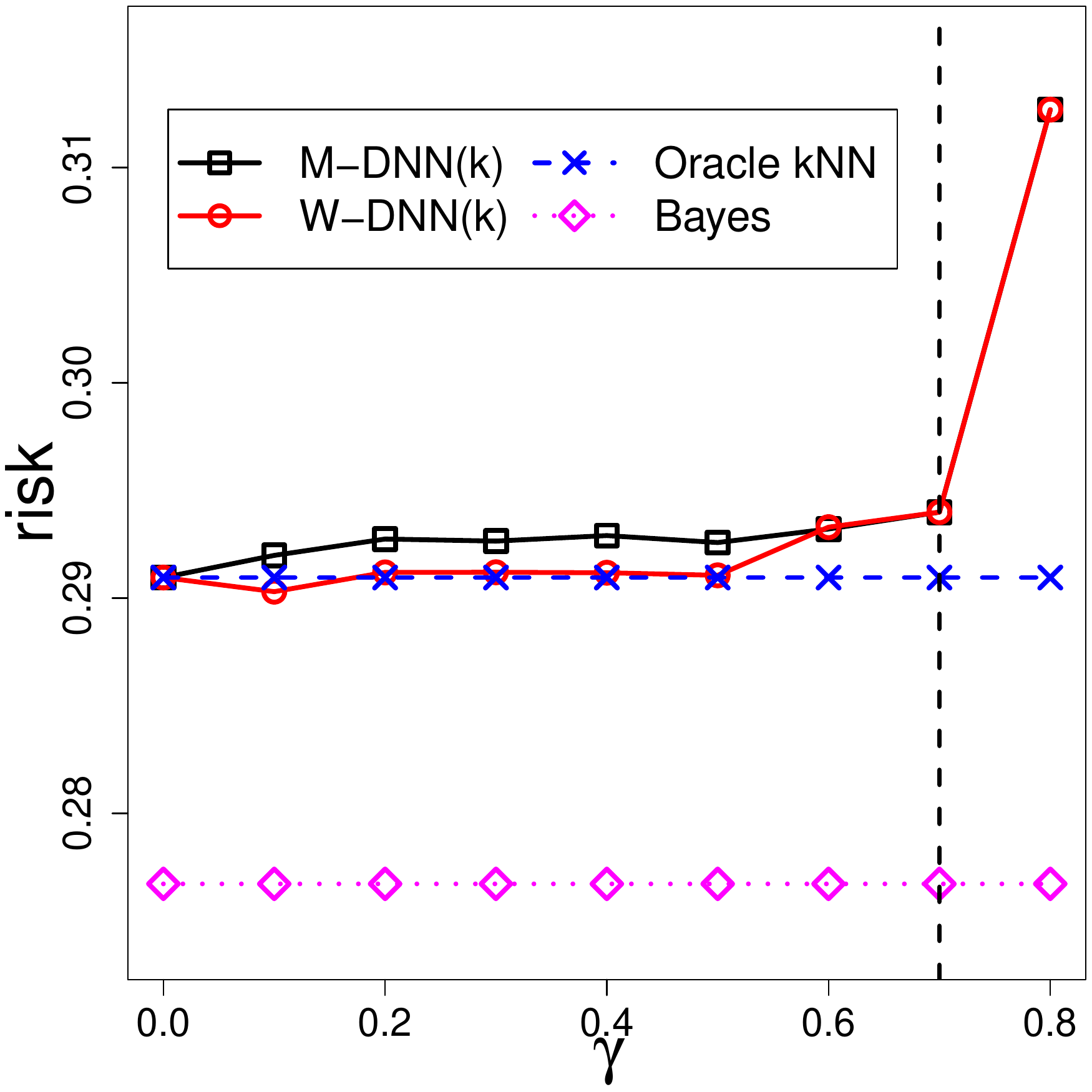}
    \includegraphics[width=0.32\textwidth,height=0.28\textwidth]{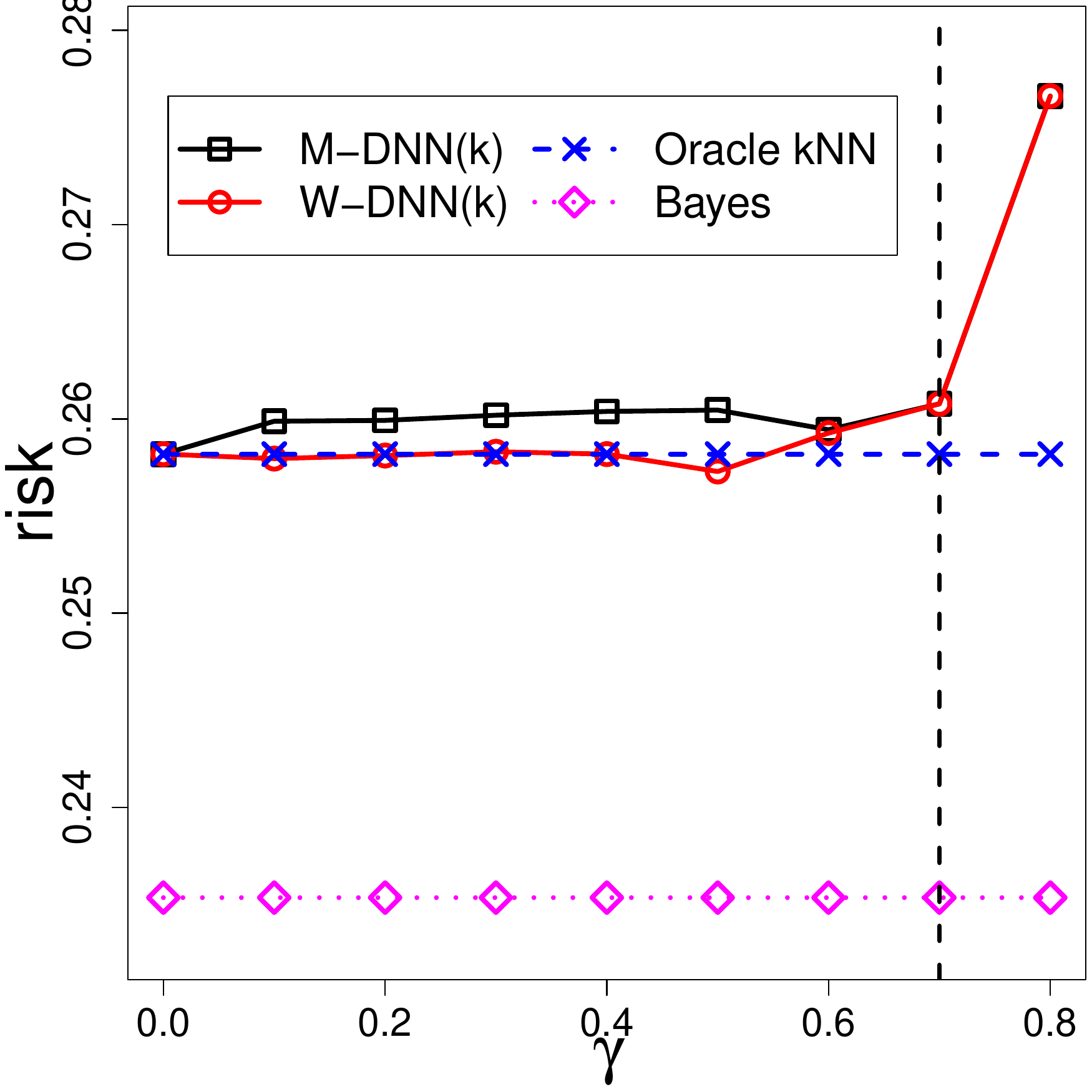}
    \includegraphics[width=0.32\textwidth,height=0.28\textwidth]{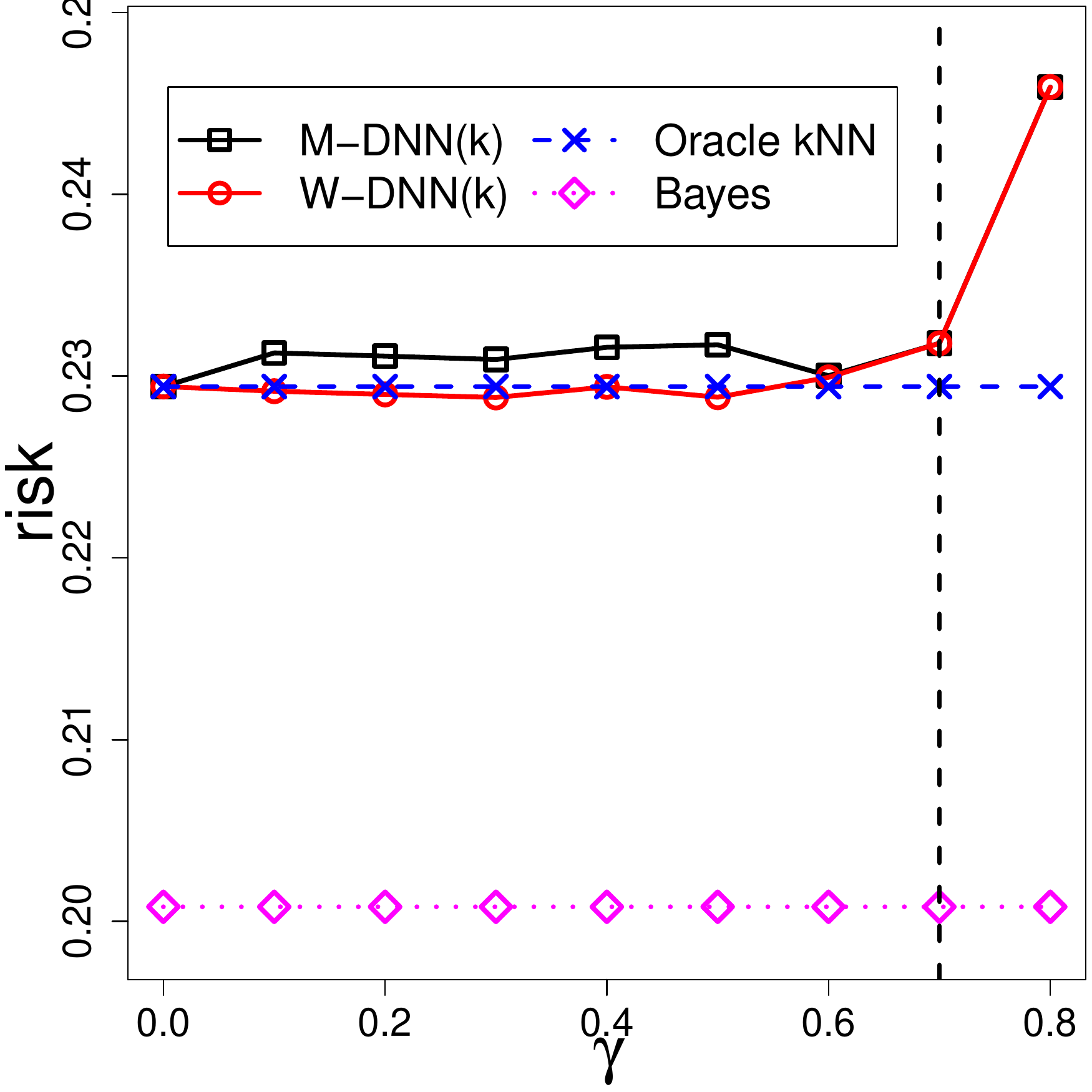}
\vspace{-1em}	
\caption{Risk of M-DNN(k), W-DNN(k), oracle $K$NN and the Bayes rule for different $\gamma$. Top/middle/bottom: Simulation $1/2/3$; left/middle/right: $d=4/6/8$. $\gamma=0.7$ is shown as a vertical line: DNN methods at or after this line have only 1 nearest neighbor at each subsample which participates in the prediction.} \label{fig:sim_risk_gamma}
\vspace{-0.9em}
\end{figure}

Figure \ref{fig:sim_time_gamma} shows that M-DNN(k) and W-DNN(k) require similar computing time, and both are significantly faster than the oracle method. As the number of subsamples increases, the running time decreases, which shows the time benefit of the distributed learning framework. The computing time comparison with the oracle OWNN is omitted since the message is the same.

The comparison between the risks of the three $k$NN methods (one oracle and two distributed) are reported in Figure \ref{fig:sim_risk_gamma}. For smaller $\gamma$ values, the risk curve for W-DNN(k) overlaps with that of the oracle kNN, while the curve for M-DNN(k) has a conceivable gap with both. These verify the main results in Theorem \ref{thm:M-DNN_rr} and Theorem \ref{thm:W-DNN_rr}.
The performance of the M-DNN(k) method starts to deviate from the oracle kNN since $\gamma=0.6$. As $s$ increases and goes beyond the threshold $\gamma=0.7$, the risk deteriorates more quickly. These may be caused by the finite (or very small) number of voting neighbors $k$ at each subsample, which means the requirements $k=\lceil (\pi/2)^{d/(d+4)}K/s \rceil \rightarrow \infty$, $k=\lceil K/s \rceil \rightarrow \infty$ suggested by Theorem \ref{thm:M-DNN_rr} and Theorem \ref{thm:W-DNN_rr} respectively are not satisfied. Specifically, when $\gamma=0.6, 0.7,0.8$, the number of voting neighbors $k$ are no more than 3 in these simulated examples. We did not tune the parameters and simply set K in the oracle $K$NN as $N^{0.7}$, since the results in Theorem \ref{thm:M-DNN_rr} and Theorem \ref{thm:W-DNN_rr} should hold for any reasonable weights (or reasonable choice of $k$), not necessarily the optimal one.

\begin{figure}[!ht]
	\centering\vspace{-0.5em}
\includegraphics[width=0.32\textwidth,height=0.28\textwidth]{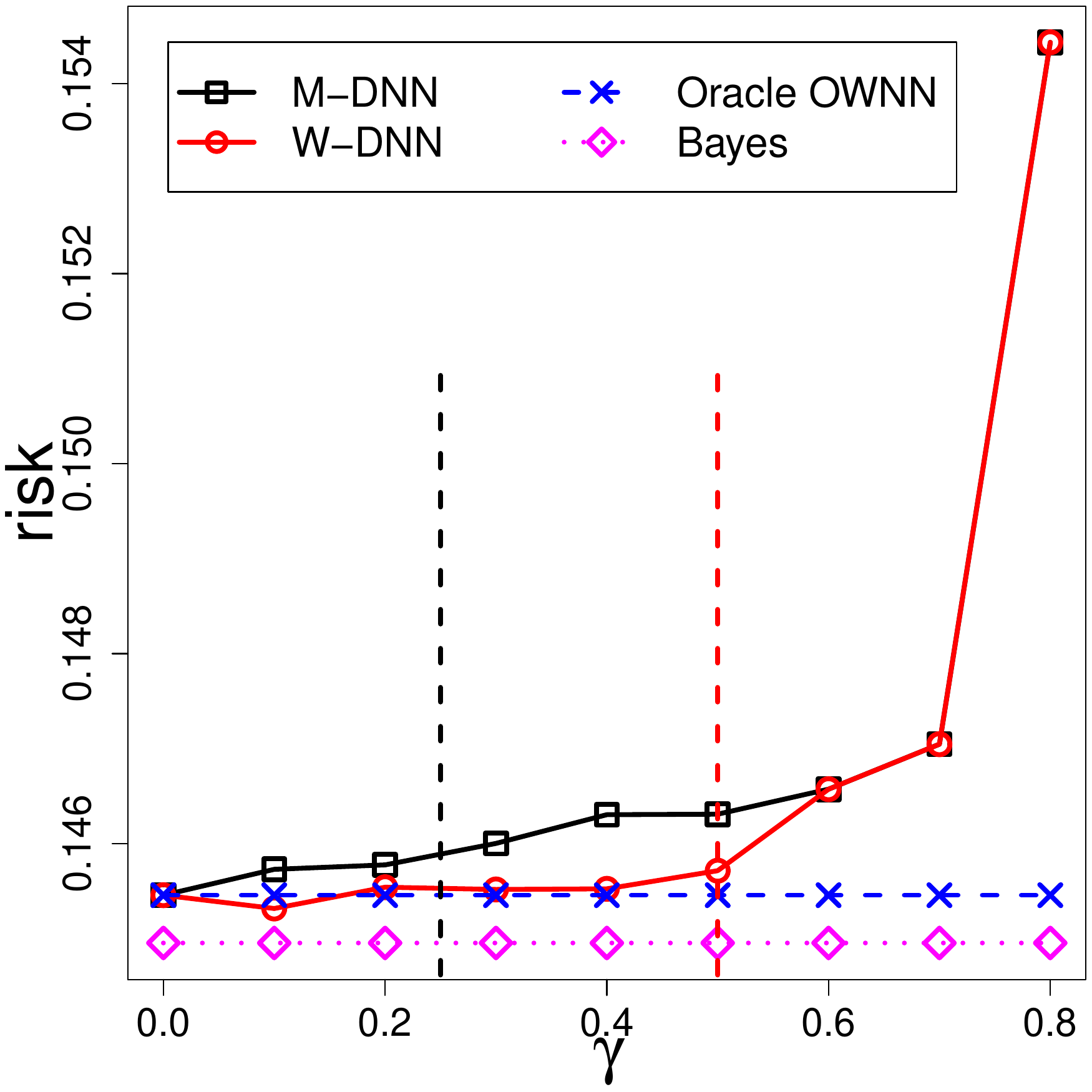}
\includegraphics[width=0.32\textwidth,height=0.28\textwidth]{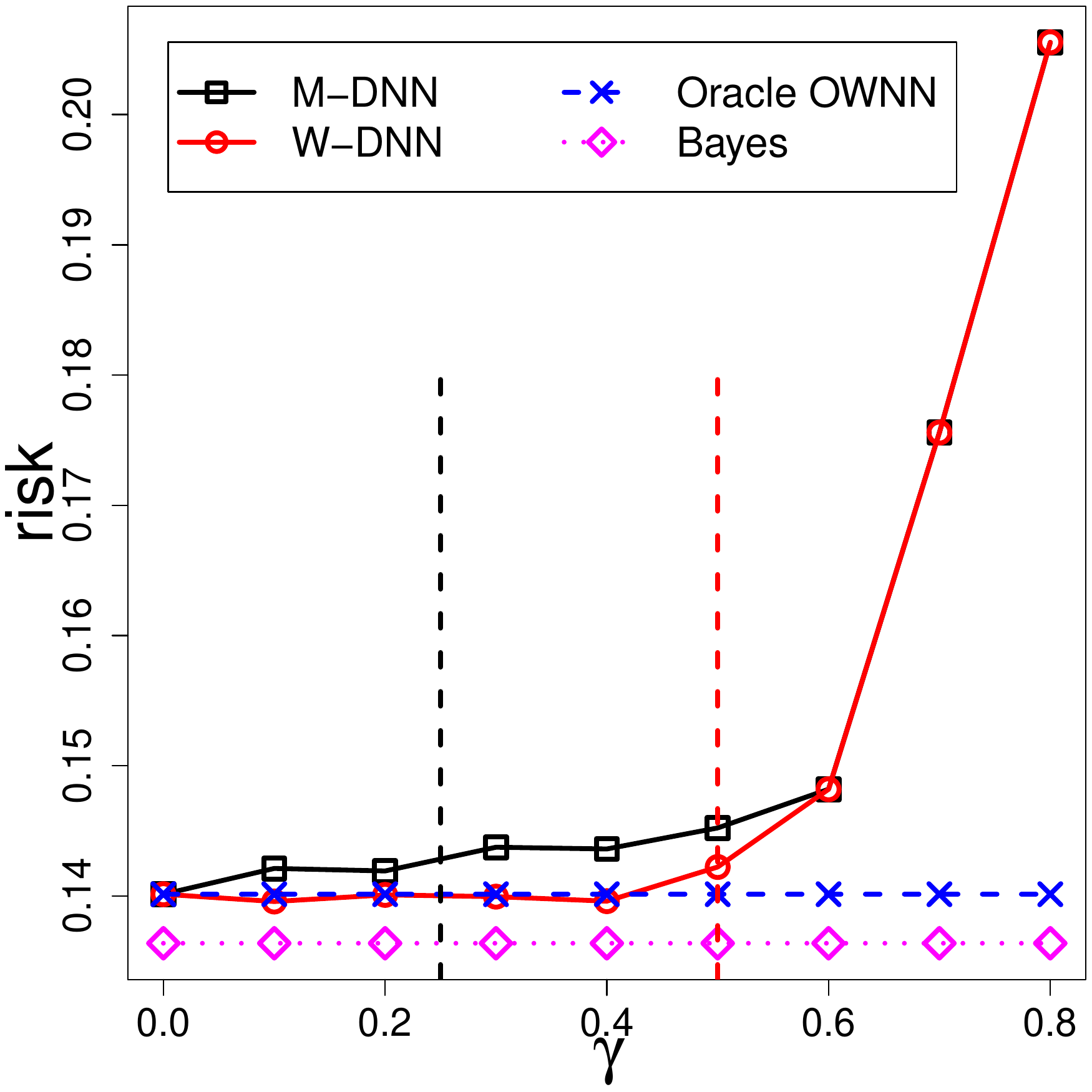}
\includegraphics[width=0.32\textwidth,height=0.28\textwidth]{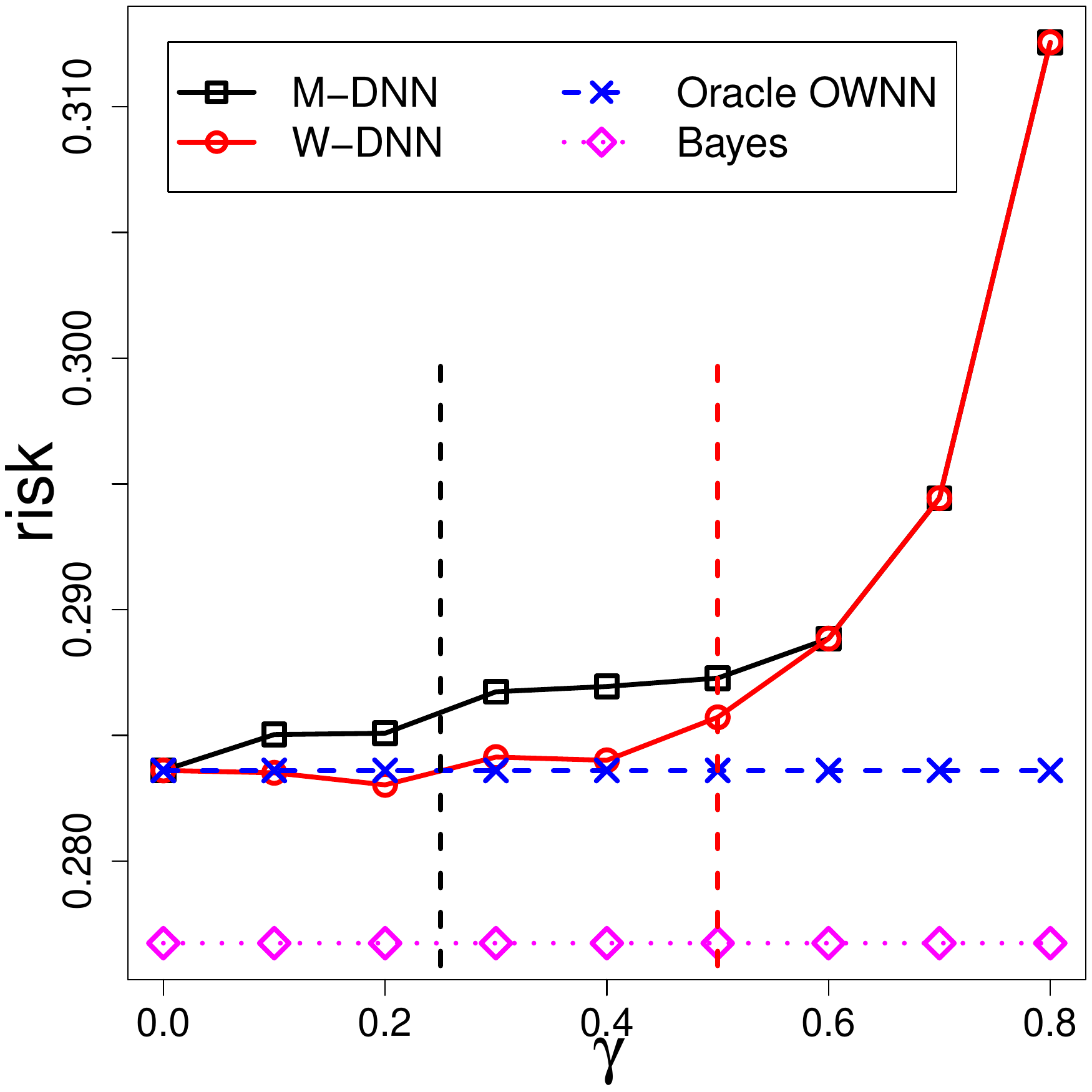}
\vspace{-1em}	
\caption{Risk of optimal M-DNN, W-DNN, oracle OWNN and the Bayes rule for different $\gamma$. Left/middle/right: Simulation $1/2/3$, $d=4$. Upper bounds for number of subsamples in optimal M-DNN ($\gamma=1/4$) and W-DNN ($\gamma=1/2$) are shown as two vertical lines.} \label{fig:sim_opt_risk_gamma}
\vspace{-0.5em}
\end{figure}

On the other hand, since the comparison with the oracle OWNN is meant to verify the sharp upper bound for $\gamma$ in the optimal weight setting (Corollary \ref{thm:opt_M-DNN} and Corollary \ref{thm:opt_W-DNN}), we carefully tune the weights in the oracle OWNN method in order to reach the optimality. Figure \ref{fig:sim_opt_risk_gamma} shows the comparison of risks for M-DNN, W-DNN and oracle OWNN methods. Our focus here is when the two DNN methods start to have significantly worse performance than the oracle OWNN, and the answers lie in the upper bounds in Corollary \ref{thm:opt_M-DNN} and Corollary \ref{thm:opt_W-DNN}. For simplicity, we set $d=4$, which leads to an upper bound of $2/(d+4)=0.25$ for the M-DNN method, and an upper bound of $4/(d+4)=0.5$ for the W-DNN method. These upper bounds are shown as vertical lines in Figure \ref{fig:sim_opt_risk_gamma}. Specifically, the M-DNN deteriorates much earlier than W-DNN with much few machines (subsamples) at its disposal. The W-DNN method performs much better than M-DNN, having almost the same performance as the OWNN method for small $\gamma$. Yet, even W-DNN has a limit at $4/(d+4)$ when comparing to OWNN. These verify the results in Corollary \ref{thm:M-DNN_rr} and Corollary \ref{thm:W-DNN_rr}.

\begin{table}[ht!]
	\caption{Risk (in $\%$) of M-DNN(k) and W-DNN(k) compared to oracle $K$NN and OWNN for $7$ real datasets. The speedup factor is defined as computing time of the oracle $K$NN divided by the time of the slower of the two DNN(k) methods.}  \label{tab:real1} 	
	\begin{center}
	\setlength\tabcolsep{3.5pt} 
	\small
	\begin{tabular}{ l r r | r r r | r r | r}	
		\hline
		Data & N & d & $\gamma$ & M-DNN($k$) & W-DNN($k$) & $k$NN & OWNN & Speedup\\
		\hline
		      &       &       & $0.1$ & $15.36$& $15.22$ &         &         & $1.19$ \\
		Musk1 & $476$ & $166$ & $0.2$ & $15.45$& $15.28$ & $15.10$ & $14.98$ & $2.23$ \\
		      &       &       & $0.3$ & $15.82$& $15.53$ &         &         & $3.30$ \\
		\hline
		        &        &        & $0.1$ & $4.01$ & $3.70$ &        &        & $2.54$ \\
		Gisette & $6000$ & $5000$ & $0.2$ & $4.18$ & $3.94$ & $3.62$ & $3.48$ & $4.55$ \\
		        &        &        & $0.3$ & $4.10$ & $3.88$ &        &        & $10.68$ \\
		\hline
		      &        &       & $0.1$ & $3.91$ & $3.78$ &        &        & $3.30$ \\
		Musk2 & $6598$ & $166$ & $0.2$ & $3.91$ & $3.75$ & $3.54$ & $3.28$ & $5.69$ \\
		      &        &       & $0.3$ & $4.23$ & $3.98$ &        &        & $15.62$ \\
		\hline
		      &         &     & $0.1$ & $2.26$ & $2.20$ &        &        & $3.27$  \\
		HTRU2 & $17898$ & $8$ & $0.2$ & $2.23$ & $2.18$ & $2.19$ & $2.12$ & $7.96$  \\
		      &         &     & $0.3$ & $2.30$ & $2.22$ &        &        & $21.90$  \\
		\hline
		      &         &     & $0.1$ & $0.69$ & $0.65$ &        &        & $3.01$ \\
		Occup & $20560$ & $6$ & $0.2$ & $0.75$ & $0.73$ & $0.65$ & $0.60$ & $7.47$ \\	
		      &         &     & $0.3$ & $0.86$ & $0.80$ &        &        & $21.25$ \\	
		\hline
			   &         &      & $0.1$  & $19.37$ & $19.28$ &         &         & $3.00$ \\
		Credit & $30000$ & $24$ & $0.2$  & $19.31$ & $19.23$ & $19.08$ & $18.96$ & $7.67$ \\
		       &         &      & $0.3$  & $19.33$ & $19.27$ &         &         & $23.57$ \\
		\hline
		     &         &      & $0.1$ & $23.58$ & $22.32$ &         &         & $4.02$ \\
		SUSY & $5000$K & $18$ & $0.2$ & $23.63$ & $22.30$ & $21.57$ & $21.11$ & $16.56$ \\
		     &         &      & $0.3$ & $23.76$ & $22.51$ &         &         & $72.78$ \\		
		\hline
	\end{tabular}
	\end{center}
\end{table}

\begin{table}[ht!]
	\caption{Risk and CIS (in $\%$) of M-DNN with optimal weights $\bw^*$ and W-DNN with optimal weights $\bw^\dag$ compared to oracle BNN for $7$ real datasets.}  \label{tab:real2} 	
	\begin{center}
	\setlength\tabcolsep{3.5pt} 
	\small
	\begin{tabular}{ l | c | r r r | r r r}
		\hline
		&&\multicolumn{3}{c|}{Risk}&\multicolumn{3}{c}{CIS}\\
		Data & $\gamma$ & M-DNN$_{\bw^*}$ & W-DNN$_{\bw^{\dag}}$ & BNN & M-DNN$_{\bw^*}$ & W-DNN$_{\bw^{\dag}}$ & BNN \\
		\hline
		      & $0.1$ & $15.11$ & $14.98$ &         & $23.63$ & $23.27$ &   \\
        Musk1 & $0.2$ & $15.34$ & $15.12$ & $14.99$ & $24.66$ & $23.99$ & $23.48$ \\
    	      & $0.3$ & $15.55$ & $15.21$ &         & $24.96$ & $24.58$ &   \\
    	\hline
		        & $0.1$ & $3.96$ & $3.75$ &        & $4.80$ & $4.15$ &   \\
		Gisette & $0.2$ & $3.99$ & $3.86$ & $3.92$ & $4.44$ & $3.95$ & $4.23$ \\
		        & $0.3$ & $4.10$ & $3.88$ &        & $4.53$ & $3.92$ &   \\
		\hline
		      & $0.1$ & $3.81$ & $3.43$ &        & $4.92$ & $4.28$ &   \\
		Musk2 & $0.2$ & $3.89$ & $3.64$ & $3.84$ & $4.65$ & $4.12$ & $4.44$ \\
		      & $0.3$ & $4.13$ & $3.72$ &        & $4.49$ & $3.85$ &    \\
		\hline
		      & $0.1$ & $2.21$ & $2.13$ &        & $0.63$ & $0.56$ &  \\ 
		HTRU2 & $0.2$ & $2.22$ & $2.11$ & $2.19$ & $0.66$ & $0.53$ & $0.58$ \\ 
		      & $0.3$ & $2.25$ & $2.16$ &        & $0.60$ & $0.52$ &   \\ 	
		\hline
		      & $0.1$ & $0.71$  & $0.63$ &        & $0.63$ & $0.50$ &  \\ 
		Occup & $0.2$ & $0.73$  & $0.66$ & $0.69$ & $0.60$ & $0.45$ & $0.58$ \\ 
		      & $0.3$ & $0.82$  & $0.69$ &        & $0.66$ & $0.43$ &   \\ 			
		\hline
		       & $0.1$ & $19.25$  & $19.10$ &         & $3.33$ & $2.60$ &  \\ 
		Credit & $0.2$ & $19.26$  & $19.18$ & $19.23$ & $3.45$ & $2.71$ & $3.14$ \\ 
		       & $0.3$ & $19.29$  & $19.22$ &         & $3.43$ & $2.78$ &   \\ 			
		\hline
		     & $0.1$ & $23.10$  & $21.56$ &         & $9.58$ & $7.16$ &  \\ 
		SUSY & $0.2$ & $23.24$  & $21.34$ & $21.88$ & $9.98$ & $7.34$ & $8.89$ \\ 
		     & $0.3$ & $23.60$  & $21.99$ &         & $9.30$ & $7.23$ &   \\ 		
	    \hline
	\end{tabular}
	\end{center}
\end{table}

\subsection{Real Examples}\label{sec:real}
This subsection serves two goals. The first goal is to empirically check how accurate the M-DNN(k) and W-DNN(k) methods are compared with the oracle $K$NN and oracle OWNN methods. The second goal is to see how DNN methods compare with the bagging approach. In particular, we compare the finite-sample accuracy and stability of M-DNN, W-DNN and the oracle BNN using real data.  We choose to seperate these comparisons because bagging has a slightly different goal of increasing stability to increase accuracy, and hence  it would be a perfect benchmark to compare the classification instability of the DNN methods with.

We have retained benchmark data sets HTRU2 \citep{Lyon2016}, Gisette \citep{Guyon2005}, Musk1 \citep{Dietterich1994}, Musk2 \citep{Dietterich1997}, Occupancy \citep{Candanedo2016}, Credit \citep{Yeh2009}, and SUSY \citep{Baldi2014}, from the UCI machine learning repository \citep{Lichman2013a}. The test sample sizes are set as $\min(1000,\mbox{total sample size}/5)$. Parameters in the oracle $K$NN, BNN and OWNN are tuned using cross-validation, and the parameter $k$ in M-DNN(k), W-DNN(k) and parameter $l$ in M-DNN, W-DNN for each subsample are set using bridging formulas stated in our theorems. The empirical risk and CIS are calculated over $1000$ replications.

In Table \ref{tab:real1}, we compare the empirical risk (test error), and the speedup factor of M-DNN(k) and W-DNN(k) relative to oracle $K$NN. The latter is defined as the computing time of the oracle $K$NN divided by the time of the slower of the two DNN(k) methods. OWNN typically has similar computing time as oracle $K$NN  and hence the speed comparison with OWNN is omitted. From Table \ref{tab:real1}, we can see that the W-DNN(k) has a similar risk as the oracle $K$NN while M-DNN(k) has a somewhat larger risk.
Compared with the oracle OWNN method, which has an optimally chosen weight function, both DNN methods have a little large risks. The DNN methods have a computational advantage over oracle $K$NN and oracle OWNN, and such an advantage increases as the overall sample size increases (for a given $\gamma$). It seems that larger $\gamma$ values may induce slightly worse performance for the DNN(k) classifiers, although such an observation is not conclusive. Lastly, we note that the parameters for the oracle $K$NN  and OWNN methods are tuned, but for the DNN methods, their parameter values are based on the tuned parameters for the oracle methods and asymptotic connection formulas from the theorems. This may also slightly adds to the disadvantage of the DNN methods.

In Table \ref{tab:real2}, we compare the empirical risk (test error) and CIS of M-DNN and W-DNN with the oracle BNN (bagging).  For both criteria, we notice that 
W-DNN performs better than bagging while M-DNN performs slightly worse. Again, we stress that BNN is tuned while the parameters in DNN is found by connection formulas suggested by the theorems.

\section{Discussions}\label{sec:conclusion}
There are a couple of interesting directions to be pursued in the future. The first two are extensions to the multicategory classification problem and to high-dimensional data. The third direction is related to a realistic attack paradigm named adversarial examples that received a lot of recent attentions \citep{szegedy2013intriguing,  papernot2016limitations}. \citet{wang2017analyzing} proposed a theoretical framework for learning robustness to adversarial examples and introduced a modified 1-nearest neighbor algorithm with good robustness. This work leaves us wonder how to take advantage of the distributed nature of DNN to deal with adversarial samples.


\newcommand\invisiblesection[1]{%
  \refstepcounter{section}%
  \addcontentsline{toc}{section}{\protect\numberline{\thesection}#1}%
  \sectionmark{#1}}
\invisiblesection{Appendices}\label{sec:appendices}

\setcounter{subsection}{0}
\renewcommand{\thesubsection}{A.\Roman{subsection}}
\setcounter{equation}{0}
\renewcommand{\theequation}{A.\arabic{equation}}
\setcounter{lemma}{0}
\renewcommand{\thelemma}{A.\arabic{lemma}}

\begin{center}
\large APPENDIX
\end{center}

\subsection{Assumptions (A1) - (A4)}\label{sec:assumptions}
For a smooth function $g$, we write $\dot{g}(x)$ for its gradient vector at $x$. The following conditions are assumed throughout this paper.

(A1) The set ${\cal R}\subset \mathbb R^d$ is a compact $d$-dimensional manifold with boundary $\partial{\cal R}$.

(A2) The set ${\cal S}=\{x\in {\cal R}: \eta(x)=1/2\}$ is nonempty. There exists an open subset $U_0$ of ${\mathbb R}^d$ which contains ${\cal S}$ such that: (1) $\eta$ is continuous on $U\backslash U_0$ with $U$ an open set containing ${\cal R}$; (2) the restriction of the conditional distributions of $X$, $P_1$ and $P_0$, to $U_0$ are absolutely continuous with respect to Lebesgue measure, with twice continuously differentiable Randon-Nikodym derivatives $f_1$ and $f_0$.

(A3) There exists $\rho >0$ such that $\int_{{\mathbb R}^d}\|x\|^{\rho} d\bar{P}(x) < \infty$. In addition, for sufficiently small $\delta>0$, $\inf_{x\in {\cal R}}\bar{P}(B_{\delta}(x))/(a_d\delta^d) \ge C_0 >0$, where $a_d=\pi^{d/2}/\Gamma(1+d/2)$, $\Gamma(\cdot)$ is gamma function, and $C_0$ is a constant independent of $\delta$.

(A4) For all $x\in {\cal S}$, we have $\dot{\eta}(x)\ne 0$, and for all $x\in {\cal S}\cap \partial{\cal R}$, we have $\dot{\partial \eta}(x)\ne 0$, where $\partial \eta$ is the restriction of $\eta$ to $\partial {\cal R}$. \hfill $\blacksquare$

\subsection{Definitions of \texorpdfstring{$a(x)$, $B_1$, $B_2$, $W_{n,\beta}$ and $W_{N,\beta}$}{Lg}}\label{sec:defwnb}
For a smooth function $g$: $\mathbb{R}^d\rightarrow \mathbb{R}$, denote $g_j(x)$ as its $j$-th partial derivative at $x$ and $g_{jk}(x)$ the $(j,k)$-th element of its Hessian matrix at $x$. Let $c_{j,d}=\int_{v:\|v\|\le 1} v_j^2 dv$, $\bar{f}=\pi_1 f_1 + (1-\pi_1)f_0$. Define
\begin{equation*}
a(x)=\sum_{j=1}^d  \frac{c_{j,d}\{\eta_{j}(x)\bar{f}_j(x) + 1/2 \eta_{jj}(x)\bar{f}(x)\}}{a_d^{1+2/d} \bar{f}(x)^{1+2/d}}.
\end{equation*}
Moreover, define two distribution-related constants
\begin{eqnarray*}
B_1 = \int_{\cal S} \frac{\bar{f}(x)}{4\|\dot{\eta}(x)\|} d \textrm{Vol}^{d-1}(x),\quad B_2 = \int_{\cal S} \frac{\bar{f}(x)}{\|\dot{\eta}(x)\|} a(x)^2 d \textrm{Vol}^{d-1}(x),
\end{eqnarray*}
where $\textrm{Vol}^{d-1}$ is the natural $(d-1)$-dimensional volume measure that ${\cal S}$ inherits as a subset of $\mathbb{R}^d$. According to Assumptions (A1)-(A4) in Appendix \ref{sec:assumptions}, $B_1$ and $B_2$ are finite with $B_1>0$ and $B_2\ge 0$, with equality only when $a(x)=0$ on ${\cal S}$.

In addition, for $\beta>0$, we define $W_{n,\beta}$ as the set of $\bw_n$ satisfying:
\begin{itemize}
    \item[(w.1)] $\sum_{i=1}^n w_{ni}^2 \le n^{-\beta}$;
    \item[(w.2)] $n^{-4/d}(\sum_{i=1}^n\alpha_iw_{ni})^2\le n^{-\beta}$, where $\alpha_i=i^{1+\frac{2}{d}}-(i-1)^{1+\frac{2}{d}}$;
     \item[(w.3)] $n^{2/d}\sum_{i=k_2+1}^n w_{ni}/\sum_{i=1}^n \alpha_iw_{ni}\le 1/\log n$ with $k_2=\lceil n^{1-\beta} \rceil$;
      \item[(w.4)] $\sum_{i=k_2+1}^n w_{ni}^2/\sum_{i=1}^n w_{ni}^2 \le 1/\log n$;
       \item[(w.5)] $\sum_{i=1}^n w_{ni}^3/(\sum_{i=1}^n w_{ni}^2)^{3/2} \le 1/\log n$.
\end{itemize}
When $n$ in (w.1)--(w.5) is replaced by $N$, we can define the set $W_{N,\beta}$. \hfill $\blacksquare$

\subsection{Proof of Theorem \ref{thm:M-DNN_re}}\label{sec:pf_thm:M-DNN_re}
For the sake of simplicity, we omit $\bw_n$ in the subscript of such notations as $\widehat{\phi}_{n,s,\bw_n}^{M}$ and $S_{n,\bw_n}^{(j)}$. Write $P^\circ=\pi_1 P_1 - (1-\pi_1) P_0$. We have 
\begin{eqnarray*}
\textrm{Regret} (\widehat{\phi}_{n,s}^{M}) 
&=& R (\widehat{\phi}_{n,s}^{M})- R(\phi^{\ast})  \\
&=& \int_{{\cal R}} \pi_1 \big[  {\mathbb P}\big(\widehat{\phi}_{n,s}^{M}(x) =0 \big) -  \indi{\phi^{\ast}(x) =0} \big]  d P_1(x) \\
&~~~& +  \int_{{\cal R}} (1-\pi_1) \big[  {\mathbb P}\big(\widehat{\phi}_{n,s}^{M}(x) =1 \big) - \indi{\phi^{\ast}(x) =1} \big]  d P_0(x) \\
&=& \int_{{\cal R}} \big[  {\mathbb P}\big(\widehat{\phi}_{n,s}^{M}(x) =0 \big) - \indi{\eta(x)<1/2} \big] d P^\circ(x).
\end{eqnarray*}

Without loss of generality, we consider the $j$-th subsample of ${\cal D}$: ${\cal D}^{(j)}=\{(X_i^{(j)},Y_i^{(j)}), i=1,\ldots,n\}$. Given $X=x$, we define $(X_{(i)}^{(j)},Y_{(i)}^{(j)})$ such that $\|X_{(1)}^{(j)}-x\|\le \|X_{(2)}^{(j)}-x\|\le \ldots \le \|X_{(n)}^{(j)}-x\|$. 
Denote the estimated regression function on the $j$-th subsample as 
\begin{equation*}
S_{n}^{(j)}(x)={\textstyle\sum}_{i=1}^n w_{ni}Y_{(i)}^{(j)}.
\end{equation*}
Denote the WNN classifier on the $j$-th subsample as
\begin{equation*}
\widehat{\phi}^{(j)}_{n}(x)=\indi{ S_{n}^{(j)}(x) \geq 1/2}. 
\end{equation*} 
For any $j$ and $x$, we have ${\mathbb P}\big(S_{n}^{(j)}(x) \ge 1/2 \big)= {\mathbb P}\big(S_{n}(x) \ge 1/2 \big)$, where $S_{n}(x)$ is a generic local WNN regression function on any subsample. Hence, $\widehat{\phi}^{(j)}_{n}(x) \;(j=1,\dots,s)$ follow i.i.d. Bernoulli distribution with success probability ${\mathbb P}\big(S_{n}(x) \ge 1/2 \big)$. In particular, we have
\begin{align*}
{\mathbb E}\{\widehat{\phi}^{(j)}_{n}(x)\} 
&= {\mathbb P}\big(S_{n}(x)\ge1/2\big), \\
Var \{\widehat{\phi}^{(j)}_{n}(x)\} 
&= {\mathbb P}\big(S_{n}(x) < 1/2 \big){\mathbb P}\big(S_{n}(x) \ge 1/2 \big).
\end{align*}

Denote the average of the predictions from $s$ subsamples as 
\begin{equation*}
S_{n,s}^{M}(x)=s^{-1}{\textstyle\sum}_{j=1}^s \widehat{\phi}^{(j)}_{n}(x).
\end{equation*}
Therefore, 
\begin{align*}
{\mathbb E}\{S_{n,s}^{M}(x)\} &= {\mathbb P}\big(S_{n}(x) \geq 1/2 \big), \\
Var \{S_{n,s}^{M}(x) \} &= s^{-1} {\mathbb P}\big(S_{n}(x) < 1/2 \big){\mathbb P}\big(S_{n}(x) \ge 1/2 \big).
\end{align*}
The M-DNN classifier is defined as
$$
\widehat{\phi}_{n,s}^{M}(x)=\indi{S_{n,s}^{M}(x) \geq 1/2}.
$$
\noindent Since ${\mathbb P}\big(\widehat{\phi}_{n,s}^{M}(x) =0 \big)  =  {\mathbb P}\big(S_{n,s}^{M}(x) < 1/2 \big)$,
the regret of M-DNN becomes
\begin{eqnarray*}
\textrm{Regret} (\widehat{\phi}_{n,s}^{M}) = \int_{{\cal R}} \big\{ {\mathbb P}\big(S_{n,s}^{M}(x) < 1/2 \big) -  \indi{\eta(x)<1/2} \big\}  d P^\circ(x).
\end{eqnarray*}

In any subsample, denote the boundary ${\cal S}=\{x\in {\cal R}: \eta(x)=1/2\}$. For $\epsilon>0$, let ${\cal S}^{\epsilon\epsilon} = \{x\in {\mathbb R}^d: \eta(x)=1/2 ~\textrm{and}~ \textrm{dist}(x, {\cal S})<\epsilon \}$, where $\textrm{dist}(x, {\cal S})=\inf_{x_0\in {\cal S}} \|x-x_0\|$. We will focus on the set
$$
{\cal S}^{\epsilon} = \Big\{x_0 + t\frac{\dot{\eta}(x_0)}{\|\dot{\eta}(x_0)\|}: x_0 \in {\cal S}^{\epsilon\epsilon}, |t| < \epsilon \Big\}.
$$

Let $\mu_n(x)={\mathbb E}\{S_{n}(x)\}$, $\sigma_n^2(x)=\textrm{Var}\{S_{n}(x)\}$, and $\epsilon_n=n^{-\beta/(4d)}$. Denote $s_n^2=\sum_{i=1}^n w_{ni}^2$ and $t_n=n^{-2/d}\sum_{i=1}^n \alpha_i w_{ni}$. \citet{S12} showed that, uniformly for $\bw_n\in W_{n,\beta}$,
\begin{eqnarray}
\sup_{x\in {\cal S}^{\epsilon_n}} |\mu_n(x) - \eta(x) - a(x)t_n| &=& o(t_n),\label{step1:tn}\\
\sup_{x\in {\cal S}^{\epsilon_n}}  \big|\sigma_n^2(x)-\frac{1}{4}s_n^2\big| &=& o(s_n^2). \label{step1:sn}
\end{eqnarray}
Let $\epsilon_{n,s}^{M} = a_0t_n+b_0\frac{\log(s)}{\sqrt{s}}s_n$, where $a_0$ and $b_0$ are constants that $a_0>\frac{2|a(x_0)|}{\|\dot{\eta}(x_0)\|}$ and $b_0>\frac{\sqrt{2\pi}}{\|\dot{\eta}(x_0)\|}$, for any $x_0\in {\cal S}$.

We organize our proof in four steps. In {\it Step 1}, we decompose the integral over ${\cal R} \cap {\cal S}^{\epsilon_{n}}$ as an integral along ${\cal S}$ and an integral in the perpendicular direction; in {\it Step 2}, we bound the contribution to regret from ${\cal R} \backslash {\cal S}^{\epsilon_{n}} $; in {\it Step 3}, we bound the contribution to regret from ${\cal S}^{\epsilon_{n}} \backslash {\cal S}^{\epsilon_{n,s}^{M}} $; {\it Step 4} combines the results in previous steps and applies the normal approximation in ${\cal S}^{\epsilon_{n,s}^{M}}$ to yield the final conclusion.

{\it Step 1}: For $x_0\in {\cal S}$ and $t\in {\mathbb R}$, denote $x_0^t=x_0+t \dot{\eta}(x_0)/\|\dot{\eta}(x_0)\|$. Denote $\psi=\pi_1 f_1 - (1-\pi_1) f_0$, $\bar{f}=\pi_1 f_1 + (1-\pi_1) f_0$ as the Radon-Nikodym derivatives with respect to Lebesgue measure of the restriction of $P^\circ$ and $\bar{P}$ to ${\cal S}^{\epsilon_n}$ for large $n$ respectively.

Similar to \citet{S12}, we consider a change of variable from $x$ to $x_0^t$.
By the theory of integration on manifolds and Weyl's tube formula \citep{G04}, we have, uniformly for $\bw_n\in W_{n,\beta}$,
\begin{align*}
&\int_{{\cal R}\cap{{\cal S}^{\epsilon_n}}} \big\{ {\mathbb P}(S_{n,s}^{M}(x) < 1/2) - \indi{\eta(x)<1/2} \big\} d P^\circ(x)  \nonumber\\
=&\int_{{\cal S}} \int_{-\epsilon_n}^{\epsilon_n}\psi(x_0^t) \big\{{\mathbb P}\big(S_{n,s}^{M}(x_0^t) < 1/2\big)- \indi{t<0}\big\} dt d\textrm{Vol}^{d-1}(x_0)\{1+o(1)\}. \nonumber
\end{align*}

{\it Step 2}: Bound the contribution to regret from ${\cal R}\backslash {\cal S}^{\epsilon_{n}}$. We show that,
\begin{align*}
\sup_{\bw_n\in W_{n,\beta}} \int_{{\cal R}\backslash {\cal S}^{\epsilon_n}} \big\{ {\mathbb P}\big(S_{n,s}^{M}(x) < 1/2\big) - \indi{\eta(x)<1/2} \big\} d P^\circ(x)=o(\frac{s_n^2}{s} +t_n^2). 
\end{align*}

According to \cite{S12}, for all $M>0$, uniformly for $\bw_n\in W_{n,\beta}$ and $x\in {\cal R}\backslash {\cal S}^{\epsilon_n}$, we have
\begin{equation*}
|{\mathbb P}\big(S_{n}(x) < 1/2\big)-\mathds{1}\{\eta(x)<1/2\}|=O(n^{-M}).
\end{equation*}
Therefore, we have for $x\in {\cal R}\backslash {\cal S}^{\epsilon_n}$,
\begin{align}
\inf_{\eta(x)<1/2} {\mathbb P}\big(S_{n}(x) < 1/2\big)-1/2 &> 1/4,   \label{eta_smaller} \\
\sup_{\eta(x) \geq 1/2}{\mathbb P}\big(S_{n}(x) < 1/2 \big)-1/2 &<-1/4. \label{eta_larger}
\end{align}

Applying Hoeffding's inequality to $S_{n,s}^{M}(x)$, along with \eqref{eta_smaller} and \eqref{eta_larger}, we have
\begin{eqnarray}
&&|{\mathbb P}\big(S_{n,s}^{M}(x) < 1/2\big) -  \indi{\eta(x)<1/2} | 
\leq \exp\Big(  \frac{-2({\mathbb E}\{S_{n,s}^{M}(x)\}-1/2)^2}{\sum_{j=1}^s (1/s-0)^2} \Big) \nonumber \\
&=& \exp\Big( \frac{-2({\mathbb P}(S_{n}(x) < 1/2)-1/2)^2}{1/s} \Big) 
= o(\frac{s_n^2}{s}+t_n^2), \nonumber 
\end{eqnarray}
uniformly for $\bw_n\in W_{n,\beta}$ and $x\in {\cal R}\backslash {\cal S}^{\epsilon_n}$. This completes {\it Step 2}.

{\it Step 3}: Bound the contribution to regret from ${\cal S}^{\epsilon_{n}} \backslash {\cal S}^{\epsilon_{n,s}^{M}}$. We show that
\begin{align*}
&\sup_{\bw_n\in W_{n,\beta}}\int_{{\cal S}} \int_{(-\epsilon_n,\epsilon_n)\backslash(-\epsilon_{n,s}^{M}, \epsilon_{n,s}^{M})}\psi(x_0^t) \big\{{\mathbb P}\big(S_{n,s}^{M}(x_0^t) < 1/2\big) 
\nonumber\\
&\qquad\qquad\qquad\qquad
- \indi{t<0}\big\} dt d\textrm{Vol}^{d-1}(x_0)
= o(s_n^2/s+t_n^2). 
\end{align*}

For $x_0^t \in {\cal S}^{\epsilon_{n}}\backslash {\cal S}^{\epsilon_{n,s}^{M}}$, we have $t \not\in
 (-\epsilon_{n,s}^{M}, \epsilon_{n,s}^{M})$.
By \eqref{step1:tn}, \eqref{step1:sn} and Taylor expansion, we have
\begin{eqnarray*}
\frac{1/2 - \mu_{n}(x_0^t)}{\sigma_{n}(x_0^t)} &=& \frac{-t\|\dot{\eta}(x_0)\|- a(x_0)t_n+o(t_n)}{s_n/2+o(s_n)}. 
\end{eqnarray*}
Since $t \not\in (-\epsilon_{n,s}^{M}, \epsilon_{n,s}^{M})$ , $|t|>\epsilon_{n,s}^{M} = a_0t_n+b_0\frac{\log(s)}{\sqrt{s}}s_n$, we have
$$
\big| \frac{1/2 - \mu_{n}(x_0^t)}{\sigma_{n}(x_0^t)} \big| > \sqrt{2\pi}\frac{\log(s)}{\sqrt{s}}.
$$

If in addition, $\big|\frac{1/2 - \mu_{n}(x_0^t)}{\sigma_{n}(x_0^t)} \big| = o(1)$, then by Lemma \ref{lemma:Phi}, we have 
$$
\big| \Phi\Big(\frac{1/2 - \mu_{n}(x_0^t)}{\sigma_{n}(x_0^t)}\Big) -1/2 \big| > \frac{\log(s)}{\sqrt{s}},
$$
where $\Phi$ is the standard normal distribution function. 

Otherwise if $ \big| \frac{1/2 - \mu_{n}(x_0^t)}{\sigma_{n}(x_0^t)} \big| > c_{10}
$, where $c_{10}$ is a positive constant, then we have, for $s$ large enough,
$$
\big|  \Phi\Big(\frac{1/2 - \mu_{n}(x_0^t)}{\sigma_{n}(x_0^t)}\Big) -1/2 \big| > \frac{\log(s)}{\sqrt{s}}.
$$

In summary, for $x_0^t \in {\cal S}^{\epsilon_{n}}\backslash {\cal S}^{\epsilon_{n,s}^{M}}$,
\begin{equation}
\big|  \Phi\Big(\frac{1/2 - \mu_{n}(x_0^t)}{\sigma_{n}(x_0^t)}\Big) -1/2 \big| > \frac{\log(s)}{\sqrt{s}}. \label{phi_greater}
\end{equation}

Let $Z_i= (w_{ni}Y_{(i)} - w_{ni} \mathbb{E} \left[Y_{(i)}\right])/\sigma_{n}(x)$ and $W=\sum_{i=1}^n Z_i$. Note that $\mathbb{E}(Z_i)=0$, $\textrm{Var}(Z_i)<\infty$ and $\textrm{Var}(W)=1$. The nonuniform Berry-Esseen Theorem \citep{GP12} implies that there exists a constant $c_{11}>0$ such that 
$$
\Big|\mathbb{P}(W\le By) - \Phi(y)\Big| \le \frac{c_{11}A}{B^3( 1 + |y|^3)},
$$
where $A={\textstyle\sum}_{i=1}^n E|Z_i|^3 \;\;{\rm and}\;\; B=\big({\textstyle\sum}_{i=1}^n E|Z_i|^2)^{1/2}$. In our case, 
\begin{align*}
A&=\sum_{i=1}^n \mathbb{E}|\frac{w_{ni}Y_i - w_{ni} \mathbb{E} [Y_i]}{\sigma_{n}(x)}|^3 \leq \sum_{i=1}^n \frac{2|w_{ni}|^3}{(s_n/2)^3}=\frac{16\sum_{i=1}^n w_{ni}^3}{(\sum_{i=1}^n w_{ni}^2)^{3/2}}, \\
B&=({\textstyle\sum}_{i=1}^n \textrm{Var}(Z_i))^{1/2}=\sqrt{ \textrm{Var}(W)}=1.
\end{align*}
Let $c_{12}=16c_{11}$, we have
\begin{align*}
\sup_{x_0\in {\cal S}}\sup_{t\in [-\epsilon_n,\epsilon_n]}\Big|\mathbb{P}\Big(\frac{S_{n}(x_0^t)-\mu_{n}(x_0^t)}{\sigma_{n}(x_0^t)}\le y \Big) - \Phi(y) \Big|\le \frac{\sum_{i=1}^n w_{ni}^3}{(\sum_{i=1}^n w_{ni}^2)^{3/2}}\frac{c_{12}}{1 + |y|^3}.
\end{align*}
Setting $y=\frac{1/2 - \mu_{n}(x_0^t)}{\sigma_{n}(x_0^t)}$, we have
\begin{align}
&\sup_{x_0\in {\cal S}}\sup_{t\in [-\epsilon_n,\epsilon_n]} \Big|{\mathbb P}\big(S_{n}(x_0^t) < 1/2\big)- \Phi\Big(\frac{1/2 - \mu_{n}(x_0^t)}{\sigma_{n}(x_0^t)}\Big) \Big| \label{BE_thm_re_M}\\
\le& \frac{c_{12}\sum_{i=1}^n w_{ni}^3}{(\sum_{i=1}^n w_{ni}^2)^{3/2}}=o\Big(\frac{1}{\sqrt{s}(\log(s))^2}\Big). \nonumber
\end{align}
The last equality holds by \eqref{eq:extra_condition}. 

By \eqref{phi_greater} and \eqref{BE_thm_re_M}, we have, when $x_0^t \in {\cal S}^{\epsilon_{n}}\backslash {\cal S}^{\epsilon_{n,s}^{M}}$,
\begin{eqnarray}
\big|{\mathbb P}\big(S_{n}(x_0^t) < 1/2\big) -1/2 \big| > \log(s)/\sqrt{s}.
\end{eqnarray}

Applying Hoeffding's inequality to $S_{n,s}^{M}(x_0^t)$, we have
\begin{align}
&|{\mathbb P}(S_{n,s}^{M}(x_0^t) < 1/2) -  \indi{t<0}| 
\leq \exp\Big[\frac{-2({\mathbb E}\{S_{n,s}^{M}(x_0^t)\}-1/2)^2}{\sum_{j=1}^s (1/s-0)^2} \Big] \nonumber\\ 
=&\exp\big[ -2s({\mathbb P}(S_{n}(x_0^t) < 1/2)-1/2)^2 \big]<s^{-2\log(s)}
=o(s_n^2/s+t_n^2), \nonumber
\end{align} 
uniformly for $\bw_n\in W_{n,\beta}$ and $x_0^t\in {\cal S}^{\epsilon_{n}}\backslash {\cal S}^{\epsilon_{n,s}^M} $. This completes {\it Step 3}.

{\it Step 4}: In the end, we will show
\begin{eqnarray*}
&&\int_{{\cal S}} \int_{-\epsilon_{n,s}^{M}}^{\epsilon_{n,s}^{M}} \psi(x_0^t) \big\{{\mathbb P}\big(S_{n,s}^{M}(x_0^t) < 1/2\big) - \indi{t<0}\big\} dt d\textrm{Vol}^{d-1}(x_0) \nonumber\\
&=& B_1 \frac{\pi}{2s}s_n^2 + B_2 t_n^2 + o(\frac{s_n^2}{s}+t_n^2).
\end{eqnarray*}

Applying Taylor expansion, we have, for $x_0 \in {\cal S}$,
\begin{eqnarray}
\psi(x_0^t) &=& \psi(x_0) + \dot{\psi}(x_0)^T(x_0^t-x_0) + o(x_0^t-x_0) \label{eq:psi_taylor} \\
&=& \dot{\psi}(x_0)^T\frac{\dot{\eta}(x_0)}{\|\dot{\eta}(x_0)\|} t + o(t) \nonumber\\
&=& \|\dot{\psi}(x_0)\| t +o(t), \nonumber
\end{eqnarray}
where the above second equality holds by definition of $x_0^t$, and the third equality holds by Lemma \ref{lemma:dot}. Hence, 
\begin{eqnarray}
&&\int_{{\cal S}} \int_{-\epsilon_{n,s}^{M}}^{\epsilon_{n,s}^{M}} \psi(x_0^t) \big\{{\mathbb P}\big(S_{n,s}^{M}(x_0^t) < 1/2\big) - \indi{t<0}\big\} dt d\textrm{Vol}^{d-1}(x_0) \label{Taylor11_re_M}\\
&=& \int_{{\cal S}} \int_{-\epsilon_{n,s}^{M}}^{\epsilon_{n,s}^{M}} t\|\dot{\psi}(x_0)\|  \big\{{\mathbb P}\big(S_{n,s}^{M}(x_0^t) < 1/2\big) \nonumber\\
&&\qquad\qquad\qquad\qquad
 - \indi{t<0}\big\} dt d\textrm{Vol}^{d-1}(x_0) \{1+o(1) \}. \nonumber
\end{eqnarray}


Next, we decompose
\begin{align}
&\int_{{\cal S}} \int_{-\epsilon_{n,s}^{M}}^{\epsilon_{n,s}^{M}}t\|\dot{\psi}(x_0)\| \big\{ {\mathbb P}\big(S_{n,s}^{M}(x_0^t) < 1/2\big) - \indi{t<0} \big\} dt d\textrm{Vol}^{d-1}(x_0)\label{eq:decompose_11_re_M} \\
=& \int_{{\cal S}} \int_{-\epsilon_{n,s}^{M}}^{\epsilon_{n,s}^{M}} t\|\dot{\psi}(x_0)\| \big\{ \Phi\Big[\frac{\sqrt{s}\big(1/2-{\mathbb P}\big(S_{n}(x_0^t)\geq1/2\big)\big)}{\sqrt{{\mathbb P}\big(S_{n}(x_0^t)<1/2\big){\mathbb P}\big(S_{n}(x_0^t) \ge 1/2\big)}} \Big] \nonumber \\
& \qquad\qquad\qquad\qquad\qquad
- \indi{t<0}  \big\} dt d\textrm{Vol}^{d-1}(x_0) +R_{11}. \nonumber
\end{align}
If $|1/2-{\mathbb P}\big(S_{n}(x_0^t)<1/2\big)| \le \log(s)/\sqrt{s} $, by the uniform Berry-Esseen Theorem \citep{lehmann2004elements}, there exists a constant $c_{13}>0$ such that
\begin{eqnarray*}
&&\Big|{\mathbb P}\Big(\frac{\sqrt{s}\big[S_{n,s}^{M}(x_0^t)- {\mathbb P}\big(S_{n}(x_0^t)\ge 1/2\big)\big]}{\sqrt{{\mathbb P}\big(S_{n}(x_0^t)<1/2\big){\mathbb P}\big(S_{n}(x_0^t) \ge 1/2\big)}} < y\Big)-\Phi\big(y \big)\Big| \\
&\leq& \frac{c_{13}}{\sqrt{s}}\frac{{\mathbb E}\big|\widehat{\phi}^{(j)}_{n}(x_0^t)-{\mathbb P}\big(S_{n}(x_0^t)\ge 1/2\big)\big|^3  }{\big[{\mathbb P}\big(S_{n}(x_0^t)<1/2\big){\mathbb P}\big(S_{n}(x_0^t) \ge 1/2\big)\big]^{3/2}} \le \frac{8c_{13}}{\sqrt{s}} =O\Big(\frac{1}{\sqrt{s}}\Big).\\
\end{eqnarray*}
Setting $y = \frac{\sqrt{s}(1/2-{\mathbb P}(S_{n}(x_0^t)\geq1/2))}{\sqrt{{\mathbb P}(S_{n}(x_0^t)<1/2){\mathbb P}(S_{n}(x_0^t) \ge 1/2)}}$, we have
\begin{align*}
\Big| {\mathbb P}\big(S_{n,s}^{M}(x_0^t) < 1/2\big)- \Phi\Big[\frac{\sqrt{s}\big(1/2-{\mathbb P}\big(S_{n}(x_0^t)\geq1/2\big)\big)}{\sqrt{{\mathbb P}\big(S_{n}(x_0^t)<1/2\big){\mathbb P}\big(S_{n}(x_0^t) \ge 1/2\big)}} \Big] \Big| = O\Big(\frac{1}{\sqrt{s}}\Big). \nonumber
\end{align*}
In addition, if $|1/2-{\mathbb P}\big(S_{n}(x_0^t)\geq1/2\big)| > \frac{\log(s)}{\sqrt{s}} $, applying Hoeffding's inequality and Lemma \ref{lemma:feller} , we have
\begin{align*}
&|{\mathbb P}(S_{n,s}^{M}(x_0^t) < 1/2) - \indi{{\mathbb P}\big(S_{n}(x_0^t)\geq1/2\big)<1/2}|  \nonumber\\
=&\exp\big[ -2s(1/2-{\mathbb P}(S_{n}(x_0^t)\ge 1/2))^2 \big]\le \exp(-2[\log(s)]^2) =o\Big(\frac{1}{\sqrt{s}}\Big), \nonumber  \\
&\Big| \Phi\Big[\frac{\sqrt{s}\big(1/2-{\mathbb P}\big(S_{n}(x_0^t)\geq1/2\big)\big)}{\sqrt{{\mathbb P}\big(S_{n}(x_0^t)<1/2\big){\mathbb P}\big(S_{n}(x_0^t) \ge 1/2\big)}} \Big] -  \indi{{\mathbb P}\big(S_{n}(x_0^t)\geq1/2\big)<1/2}  \Big| \nonumber\\
\le& \frac{1}{2\log(s)}\frac{e^{-[2\log(s)]^2/2}}{\sqrt{2\pi}}=o\Big(\frac{1}{\sqrt{s}}\Big). \nonumber
\end{align*}
In this case, we have
\begin{align*}
&\Big| {\mathbb P}\big(S_{n,s}^{M}(x_0^t) < 1/2\big)- \Phi\Big[\frac{\sqrt{s}\big(1/2-{\mathbb P}\big(S_{n}(x_0^t)\geq1/2\big)\big)}{\sqrt{{\mathbb P}\big(S_{n}(x_0^t)<1/2\big){\mathbb P}\big(S_{n}(x_0^t) \ge 1/2\big)}} \Big] \Big|\\
\le &\big| {\mathbb P}\big(S_{n,s}^{M}(x_0^t) < 1/2\big)-\indi{{\mathbb P}\big(S_{n}(x_0^t)\geq1/2\big)<1/2}\big| \\
&+ \Big| \Phi\Big[\frac{\sqrt{s}\big(1/2-{\mathbb P}\big(S_{n}(x_0^t)\geq1/2\big)\big)}{\sqrt{{\mathbb P}\big(S_{n}(x_0^t)<1/2\big){\mathbb P}\big(S_{n}(x_0^t) \ge 1/2\big)}} \Big] -\indi{{\mathbb P}\big(S_{n}(x_0^t)\geq1/2\big)<\frac{1}{2}}\Big| \\
=& o\Big(\frac{1}{\sqrt{s}}\Big).  
\end{align*}
In summary, we have
\begin{eqnarray}
&&\sup_{x_0\in {\cal S}}\sup_{t\in [-\epsilon_{n,s}^{M},\epsilon_{n,s}^{M}]} \Big| {\mathbb P}\big(S_{n,s}^{M}(x_0^t) < 1/2\big) \label{Edgeworth}\\
&&\qquad\qquad
- \Phi\Big[\frac{\sqrt{s}\big(1/2-{\mathbb P}\big(S_{n}(x_0^t)\geq1/2\big)\big)}{\sqrt{{\mathbb P}\big(S_{n}(x_0^t)<1/2\big){\mathbb P}\big(S_{n}(x_0^t) \ge 1/2\big)}} \Big] \Big| = O\Big(\frac{1}{\sqrt{s}}\Big).  \nonumber
\end{eqnarray}
Thus, we have 
\begin{align*}
|R_{11}|\le&\int_{{\cal S}} \int_{-\epsilon_{n,s}^{M}}^{\epsilon_{n,s}^{M}}|t|\|\dot{\psi}(x_0)\| \Big| {\mathbb P}\big(S_{n,s}^{M}(x_0^t) < 1/2\big)\\
&
-\Phi\Big[\frac{\sqrt{s}\big(1/2-{\mathbb P}\big(S_{n}(x_0^t)\geq1/2\big)\big)}{\sqrt{{\mathbb P}\big(S_{n}(x_0^t)<1/2\big){\mathbb P}\big(S_{n}(x_0^t) \ge 1/2\big)}} \Big] \Big|  dt d\textrm{Vol}^{d-1}(x_0) \nonumber\\
\leq &  O(\frac{1}{\sqrt{s}}) \int_{{\cal S}} \int_{-\epsilon_{n,s}^{M}}^{\epsilon_{n,s}^{M}}|t|\|\dot{\psi}(x_0)\|  dt d\textrm{Vol}^{d-1}(x_0) = o(t_n^2+\frac{1}{s}s_n^2).
\end{align*}

Next, we decompose
\begin{align}
&\int_{{\cal S}}\int_{-\epsilon_{n,s}^{M}}^{\epsilon_{n,s}^{M}}t\|\dot{\psi}(x_0)\| \big\{ \Phi\Big[\frac{\sqrt{s}\big(1/2-{\mathbb P}\big(S_{n}(x_0^t)\ge1/2\big)\big)}{\sqrt{{\mathbb P}\big(S_{n}(x_0^t)<1/2\big){\mathbb P}\big(S_{n}(x_0^t) \ge 1/2\big)}} \Big] \label{eq:decompose_12_re_M}\\
&\qquad\qquad\qquad\qquad
- \indi{t<0} \big\} dt d\textrm{Vol}^{d-1}(x_0)\nonumber\\
=& \int_{{\cal S}} \int_{-\epsilon_{n,s}^{M}}^{\epsilon_{n,s}^{M}}t\|\dot{\psi}(x_0)\|\big \{\Phi\big[2\sqrt{s}\big({\mathbb P}\big(S_{n}(x_0^t)<1/2\big)-1/2\big) \big]\nonumber  \\
&\qquad\qquad\qquad\qquad
-\indi{t<0} \big\} dt d\textrm{Vol}^{d-1}(x_0) +R_{12}. \nonumber
\end{align}
If $|{\mathbb P}\big(S_{n}(x_0^t)<1/2\big)-1/2| \le \log(s)/\sqrt{s}$, along with Lemma \ref{lemma:mean_value}, we have
\begin{eqnarray}
&&\big| \Phi\Big[\frac{\sqrt{s}\big({\mathbb P}\big(S_{n}(x_0^t)<1/2\big)-1/2\big)}{\sqrt{{\mathbb P}\big(S_{n}(x_0^t)<1/2\big){\mathbb P}\big(S_{n}(x)\ge1/2\big)}} \Big] \nonumber\\ 
&&\qquad\qquad\qquad\qquad
-\Phi\big[2\sqrt{s}\big({\mathbb P}\big(S_{n}(x_0^t)<1/2\big)-1/2\big) \big] \big|  \nonumber\\
&\leq& \sqrt{s}\big|{\mathbb P}\big(S_{n}(x_0^t)<1/2\big)-1/2\big| \Big| \frac{1}{\sqrt{1/4+O(\frac{\log(s)}{\sqrt{s}})}}-2 \Big| \nonumber\\
&\le&  \sqrt{s}\frac{\log(s)}{\sqrt{s}} \frac{1-2\sqrt{1/4+O(\frac{\log(s)}{\sqrt{s}})}}{\sqrt{1/4+O(\frac{\log(s)}{\sqrt{s}})}}=O(\frac{(\log(s))^2}{\sqrt{s}}). \nonumber
\end{eqnarray} 
In addition, if $|{\mathbb P}\big(S_{n}(x_0^t)\geq1/2\big)-1/2| > \frac{\log(s)}{\sqrt{s}} $, applying Lemma \ref{lemma:feller}, we have
\begin{align*}
&\Big| 1-  \Phi\Big[\frac{\sqrt{s}\big|{\mathbb P}\big(S_{n}(x_0^t)<1/2\big)-1/2\big|}{\sqrt{{\mathbb P}\big(S_{n}(x_0^t)<1/2\big){\mathbb P}\big(S_{n}(x_0^t) \ge 1/2\big)}}\Big]  \Big| \nonumber\\
\le& \Big| 1-  \Phi\big(2\log(s) \big)  \Big| \le \frac{1}{2\log(s)}\frac{e^{-[2\log(s)]^2/2}}{\sqrt{2\pi}}=o\Big(\frac{1}{\sqrt{s}}\Big), \nonumber \\
&\big| 1- \Phi\big[2\sqrt{s}\big|{\mathbb P}\big(S_{n}(x_0^t)<1/2\big)-1/2\big| \big] \big|  \\
\le& \big| 1-  \Phi\big(2\log(s) \big)  \big| \le \frac{1}{2\log(s)}\frac{e^{-[2\log(s)]^2/2}}{\sqrt{2\pi}}=o\Big(\frac{1}{\sqrt{s}}\Big). \nonumber 
\end{align*}
In this case, we have
\begin{align*}
&\big| \Phi\Big[\frac{\sqrt{s}\big({\mathbb P}\big(S_{n}(x_0^t)<1/2\big)-1/2\big)}{\sqrt{{\mathbb P}\big(S_{n}(x_0^t)<1/2\big){\mathbb P}\big(S_{n}(x)\ge1/2\big)}} \Big] \\
&\qquad\qquad
-\Phi\big[2\sqrt{s}\big({\mathbb P}\big(S_{n}(x_0^t)<1/2\big)-1/2\big) \big] \big| \nonumber \\
\le& \big| 1-  \Phi\Big[\frac{\sqrt{s}\big|{\mathbb P}\big(S_{n}(x_0^t)<1/2\big)-1/2\big|}{\sqrt{{\mathbb P}\big(S_{n}(x_0^t)<1/2\big){\mathbb P}\big(S_{n}(x_0^t) \ge 1/2\big)}}\Big]  \big| \\
&\qquad\qquad
+\big| 1- \Phi\big[2\sqrt{s}\big|{\mathbb P}\big(S_{n}(x_0^t)<1/2\big)-1/2\big| \big] \big|=o(1/\sqrt{s}).
\end{align*}  
In summary, we have
\begin{eqnarray}
&&\sup_{x_0\in {\cal S}}\sup_{t\in [-\epsilon_{n,s}^{M},\epsilon_{n,s}^{M}]}\big| \Phi\Big[\frac{\sqrt{s}\big({\mathbb P}\big(S_{n}(x_0^t)<1/2\big)-1/2\big)}{\sqrt{{\mathbb P}\big(S_{n}(x_0^t)<1/2\big){\mathbb P}\big(S_{n}(x)\ge1/2\big)}} \Big] \label{P_P}\\
&&\qquad\qquad
-\Phi\big[2\sqrt{s}\big({\mathbb P}\big(S_{n}(x_0^t)<1/2\big)-1/2\big) \big] \big| =O((\log(s))^2/\sqrt{s}). \nonumber
\end{eqnarray}  
Therefore,
\begin{align*}
|R_{12}| \le & \int_{{\cal S}}\int_{-\epsilon_{n,s}^{M}}^{\epsilon_{n,s}^{M}}|t|\|\dot{\psi}(x_0)\| \big| \Phi\big[\frac{\sqrt{s}\big(1/2-{\mathbb P}\big(S_{n}(x_0^t) \geq 1/2\big)\big)}{\sqrt{{\mathbb P}\big(S_{n}(x_0^t)<1/2\big){\mathbb P}\big(S_{n}(x_0^t)\ge1/2\big)}} \big] \nonumber\\
&\qquad 
-\Phi\big[2\sqrt{s}\big({\mathbb P}\big(S_{n}(x_0^t)<1/2\big)-1/2\big) \big] \big| dt d\textrm{Vol}^{d-1}(x_0)\nonumber \\
=& \int_{{\cal S}}\int_{-\epsilon_{n,s}^{M}}^{\epsilon_{n,s}^{M}}|t|\|\dot{\psi}(x_0)\| \big| \Phi\big[\frac{\sqrt{s}\big({\mathbb P}\big(S_{n}(x_0^t)<1/2\big)-1/2\big)}{\sqrt{{\mathbb P}\big(S_{n}(x_0^t)<1/2\big){\mathbb P}\big(S_{n}(x_0^t)\ge1/2\big)}} \big] \nonumber\\
&\qquad 
-\Phi\big[2\sqrt{s}\big({\mathbb P}\big(S_{n}(x_0^t)<1/2\big)-1/2\big) \big] \big| dt d\textrm{Vol}^{d-1}(x_0) \nonumber\\
\leq&  O(\frac{(\log(s))^2}{\sqrt{s}})\int_{{\cal S}}\int_{-\epsilon_{n,s}^{M}}^{\epsilon_{n,s}^{M}}|t|\|\dot{\psi}(x_0)\| dt d\textrm{Vol}^{d-1}(x_0) =o(t_n^2+\frac{1}{s}s_n^2). \nonumber
\end{align*}

Next, we decompose
\begin{align}
&\int_{{\cal S}} \int_{-\epsilon_{n,s}^{M}}^{\epsilon_{n,s}^{M}}t\|\dot{\psi}(x_0)\| \big\{ \Phi\big[2\sqrt{s}\big({\mathbb P}\big(S_{n}(x_0^t)<1/2\big)-1/2\big) \big] \label{eq:decompose_13_re_M}\\
&\qquad\qquad\qquad\qquad
- \indi{t<0} \big\} dt d\textrm{Vol}^{d-1}(x_0) \nonumber\\
=& \int_{{\cal S}} \int_{-\epsilon_{n,s}^{M}}^{\epsilon_{n,s}^{M}}t\|\dot{\psi}(x_0)\| \big\{ \Phi\big[2\sqrt{s}\big(\Phi\big(\frac{1/2-\mu_{n}(x_0^t)}{\sigma_{n}(x_0^t)}\big)-1/2\big)\big]  \nonumber\\
&\qquad\qquad\qquad\qquad
- \indi{t<0} \big\} dt d\textrm{Vol}^{d-1}(x_0) +R_{13}. \nonumber
\end{align}
Applying Lemma \ref{lemma:mean_value} and \eqref{BE_thm_re_M}, we have
\begin{align}
&\sup_{x_0\in {\cal S}}\sup_{t\in [-\epsilon_{n,s}^{M},\epsilon_{n,s}^{M}]} \big|\Phi\big[2\sqrt{s}\big({\mathbb P}\big(S_{n}(x_0^t)<1/2\big)-1/2\big) \big]  \label{R13_temp}\\
&\qquad\qquad\qquad\qquad
-\Phi\big[2\sqrt{s}\big(\Phi\big(\frac{1/2-\mu_{n}(x_0^t)}{\sigma_{n}(x_0^t)}\big)-1/2\big) \big] \big|\nonumber\\
\leq& \sup_{x_0\in {\cal S}}\sup_{t\in [-\epsilon_{n,s}^{M},\epsilon_{n,s}^{M}]} \sqrt{s} \big| {\mathbb P}\big(S_{n}(x_0^t)<1/2\big)-\Phi\big(\frac{1/2-\mu_{n}(x_0^t)}{\sigma_{n}(x_0^t)}\big)  \big| \nonumber\\
\le& o\Big(\sqrt{s}\frac{1}{\sqrt{s}(\log(s))^2}\Big) = o((\log(s))^{-2}). \nonumber
\end{align}
Hence,
\begin{align*}
|R_{13}|
\le&  \int_{{\cal S}} \int_{-\epsilon_{n,s}^{M}}^{\epsilon_{n,s}^{M}}|t|\|\dot{\psi}(x_0)\|\Big|\Phi\Big[2\sqrt{s}\big({\mathbb P}\big(S_{n}(x_0^t)<1/2\big)-1/2\big) \Big] \\
&\qquad\qquad\qquad
-\Phi\Big[2\sqrt{s}\big(\Phi\big(\frac{1/2-\mu_{n}(x_0^t)}{\sigma_{n}(x_0^t)}\big)-1/2\big) \big) \Big|  dt d\textrm{Vol}^{d-1}(x_0)\\
=&o((\log(s))^{-2}) \int_{{\cal S}} \int_{-\epsilon_{n,s}^{M}}^{\epsilon_{n,s}^{M}} |t|\|\dot{\psi}(x_0)\| dt d\textrm{Vol}^{d-1}(x_0)
=o(s_n^2/s+t_n^2).
\end{align*}

Next, we decompose
\begin{align}
&\int_{{\cal S}} \int_{-\epsilon_{n,s}^{M}}^{\epsilon_{n,s}^{M}}t\|\dot{\psi}(x_0)\| \big\{ \Phi\big[2\sqrt{s}\big(\Phi\big(\frac{1/2-\mu_{n}(x_0^t)}{\sigma_{n}(x_0^t)}\big)-1/2\big) \big] \label{eq:decompose_14_re_M}\\
&\qquad\qquad\qquad
- \indi{t<0}
 \big\} dt d\textrm{Vol}^{d-1}(x_0) \nonumber\\
=& \int_{{\cal S}} \int_{-\epsilon_{n,s}^{M}}^{\epsilon_{n,s}^{M}}t\|\dot{\psi}(x_0)\|\big\{\Phi\big(\frac{1/2-\mu_{n}(x_0^t)}{\sqrt{\pi/(2s)}\sigma_{n}(x_0^t)} \big) \nonumber\\
&\qquad\qquad\qquad
-\indi{t<0}
 \big\} dt d\textrm{Vol}^{d-1}(x_0) +R_{14}. \nonumber
\end{align}
If $|\frac{1/2-\mu_{n}(x_0^t)}{\sigma_{n}(x_0^t)}|\le \frac{\log(s)}{\sqrt{s}}$, applying Lemma \ref{lemma:Phi} and Lemma \ref{lemma:mean_value}, we have, for large $s$,
\begin{align*}
&\Big|\Phi\Big[2\sqrt{s}\big(\Phi\big(\frac{1/2-\mu_{n}(x_0^t)}{\sigma_{n}(x_0^t)}\big)-1/2\big)\Big]
-\Phi\Big(\frac{1/2-\mu_{n}(x_0^t)}{\sqrt{\pi/(2s)}\sigma_{n}(x_0^t)}\Big) \Big| \\
\le& \sqrt{s}\Big|\Phi\Big(\frac{1/2-\mu_{n}(x_0^t)}{\sigma_{n}(x_0^t)}\Big)-1/2-\frac{1}{\sqrt{2\pi}}\Big(\frac{1/2-\mu_{n}(x_0^t)}{\sigma_{n}(x_0^t)}\Big)\Big| \\
=&  O\Big(\sqrt{s}\big(\frac{1/2-\mu_{n}(x_0^t)}{\sigma_{n}(x_0^t)}\big)^3\Big)   =O\Big(\sqrt{s}(\frac{\log(s)}{\sqrt{s}})^3\Big)=o\Big(\frac{1}{\sqrt{s}}\Big).
\end{align*}
In addition, if $|\frac{1/2-\mu_{n}(x_0^t)}{\sigma_{n}(x_0^t)}| >\frac{\log(s)}{\sqrt{s}}$, applying mean value theorem, there exists $x_0 \in (0,\frac{\log(s)}{\sqrt{s}})$ such that, for large $s$
\begin{eqnarray*}
&&\Phi\Big(\Big|\frac{1/2-\mu_{n}(x_0^t)}{\sigma_{n}(x_0^t)}\Big|\Big)-1/2> \Phi\Big(\frac{\log(s)}{\sqrt{s}}\Big)-\Phi(0) \\
&=& \frac{\log(s)}{\sqrt{s}} \frac{1}{\sqrt{2\pi}}\exp(-x_0^2/2)
> \frac{\log(s)}{\sqrt{s}} \frac{1}{\sqrt{2\pi}}\exp\Big(-(\frac{\log(s)}{\sqrt{s}})^2/2\Big)  > \frac{\log(s)}{4\sqrt{s}}.
\end{eqnarray*}
In this case, applying Lemma \ref{lemma:feller}, we have for large $s$
\begin{align*}
&1-\Phi\Big(\Big|\frac{1/2-\mu_{n}(x_0^t)}{\sqrt{\pi/(2s)}\sigma_{n}(x_0^t)}\Big|\Big) < 1-\Phi\big(\sqrt{2/\pi}\log(s)\big)  \nonumber\\ 
\le& \frac{1}{\sqrt{2/\pi}\log(s)}\frac{e^{-[\sqrt{2/\pi}\log(s)]^2/2}}{\sqrt{2\pi}}=o\Big(\frac{1}{\sqrt{s}}\Big)\;\;{\rm and}\;\;\\ 
&1-\Phi\Big[2\sqrt{s}\Big(\Phi\Big(\Big|\frac{1/2-\mu_{n}(x_0^t)}{\sigma_{n}(x_0^t)}\Big|\Big)-1/2\Big)\Big] 
< 1-\Phi\big((1/2)\log(s)\big) \nonumber \\ 
\le& \frac{1}{(1/2)\log(s)}\frac{e^{-[(1/2)\log(s)]^2/2}}{\sqrt{2\pi}}=o\Big(\frac{1}{\sqrt{s}}\Big).  
\end{align*}
Therefore,
\begin{eqnarray*}
&&\Big|\Phi\Big[2\sqrt{s}\big(\Phi\big(\frac{1/2-\mu_{n}(x_0^t)}{\sigma_{n}(x_0^t)}\big)-1/2\big)\Big]-\Phi\Big(\frac{1/2-\mu_{n}(x_0^t)}{\sqrt{\pi/(2s)}\sigma_{n}(x_0^t)}\Big) \Big| \nonumber\\
&=&\Big|\Phi\Big[2\sqrt{s}\Big(\Phi\Big(\Big|\frac{1/2-\mu_{n}(x_0^t)}{\sigma_{n}(x_0^t)}\Big|\Big)-1/2\Big)\Big]-\Phi\Big(\Big|\frac{1/2-\mu_{n}(x_0^t)}{\sqrt{\pi/(2s)}\sigma_{n}(x_0^t)}\Big|\Big) \Big| \nonumber\\
&\le &  \Big|1-\Phi\Big[2\sqrt{s}\Big(\Phi\Big(\Big|\frac{1/2-\mu_{n}(x_0^t)}{\sigma_{n}(x_0^t)}\Big|\Big)-1/2\Big)\Big] \Big| \nonumber\\
&& \qquad\qquad\qquad
+ \Big|1-\Phi\Big(\Big|\frac{1/2-\mu_{n}(x_0^t)}{\sqrt{\pi/(2s)}\sigma_{n}(x_0^t)}\Big|\Big) \Big|
= o\Big(\frac{1}{\sqrt{s}}\Big).
\end{eqnarray*}
In summary, we have,
\begin{align}
&\sup_{x_0\in {\cal S}}\sup_{t\in [-\epsilon_{n,s}^{M},\epsilon_{n,s}^{M}]} \Big|\Phi\Big[2\sqrt{s}\big(\Phi\big(\frac{1/2-\mu_{n}(x_0^t)}{\sigma_{n}(x_0^t)}\big)-1/2\big)\Big]   \label{Normal_CDF_at_0} \\
&\qquad\qquad\qquad\qquad\qquad
-\Phi\Big(\frac{1/2-\mu_{n}(x_0^t)}{\sqrt{\pi/(2s)}\sigma_{n}(x_0^t)}\Big) \Big| = o\Big(\frac{1}{\sqrt{s}}\Big).\nonumber
\end{align}
Therefore, 
\begin{align*}
|R_{14}|
\le&\int_{{\cal S}} \int_{-\epsilon_{n,s}^{M}}^{\epsilon_{n,s}^{M}}|t|\|\dot{\psi}(x_0)\|\Big|\Phi\Big[2\sqrt{s}\big(\Phi\big(\frac{1/2-\mu_{n}(x_0^t)}{\sigma_{n}(x_0^t)}\big)-1/2\big)\Big]\\
&\qquad\qquad
-\Phi\Big(\frac{1/2-\mu_{n}(x_0^t)}{\sqrt{\pi/(2s)}\sigma_{n}(x_0^t)}\Big) \Big| dt d\textrm{Vol}^{d-1}(x_0) \\
\leq& o\Big(\frac{1}{\sqrt{s}}\Big)\int_{{\cal S}} \int_{-\epsilon_{n,s}^{M}}^{\epsilon_{n,s}^{M}}|t|\|\dot{\psi}(x_0)\| dt d\textrm{Vol}^{d-1}(x_0) =o(t_n^2+\frac{1}{s}s_n^2). 
\end{align*} 

Next, we decompose
\begin{eqnarray}
&&\int_{{\cal S}} \int_{-\epsilon_{n,s}^{M}}^{\epsilon_{n,s}^{M}} t\|\dot{\psi}(x_0)\| \big\{\Phi\Big(\frac{1/2 - \mu_{n}(x_0^t)}{\sqrt{\pi/(2s)}\sigma_{n}(x_0^t)} \Big) - \indi{t<0}\big\} dt d\textrm{Vol}^{d-1}(x_0) \label{eq:decompose_15_re_M} \\
&=& \int_{{\cal S}} \int_{-\epsilon_{n,s}^{M}}^{\epsilon_{n,s}^{M}} t\|\dot{\psi}(x_0)\| \big\{\Phi\Big(\frac{-2t\|\dot{\eta}(x_0)\|- 2a(x_0)t_n}{\sqrt{\pi/(2s)}s_{n}} \Big) \nonumber \\
&&\qquad\qquad\qquad\qquad\qquad\qquad\qquad
- \indi{t<0}\big\} dt d\textrm{Vol}^{d-1}(x_0) + R_{15}. \nonumber
\end{eqnarray}
Denote $r=t/s_n$ and $r_{x_0}=\frac{-a(x_0)t_n}{\|\dot{\eta}(x_0)s_n\|}$. According to $(\ref{step1:tn})$ and $(\ref{step1:sn})$, for a sufficiently small $\epsilon\in (0,\inf_{x_0\in {\cal S}} \|\dot{\eta}(x_0)\|)$ and a large $n$, for all $\bw_n\in W_{n,\beta}$, $x_0\in {\cal S}$ and $r\in[-\epsilon_n/s_n,\epsilon_n/s_n]$, \citet{S12} showed that
$$
\Big| \frac{1/2-\mu_n(x_0^{rs_n})}{\sigma_n(x_0^{rs_n})} -[-2\|\dot{\eta}(x_0)\|(r-r_{x_0})] \Big| \le \epsilon^2(|r|+t_n/s_n).
$$ 
To adapt this to our setting, we need to scale some terms properly.
Let $r^{M}=r\sqrt{2s/\pi}$, $s_{n,s}^{M}=s_n\sqrt{\pi/(2s)}$ and $r^{M}_{x_0}=r_{x_0}\sqrt{2s/\pi}=\sqrt{\frac{2s}{\pi}}\frac{-a(x_0)t_n}{\|\dot{\eta}(x_0)s_n\|}$, we have, when $r^{M}\in [-\epsilon_n/s_{n,s}^{M},\epsilon_n/s_{n,s}^{M}]$,
\begin{eqnarray*}
&&\Big| \frac{1/2-\mu_{n}(x_0^{r^{M}s_{n,s}^{M}}) }{\sqrt{\pi/(2s)}\sigma_{n}(x_0^{r^{M}s_{n,s}^{M}})} -[-2\|\dot{\eta}(x_0)\|(r^{M}-r^{M}_{x_0})] \Big| \\
&=& \sqrt{2s/\pi}\Big| \frac{1/2-\mu_n(x_0^{rs_n})}{\sigma_n(x_0^{rs_n})} -[-2\|\dot{\eta}(x_0)\|(r-r_{x_0})] \Big| \\
&\le& \sqrt{2s/\pi}\epsilon^2(|r|+t_n/s_n) =  \epsilon^2(|r^{M}|+t_n/s_{n,s}^{M}).
\end{eqnarray*}
In addition, when $|r^{M}|\le \epsilon t_n/s_{n,s}^{M}$,
\begin{eqnarray*}
\Big| \Phi\Big(\frac{1/2-\mu_{n}(x_0^{r^{M}s_{n,s}^{M}})}{\sqrt{\pi/(2s)}\sigma_{n}(x_0^{r^{M}s_{n,s}^{M}})}\Big) -\Phi\big(-2\|\dot{\eta}(x_0)\|(r^{M}-r^{M}_{x_0})\big) \Big|\le 1
\end{eqnarray*}
and when $\epsilon t_n/s_{n,s}^{M} < |r^{M}| < \epsilon_n/s_{n,s}^{M}$,
\begin{eqnarray*}
&&\Big| \Phi\Big(\frac{1/2-\mu_{n}(x_0^{r^{M}s_{n,s}^{M}})}{\sqrt{\pi/(2s)}\sigma_{n}(x_0^{r^{M}s_{n,s}^{M}})}\Big) -\Phi\big(-2\|\dot{\eta}(x_0)\|(r^{M}-r^{M}_{x_0})\big) \Big| \\
&\le& \epsilon^2(|r^{M}|+t_n/s_{n,s}^{M})\phi(\|\dot{\eta}(x_0)\||r^{M}-r^{M}_{x_0}|),
\end{eqnarray*}
where $\phi$ is the density function of standard normal distribution.

Therefore, after substituting $ t=r^{M} s_{n,s}^{M}$, we have
\begin{align*}
&\int_{-\epsilon_{n,s}^{M}}^{\epsilon_{n,s}^{M}} |t|\|\dot{\psi}(x_0)\| \Big| \Phi\Big(\frac{1/2 - \mu_{n}(x_0^t)}{\sqrt{\pi/(2s)}\sigma_{n}(x_0^t)}\Big) - \Phi\Big(\frac{-2t\|\dot{\eta}(x_0)\|- 2a(x_0)t_n}{\sqrt{\pi/(2s)}s_n} \Big)\Big| dt \nonumber\\
=& \|\dot{\psi}(x_0)\| (s_{n,s}^{M})^2 \int_{-\epsilon_{n,s}^{M}/s_{n,s}^{M}}^{\epsilon_{n,s}^{M}/s_{n,s}^{M}} |r^{M}| \Big| \Phi\Big(\frac{1/2-\mu_{n}(x_0^{r^{M}s_{n,s}^{M}})}{\sqrt{\pi/(2s)}\sigma_{n}(x_0^{r^{M}s_{n,s}^{M}})}\Big) \nonumber \\
&\qquad\qquad\qquad\qquad\qquad
-\Phi\big(-2\|\dot{\eta}(x_0)\|(r^{M}-r^{M}_{x_0})\big)  \Big| dr^{M} \nonumber \\
\le & \|\dot{\psi}(x_0)\|(s_{n,s}^{M})^2 \Big[ \int_{|r^{M}|\le \epsilon t_n/s_{n,s}^{M}} |r^{M}| dr^{M} \nonumber\\
&
+ \epsilon^2 \int_{-\infty}^{\infty} |r^{M}|(|r^{M}|+t_n/s_{n,s}^{M})\phi(\|\dot{\eta}(x_0)\||r^{M}-r^{M}_{x_0}|) dr^{M}\Big]= o(\frac{s_n^2}{s}+t_n^2).  \nonumber 
\end{align*}
The inequality above leads to $R_{15}=o(s_n^2/s+t_n^2)$.

Combining \eqref{Taylor11_re_M}, \eqref{eq:decompose_11_re_M},  \eqref{eq:decompose_12_re_M}, \eqref{eq:decompose_13_re_M}, \eqref{eq:decompose_14_re_M} and \eqref{eq:decompose_15_re_M}, we have
\begin{eqnarray}
&&\int_{{\cal S}} \int_{-\epsilon_{n,s}^{M}}^{\epsilon_{n,s}^{M}} \psi(x_0^t) \big\{{\mathbb P}\big(S_{n,s}^{M}(x_0^t) < 1/2\big) - 
\indi{t<0}\big\} dt d\textrm{Vol}^{d-1}(x_0) \label{exp1_re_M}\\
\quad &=& \int_{{\cal S}} \int_{-\epsilon_{n,s}^{M}}^{\epsilon_{n,s}^{M}} t\|\dot{\psi}(x_0)\|  \big\{\Phi\Big(\frac{-2t\|\dot{\eta}(x_0)\|- 2a(x_0)t_n}{\sqrt{\pi/(2s)}s_n} \Big) \nonumber\\
&&\qquad\qquad\qquad\qquad
- \indi{t<0}\big\} dt d\textrm{Vol}^{d-1}(x_0) + o(\frac{s_n^2}{s}+t_n^2). \nonumber
\end{eqnarray}

Finally, after substituting $t=\sqrt{\pi/(2s)}us_n/2$ in \eqref{exp1_re_M}, we have, up to $o(s_n^2/s+t_n^2)$ difference,
\begin{align}
{\rm Regret}(\widehat{\phi}_{n,s}^{M})
=&\frac{\pi}{8s}s_n^2\int_{{\cal S}} \int_{-\infty}^{\infty} u\|\dot{\psi}(x_0)\|  \big\{\Phi\big(-\|\dot{\eta}(x_0)\|u-\frac{2a(x_0)t_n}{\sqrt{\pi/(2s)}s_n} \big) \nonumber\\
&\qquad\qquad\qquad\qquad\qquad
- \indi{u<0}\big\} du d\textrm{Vol}^{d-1}(x_0) \nonumber\\
=&\frac{\pi}{4s}s_n^2\int_{{\cal S}} \int_{-\infty}^{\infty} u\|\dot{\eta}(x_0)\| \bar{f}(x_0)  \big\{\Phi\big(-\|\dot{\eta}(x_0)\|  u-\frac{2a(x_0)t_n}{\sqrt{\pi/(2s)}s_n} \big)   \label{psi_eta_2_re_M}\\
&\qquad\qquad\qquad\qquad\qquad
 -\indi{u<0} \big\} du d\textrm{Vol}^{d-1}(x_0) \nonumber\\
=& B_1 \frac{\pi}{2s}s_n^2+B_2 t_n^2.  \label{II_re_M}
\end{align}
\eqref{psi_eta_2_re_M} holds by Lemma \ref{lemma:dot}, and \eqref{II_re_M} can be calculated by applying Lemma \ref{lemma:G}. This concludes the proof of Theorem \ref{thm:M-DNN_re}. \hfill $\blacksquare$

\subsection{Proof of Theorem \ref{thm:W-DNN_re}}\label{sec:pf_thm:W-DNN_re}
In this section, we apply similar notations as those in Section \ref{sec:pf_thm:M-DNN_re}. For the sake of simplicity, we omit $\bw_n$ in the subscript of such notations as $\widehat{\phi}_{n,s,\bw_n}^{W}$ and $S_{n,s,\bw_n}^{W}$. We have
\begin{eqnarray*}
\textrm{Regret} (\widehat{\phi}_{n,s}^{W})
&=& \int_{{\cal R}} \big[  {\mathbb P}\big(\widehat{\phi}_{n,s}^{W}(x) =0 \big) -  \indi{\eta(x)<1/2} \big]  d P^\circ(x).
\end{eqnarray*}

Denote the average of estimated regression function from $s$ subsamples as 
$$
S_{n,s}^{W}(x)=s^{-1}{\textstyle\sum}_{j=1}^s S_{n}^{(j)}(x).
$$
We can also write $S_{n,s}^{W}(x)$ as
\begin{eqnarray*}
S_{n,s}^{W}(x)=s^{-1}{\textstyle\sum}_{j=1}^s {\textstyle\sum}_{i=1}^n w_{ni}Y_{(i)}^{(j)} ={\textstyle\sum}_{l=1}^N  w_{Nl}Y_l,
\end{eqnarray*}
where 
\begin{align*}
\{Y_1,Y_2,\ldots Y_N\} =&\{Y_{(1)}^{(1)},Y_{(1)}^{(2)},\ldots,Y_{(1)}^{(s)},\ldots,Y_{(n)}^{(1)},Y_{(n)}^{(2)},\ldots,Y_{(n)}^{(s)} \},\\
\{w_{N1},w_{N2},\ldots w_{NN}\}
=&\{\frac{w_{n1}}{s},\frac{w_{n1}}{s}\ldots,\frac{w_{n1}}{s},\ldots,\frac{w_{nn}}{s},\frac{w_{nn}}{s},\ldots,\frac{w_{nn}}{s} \}.
\end{align*}
The W-DNN classifier is defined as
$$
\widehat{\phi}_{n,s}^{W}(x)=\indi{S_{n,s}^{W}(x)\ge1/2}.
$$
Since ${\mathbb P}\big(\widehat{\phi}_{n,s}^{W}(x)=0\big)= {\mathbb P}\big(S_{n,s}^{W}(x)<1/2\big)$, the regret of W-DNN becomes
\begin{align*}
{\rm Regret}(\widehat{\phi}_{n,s}^{W}) &= \int_{{\cal R}} \big\{ {\mathbb P}(S_{n,s}^{W}(x) < 1/2 ) - \indi{\eta(x)<1/2}\big\} dP^\circ(x). 
\end{align*}

Let $\mu_{n,s}(x)={\mathbb E}\{S_{n,s}^{W}(x)\}$, $\sigma_{n,s}^2(x)=\textrm{Var}\{S_{n,s}^{W}(x)\}$. We have 
\begin{align*}
\mu_{n,s}(x)&={\mathbb E}\{S_{n,s}^{W}(x)\} = {\mathbb E}\{ s^{-1}{\textstyle\sum}_{j=1}^s S_{n}^{(j)}(x)\} =\mu_n(x), \\
\sigma_{n,s}^2(x)&=\textrm{Var}\{S_{n,s}^{W}(x)\} = \textrm{Var}\{  s^{-1}{\textstyle\sum}_{j=1}^s S_{n}^{(j)}(x)\} =s^{-1}\sigma_n^2(x).
\end{align*}
Denote $\epsilon_{n,s}^{W}=\epsilon_{n}/\sqrt{s}$, $s_{n,s}^2=s_{n}^2/s$ and $t_{n,s}=t_{n}$. We have, uniformly for $\bw_n\in W_{n,\beta}$,
\begin{eqnarray*}
\sup_{x\in {\cal S}^{\epsilon_n}} |\mu_{n,s}(x) - \eta(x) - a(x)t_{n,s}|&=&\sup_{x\in {\cal S}^{\epsilon_n}} |\mu_{n}(x) - \eta(x) - a(x)t_n| \\ 
&=& o(t_n) = o(t_{n,s}),
\\
\sup_{x\in {\cal S}^{\epsilon_n}} \big|\sigma_{n,s}^2(x)-\frac{1}{4}s_{n,s}^2\big|&=&\sup_{x\in {\cal S}^{\epsilon_n}} \Big|\frac{\sigma_{n}^2(x)}{s}-\frac{1}{4}\frac{s_n^2}{s}\Big| \\
&=& o(s_n^2/s) = o(s_{n,s}^2). 
\end{eqnarray*}

We organize our proof in three steps. In {\it Step 1}, we decompose the integral over ${\cal R} \cap {\cal S}^{\epsilon_{n}}$ as an integral along ${\cal S}$ and an integral in the perpendicular direction; in {\it Step 2}, we focus on the complement set ${\cal R} \backslash {\cal S}^{\epsilon_n} $; {\it Step 3} combines the results and applies a normal approximation in ${\cal S}^{\epsilon_n}$ to yield the final conclusion. 

{\it Step 1}: Similarly to {\it Step 1} in Section \ref{sec:pf_thm:M-DNN_re}, we have 
\begin{align}
&\int_{{\cal R}\cap{{\cal S}^{\epsilon_n}}} \big\{ {\mathbb P}(S_{n,s}^{W}(x) < 1/2) - \indi{\eta(x)<1/2} \big\} d P^\circ(x) \nonumber\\  
=&\int_{{\cal S}} \int_{-\epsilon_n}^{\epsilon_n}\psi(x_0^t) \big\{{\mathbb P}\big(S_{n,s}^{W}(x_0^t) < 1/2\big) 
- \indi{t<0}\big\} dt d\textrm{Vol}^{d-1}(x_0)\{1+o(1)\}, \nonumber
\end{align}
uniformly for $\bw_n\in W_{n,\beta}$.

{\it Step 2}: Bound the contribution to regret from ${\cal R}\backslash {\cal S}^{\epsilon_n}$. We show that 
\begin{align*}
\sup_{\bw_n\in W_{n,\beta}} \int_{{\cal R}\backslash {\cal S}^{\epsilon_n}} \big\{ {\mathbb P}\big(S_{n,s}^{W}(x) < 1/2\big) 
- \indi{\eta(x)<1/2} \big\} d P^\circ(x) = o(\frac{s_n^2}{s}+t_n^2).
\end{align*}

\citet{S12} showed that, in any subsample, there exists a constant $c_{30}>0$ such that, for a sufficiently large $n$, 
\begin{eqnarray*}
\inf_{x\in {\cal R}\backslash {\cal S}^{\epsilon_n} } \big| \mu_n(x) -1/2 \big| \geq c_{30} \epsilon_n/4.  
\end{eqnarray*}
Applying Hoeffding's inequality to $S_{n,s}^{W}(x)$, we have
\begin{align*}
&|{\mathbb P}(S_{n,s}^{W}(x) < 1/2) - \indi{\eta(x)<1/2}| 
\leq \exp\Big(  \frac{-2(\mu_{n,s}(x)-1/2)^2}{\sum_{l=1}^N (w_{Nl}-0)^2} \Big) \\
=&\exp\Big(  \frac{-2(\mu_{n}(x)-1/2)^2}{s_n^2/s} \Big) \leq \exp\Big(  \frac{-2s(c_{30}\epsilon_n/4)^2}{n^{-\beta}} \Big) = o(\frac{s_n^2}{s}+t_n^2),
\end{align*}
uniformly for $\bw_n\in W_{n,\beta}$ and $x\in {\cal R}\backslash {\cal S}^{\epsilon_n}$.

{\it Step 3}: In the end, we will show
\begin{eqnarray*}
&&\int_{{\cal S}} \int_{-\epsilon_n}^{\epsilon_n} \psi(x_0^t) \big\{{\mathbb P}\big(S_{n,s}^{W}(x_0^t) < 1/2\big) - \indi{t<0}\big\} dt d\textrm{Vol}^{d-1}(x_0) \nonumber\\
&=& B_1 \frac{s_n^2}{s} + B_2 t_n^2 + o(\frac{s_n^2}{s}+t_n^2).
\end{eqnarray*}

According to \eqref{eq:psi_taylor}, we have
\begin{eqnarray}
&&\int_{{\cal S}} \int_{-\epsilon_n}^{\epsilon_n} \psi(x_0^t) \big\{{\mathbb P}\big(S_{n,s}^{W}(x_0^t) < 1/2\big)  - \indi{t<0}\big\} dt d\textrm{Vol}^{d-1}(x_0)\label{Taylor1_re_W} \\
&=& \int_{{\cal S}} \int_{-\epsilon_n}^{\epsilon_n} t\|\dot{\psi}(x_0)\|  \big\{{\mathbb P}\big(S_{n,s}^{W}(x_0^t) < 1/2\big)  \nonumber\\
&&\qquad\qquad\qquad\qquad\qquad
 - \indi{t<0}\big\} dt d\textrm{Vol}^{d-1}(x_0) \{1+o(1) \}. \nonumber
\end{eqnarray}

Next, we decompose
\begin{eqnarray}
&&\int_{{\cal S}} \int_{-\epsilon_n}^{\epsilon_n} t\|\dot{\psi}(x_0)\| \big\{{\mathbb P}\big(S_{n,s}^{W}(x_0^t) < 1/2\big)- \indi{t<0}\big\} dt d\textrm{Vol}^{d-1}(x_0) \label{eq:decompose_31_re_W} \\
&=& \int_{{\cal S}}\int_{-\epsilon_n}^{\epsilon_n}t\|\dot{\psi}(x_0)\| \big\{\Phi\big(\frac{1/2 - \mu_{n,s}(x_0^t)}{\sigma_{n,s}(x_0^t)} \big) \nonumber \\
&&\qquad\qquad\qquad\qquad\qquad
- \indi{t<0}\big\} dt d\textrm{Vol}^{d-1}(x_0) + R_{31}. \nonumber
\end{eqnarray}

Let $Z_l= (w_{Nl}Y_l - w_{Nl} \mathbb{E} [Y_l])/\sigma_{n,s}(x)$ and $V=\sum_{l=1}^N Z_l$. Note that $\mathbb{E}(Z_l)=0$, $\textrm{Var}(Z_l)<\infty$, and $\textrm{Var}(V)=1$. 
The nonuniform Berry-Esseen Theorem \citep{GP12} implies that there exists  a constant $c_{31}>0$, such that
$$
\Big|\mathbb{P}(V\le By) - \Phi(y)\Big| \le \frac{c_{31}A}{B^3(1 + |y|^3)},
$$
where $A=\sum_{l=1}^N E|Z_l|^3$ and $\big(\sum_{l=1}^N E|Z_l|^2)^{1/2}$. In the case of W-DNN, 
\begin{align*}
A&=\sum_{l=1}^N \mathbb{E}|\frac{w_{Nl}Y_l - w_{Nl} \mathbb{E} [Y_l]}{\sigma_{n,s}^3(x)}|^3 \leq \sum_{l=1}^N \frac{16|w_{Nl}|^3}{s_{n,s}^3}=\frac{16\sum_{l=1}^N w_{Nl}^3}{s_{n,s}^3}, \\
B&=({\textstyle\sum}_{l=1}^N \textrm{Var}(Z_l))^{1/2}=\sqrt{ \textrm{Var}(V)}=1.
\end{align*}
Denote $c_{32}=16c_{31}$, we have
\begin{eqnarray}
&&\sup_{x_0\in {\cal S}}\sup_{t\in [-\epsilon_n,\epsilon_n]}\Big|\mathbb{P}\Big(\frac{S_{n,s}^{W}(x_0^t)-\mu_{n,s}(x_0^t)}{\sigma_{n,s}(x_0^t)}\le y \Big)-\Phi(y)\Big| \label{BE_thm_re_W}\\
&\le& \frac{\sum_{l=1}^Nw_{Nl}^3}{s_{n,s}^3}\frac{c_{32}}{ 1 + |y|^3}. \nonumber
\end{eqnarray}

\cite{S12} showed that, there exists constants $c_{33},c_{34}>0$ such that, uniformly for $\bw_n\in W_{n,\beta}$,
\begin{align*}
\inf_{x_0 \in {\cal S}}\inf_{c_{33} t_n\leq|t|\leq \epsilon_n} \Big| \frac{1/2-\mu_{n}(x_0^t)}{\sigma_{n}(x_0^t)} \Big| \geq \frac{c_{34}|t|}{s_n}. 
\end{align*}
Hence, 
\begin{equation}
\inf_{x_0 \in {\cal S}}\inf_{c_{33} t_n\leq|t|\leq \epsilon_n} \Big| \frac{1/2-\mu_{n,s}(x_0^t)}{\sigma_{n,s}(x_0^t)} \Big| \geq \frac{c_{34}|t|}{s_n/\sqrt{s}}=\frac{c_{34}|t|}{s_{n,s}}. \label{BE_greater_re_W}
\end{equation}
Therefore,
\begin{align*}
& \int_{-\epsilon_n}^{\epsilon_n}  |t|\|\dot{\psi}(x_0)\| \Big| {\mathbb P}\big(S_{n,s}^{W}(x_0^t) < 1/2\big) - \Phi\Big(\frac{1/2 - \mu_{n,s}(x_0^t)}{\sigma_{n,s}(x_0^t)} \Big)\Big| dt \\
\leq &  \int_{|t| \leq c_{33} t_n} |t| \|\dot{\psi}(x_0)\| \frac{ c_{32}\sum_{l=1}^N w_{Nl}^3}{s_{n,s}^3} dt  \\
& +  \int_{c_{33} t_n \leq |t| \leq \epsilon_n} \frac{ c_{32}\sum_{l=1}^N w_{Nl}^3}{s_{n,s}^3}\frac{|t| \|\dot{\psi}(x_0)\|}{1+ c_{34}^3 |t|^3/s_{n,s}^3}  dt  \\
\leq &  \frac{ c_{32}\sum_{i=1}^n w_{ni}^3}{\sqrt{s}s_n^3} \int_{|t| \leq c_{33} t_n} |t|\|\dot{\psi}(x_0)\|  dt  \\
& + \frac{c_{32}\sum_{i=1}^n w_{ni}^3}{\sqrt{s}s_n^3} \int_{c_{33} t_n \leq |t| \leq \epsilon_n}  \frac{\|\dot{\psi}(x_0)\||t|}{ c_{34}^2 s |t|^2/s_n^2}  dt=o(\frac{s_n^2}{s}+t_n^2). 
\end{align*}
The inequality above leads to $|R_{31}|=o(s_n^2/s+t_n^2)$.

Next, we decompose 
\begin{align}
&\int_{{\cal S}}  \int_{-\epsilon_n}^{\epsilon_n} t\|\dot{\psi}(x_0)\| \big\{\Phi\big(\frac{1/2 - \mu_{n,s}(x_0^t)}{\sigma_{n,s}(x_0^t)} \big) - \indi{t<0}\big\} dt d\textrm{Vol}^{d-1}(x_0) \label{eq:decompose_32_re_W}\\
=& \int_{{\cal S}} \int_{-\epsilon_n}^{\epsilon_n} t\|\dot{\psi}(x_0)\| \big\{\Phi\big(\frac{-2t\|\dot{\eta}(x_0)\|- 2a(x_0)t_n}{s_n/\sqrt{s}} \big) \nonumber \\
&\qquad\qquad\qquad\qquad\qquad\qquad
- \indi{t<0}\big\} dt d\textrm{Vol}^{d-1}(x_0) + R_{32}. \nonumber
\end{align}
Denote $r^{W}=r\sqrt{s}$ and $r^{W}_{x_0}=r_{x_0}\sqrt{s}$. Similarly to bounding $R_{15}$ in \eqref{eq:decompose_15_re_M}, we have
\begin{align} 
&\int_{-\epsilon_n}^{\epsilon_n} |t|\|\dot{\psi}(x_0)\| \Big| \Phi\Big(\frac{1/2 - \mu_{n,s}(x_0^t)}{\sigma_{n,s}(x_0^t)}\Big) 
- \Phi\Big(\frac{-2t\|\dot{\eta}(x_0)\|- 2a(x_0)t_n}{s_n/\sqrt{s}} \Big)\Big| dt \nonumber\\
=& \|\dot{\psi}(x_0)\| s_{n,s}^2  \int_{-\epsilon_n/s_{n,s}}^{\epsilon_n/s_{n,s}} |r^{W}| \big|\Phi\Big(\frac{1/2-\mu_{n}(x_0^{r^{W}s_{n,s}})}{s^{-1/2}\sigma_{n}(x_0^{r^{W}s_{n,s}})}\Big) \nonumber \\
&\qquad\qquad\qquad\qquad\qquad
-\Phi\big(-2\|\dot{\eta}(x_0)\|(r^{W}-r^{W}_{x_0})\big)  \big| dr^{W}  \nonumber \\
\le & \| \dot{\psi}(x_0) \| s_{n,s}^2 \Big[ \int_{|r^{W}|\le \epsilon t_n/s_{n,s}} |r^{W}| dr^{W} \nonumber\\
&+ \epsilon^2 \int_{-\infty}^{\infty} |r^{W}|(|r^{W}|+t_n/s_{n,s})\phi(\|\dot{\eta}(x_0)\||r^{W}-r^{W}_{x_0}|) dr^{W} \Big] =o(\frac{s_n^2}{s}+t_n^2). \nonumber 
\end{align} 
The inequality above leads to $R_{32}=o(s_n^2/s+t_n^2)$.

By $(\ref{Taylor1_re_W})$, $(\ref{eq:decompose_31_re_W})$ and  \eqref{eq:decompose_32_re_W}, we have
\begin{eqnarray}
&&\int_{{\cal S}} \int_{-\epsilon_n}^{\epsilon_n} \psi(x_0^t) \big\{{\mathbb P}\big(S_{n,s}^{W}(x_0^t) < 1/2\big) - 
\indi{t<0}\big\} dt d\textrm{Vol}^{d-1}(x_0) \label{exp1_re_W} \\
\quad &=& \int_{{\cal S}} \int_{-\epsilon_n}^{\epsilon_n} t\|\dot{\psi}(x_0)\|  \big\{\Phi\big(\frac{-2t\|\dot{\eta}(x_0)\|- 2a(x_0)t_n}{s_n/\sqrt{s}} \big) \nonumber\\
&&\qquad\qquad\qquad\qquad
- \indi{t<0}\big\} dt d\textrm{Vol}^{d-1}(x_0) + o(s_n^2/s+t_n^2). \nonumber
\end{eqnarray}

Finally, after replacing $t=us_n/(2\sqrt{s})$ in \eqref{exp1_re_W}, we have, up to $o(s_n^2/s+t_n^2)$ difference,
\begin{align}
{\rm Regret}(\widehat{\phi}_{n,s}^{W})
&=\frac{s_n^2}{4s}\int_{{\cal S}} \int_{-\infty}^{\infty} \|\dot{\psi}(x_0)\| u \big\{\Phi\big(-\|\dot{\eta}(x_0)\|u-\frac{2a(x_0)t_n}{s_n/\sqrt{s}} \big) \nonumber\\
&\qquad\qquad\qquad\qquad\qquad
- \indi{u<0}\big\} du d\textrm{Vol}^{d-1}(x_0) \nonumber\\
&=\frac{s_n^2}{2s}\int_{{\cal S}} \int_{-\infty}^{\infty} \|\dot{\eta}(x_0)\| \bar{f}(x_0) u \big\{\Phi\big(-\|\dot{\eta}(x_0)\|  u-\frac{2a(x_0)t_n}{s_n/\sqrt{s}} \big)   \label{psi_eta_2_re_W}\\
&\qquad\qquad\qquad\qquad\qquad
 -\indi{u<0} \big\} du d\textrm{Vol}^{d-1}(x_0) \nonumber\\
&= B_1 \frac{1}{s}s_n^2+B_2 t_n^2.  \label{II_re_W}
\end{align}
\eqref{psi_eta_2_re_W} holds by Lemma \ref{lemma:dot}, and \eqref{II_re_W} can be calculated by Lemma \ref{lemma:G}. This completes the proof of  Theorem \ref{thm:W-DNN_re}.\hfill $\blacksquare$

\bibliographystyle{asa}
\bibliography{DNN}  

\newpage

\renewcommand{\theequation}{S.\arabic{equation}}
\renewcommand{\thetable}{S\arabic{table}}
\renewcommand{\thefigure}{S\arabic{figure}}
\renewcommand{\thelemma}{S.\arabic{lemma}}
\renewcommand{\thesubsection}{S.\Roman{subsection}}
\setcounter{equation}{0}
\setcounter{table}{0}
\setcounter{lemma}{0}
\setcounter{subsection}{0}

\invisiblesection{Supplementary}

\begin{center}
\Large\bf Supplementary Materials to: Distributed Nearest Neighbor Classification
\end{center}

{
\begin{center}
    Jiexin Duan, Xingye Qiao and Guang Cheng
\end{center}
}

The supplement is organized as follows:
\begin{itemize}
	\item In Section \ref{sec:pf_thm:DNN_rr}, we prove Theorem \ref{thm:M-DNN_rr} and \ref{thm:W-DNN_rr}.
    \item In Section \ref{sec:pf_thm:opt_M-DNN}, we prove Corollary \ref{thm:opt_M-DNN}.
	\item In Section \ref{sec:pf_thm:opt_W-DNN}, we prove Corollary \ref{thm:opt_W-DNN}.
	\item In Section \ref{sec:pf_thm:DNN_ci}, we prove Theorem \ref{thm:DNN_ci}.
	\item In Section \ref{sec:pf_thm:DNN_BNN}, we prove Corollary \ref{thm:DNN_BNN}.
	\item In Section \ref{sec:main_lemma}, we provide  lemmas.
\end{itemize}

\subsection{Proof of Theorem \ref{thm:M-DNN_rr} and \ref{thm:W-DNN_rr} }\label{sec:pf_thm:DNN_rr}
From Theorem \ref{thm:M-DNN_re} and Proposition \ref{thm:WNN_re}, we have, as $n,s\rightarrow \infty$, 
\begin{align}
\frac{{\rm Regret}(\widehat{\phi}_{n,s,\bw_{n}}^{M})}{{\rm Regret}(\widehat{\phi}_{N,\bw_{N}})} 
\rightarrow \frac{ B_1\frac{\pi}{2s}  \sum_{i=1}^n w_{ni}^2+ B_2 \big( \sum_{i=1}^n \frac{\alpha_i w_{ni}}{n^{2/d}}\big)^2 }{  B_1\sum_{i=1}^N w_{Ni}^2 + B_2 \big(\sum_{i=1}^N \frac{\alpha_i w_{Ni}}{N^{2/d}} \big)^2 } = (\frac{\pi}{2})^{\frac{4}{d+4}}. \nonumber
\end{align}
The last equality holds by \eqref{eq:M-DNN_rr_weight1} and \eqref{eq:M-DNN_rr_weight2}. This completes the proof of Theorem \ref{thm:M-DNN_rr}.

Similarly, from Theorem \ref{thm:W-DNN_re} and Proposition \ref{thm:WNN_re}, we have, as $n,s\rightarrow \infty$,
\begin{align}
\frac{{\rm Regret}(\widehat{\phi}_{n,s,\bw_{n}}^{W})}{{\rm Regret}(\widehat{\phi}_{N,\bw_{N}})} 
\rightarrow \frac{ B_1 \sum_{i=1}^n w_{ni}^2+ B_2 \big( \sum_{i=1}^n \frac{\alpha_i w_{ni}}{n^{2/d}}\big)^2 }{  B_1\sum_{i=1}^N w_{Ni}^2 + B_2 \big(\sum_{i=1}^N \frac{\alpha_i w_{Ni}}{N^{2/d}} \big)^2 } = 1. \nonumber
\end{align}
The last equality holds by \eqref{eq:W-DNN_rr_weight1} and \eqref{eq:W-DNN_rr_weight2}. This completes the proof of Theorem \ref{thm:W-DNN_rr}. \hfill $\blacksquare$

\subsection{Proof of Corollary \ref{thm:opt_M-DNN}}
\label{sec:pf_thm:opt_M-DNN} 
Denote $a\succeq b$ if $b=O(a)$, $a\succ b$ if $b=o(a)$, $a\asymp b$ if $a\succeq b$ and $b \succeq a$. To find the optimal value of \eqref{eq:M-DNN_re}, we write its Lagrangian as 
$$
L(\bw_n)=\big(\sum_{i=1}^n \frac{\alpha_i w_{ni}}{n^{2/d}}\big)^2 + \lambda \sum_{i=1}^n w_{ni}^2 + \nu(\sum_{i=1}^{n}w_{ni}-1),
$$
where $\lambda=(\pi B_1)/(2sB_2)$. Since all the weights are nonnegative, we denote $l^{*}=\max\{i:w_{ni}^{*}>0\}$. Setting the derivative of $L(\bw_n)$ to be $0$, we have \begin{equation}
\frac{\partial L(\bw_n)}{\partial w_{ni}} = 2 n^{-4/d}\alpha_i\sum_{i=1}^{l^{*}} \alpha_i w_{ni} + 2\lambda w_{ni} +  \nu =0. \label{derivative_M}
\end{equation}
(i) Summing $(\ref{derivative_M})$ from 1 to $l^{*}$, (ii) multiplying $(\ref{derivative_M})$ by $\alpha_i$ and then summing from 1 to $l^{*}$, (ii) multiplying $(\ref{derivative_M})$ by $w_{ni}$ and then summing from 1 to $l^{*}$, we have
\begin{eqnarray*}
2 n^{-4/d}(l^{*})^{1+2/d}{\textstyle\sum}_{i=1}^{l^{*}}\alpha_i w_{ni}  + 2\lambda + \nu l^{*} &=& 0,\\
2 n^{-4/d}{\textstyle\sum}_{i=1}^{l^{*}}\alpha_i w_{ni} {\textstyle\sum}_{i=1}^{l^{*}} \alpha_i^2 + 2\lambda {\textstyle\sum}_{i=1}^{l^{*}}\alpha_i w_{ni} + \nu (l^{*})^{1+2/d} &=& 0, \\
 2 n^{-4/d}\big({\textstyle\sum}_{i=1}^{l^{*}} \alpha_i w_{ni}\big)^2 + 2\lambda {\textstyle\sum}_{i=1}^{l^{*}} w_{ni}^2 +  \nu &=& 0.
\end{eqnarray*}
Therefore, we have
\begin{align}
&\sum_{i=1}^{l^{*}} \alpha_i w_{ni} \asymp (l^{*})^{2/d},\;\; \sum_{i=1}^{l^{*}} w_{ni}^2 \asymp \frac{1}{l^{*}} \;\;{\rm  and}\;\; \label{l_order}\\
&w_{ni}^{*} = \frac{1}{l^{*}} + \frac{(l^{*})^{4/d}-(l^{*})^{2/d}\alpha_i}{\sum_{j=1}^{l^{*}} \alpha_j^2+\lambda n^{4/d} -(l^{*})^{1+4/d}}. \label{optweight_M}
\end{align}
Here $w_{ni}^{*}$ is decreasing in $i$, since $\alpha_i$ is increasing in $i$ and $\sum_{j=1}^{l^{*}} \alpha_j^2+\lambda n^{4/d} -(l^{*})^{1+4/d}>0$ from Lemma \ref{alpha}. Next we solve for $l^{*}$. According to the definition of $l^{*}$, we only need to find the last $l$ such that $w_{nl}^{*}> 0$. Using the results from Lemma \ref{alpha}, solving this equation reduces to finding the $l^{*}$ such that
\begin{align*}
(1+\frac{2}{d})(l^{*}-1)^{2/d} \le \lambda n^{4/d}(l^{*})^{-1-2/d} + \frac{(d+2)^2}{d(d+4)}(l^{*})^{2/d}\{1+O(\frac{1}{l^{*}})\} \qquad\\ 
\qquad\qquad\qquad\qquad\qquad\qquad\qquad\qquad \le (1+\frac{2}{d})(l^{*})^{2/d}.
\end{align*}
For large $n,s$, we have
\begin{equation}
l^{*}= \Big\lceil \Big\{\frac{d(d+4)}{2(d+2)}\Big\}^{\frac{d}{d+4}}\lambda^{\frac{d}{d+4}}n^{\frac{4}{d+4}}\Big\rceil=\Big\lceil \Big\{\frac{d(d+4)}{2(d+2)}\Big\}^{\frac{d}{d+4}}\Big(\frac{\pi B_1}{2sB_2}\Big)^{\frac{d}{d+4}}n^{\frac{4}{d+4}}\Big\rceil.\nonumber
\end{equation}
Due to Assumption (w.1) in Section \ref{sec:defwnb}, we have $l^{*} \rightarrow \infty$ as $n \rightarrow \infty$. When $\gamma<2/(d+4)$, plugging $l^{*}$ and $(\ref{sum_alpha})$ into $(\ref{optweight_M})$ yields the optimal weight and \eqref{eq:opt_M-DNN}. 

Denote $H(\bw_n)$ as the Hessian matrix of $L(\bw_n)$. We have 
\begin{align*}
\frac{\partial^2 L(\bw_n)}{\partial w_{ni}^2} = 2 n^{-4/d}\alpha_i^2 + 2\lambda \;\;{\rm and }\;\; \frac{\partial^2 L(\bw_n)}{\partial w_{ni}\partial w_{nj}} = 2 n^{-4/d}\alpha_i\alpha_j. 
\end{align*}
For any nonzero vector $X_{l^{*}}=(x_1,...,x_{l^{*}})^T$, we have
\begin{align*}
X_{l^{*}}^T H(\bw_n) X_{l^{*}}=& 2n^{-4/d}\sum_{i=1}^{l^{*}} \alpha_i^2x_i^2 + 2\lambda \sum_{i=1}^{l^{*}} x_i^2 + 2 n^{-4/d} \sum_{i\neq j} \alpha_i\alpha_j x_ix_j \\
=&  2n^{-4/d} \big(\sum_{i=1}^{l^{*}} \alpha_ix_i \big)^2 +  2\lambda \sum_{i=1}^{l^{*}} x_i^2 >0.
\end{align*}
Therefore, $H(\bw_n)$ is positive definite, and this verifies that the above optimal value achieves the global minimum.

Next, we analyze the case of $\gamma\ge 2/(d+4)$. By Cauchy--Schwarz inequality, we have
\begin{align*}
    ({\textstyle\sum}_{i=1}^{l^{*}} \bw_{ni}^3) ({\textstyle\sum}_{i=1}^{l^{*}} \bw_{ni}) \ge& ({\textstyle\sum}_{i=1}^{l^{*}} \bw_{ni}^{3/2}\bw_{ni}^{1/2})^2 \\
    =& ({\textstyle\sum}_{i=1}^{l^{*}} \bw_{ni}^2)^2\ge({\textstyle\sum}_{i=1}^{l^{*}} \bw_{ni}^2)^{3/2}/\sqrt{{l^{*}}}.
\end{align*}
The above inequality, along with condition \eqref{eq:extra_condition}, suggests that ${l^{*}}\succ s$. As $\gamma \ge 2/(d+4)$, we have ${l^{*}} \succ s \succeq N^{2/(d+4)}$ and $n=O(N^{(d+2)/(d+4)})$. Applying \eqref{l_order}, we have, as $n,s\rightarrow \infty$, 
\begin{align*}
\big(\sum_{i=1}^n \frac{\alpha_i w_{ni}}{n^{2/d}}\big)^2  \asymp  ({l^{*}}/n)^{4/d} \succ N^{-4/(d+4)}. 
\end{align*}
\citet{S12} showed that 
\begin{align}
{\rm Regret }(\widehat{\phi}_{N,\bw_N^*}) \asymp N^{-4/(d+4)}. \label{OWNN_re}
\end{align}
Therefore, we have, as $n,s\rightarrow \infty$, 
\begin{align*}
\frac{{\rm Regret}(\widehat{\phi}_{n,s,w_n}^{M})}{{\rm Regret}(\widehat{\phi}_{N,\bw_N^*})} \asymp& \frac{ B_1\frac{\pi}{2s}  \sum_{i=1}^n w_{ni}^2 + B_2 \big(\sum_{i=1}^n \frac{\alpha_i w_{ni}}{n^{2/d}}\big)^2 }{N^{-4/(d+4)}} \\
\succeq &\frac{B_2 \big(\sum_{i=1}^n \frac{\alpha_i w_{ni}}{n^{2/d}}\big)^2}{N^{-4/(d+4)}} \rightarrow \infty. \nonumber  
\end{align*} 
This completes the proof of Corollary \ref{thm:opt_M-DNN}. \hfill $\blacksquare$

\subsection{Proof of Corollary \ref{thm:opt_W-DNN}}
\label{sec:pf_thm:opt_W-DNN}
To find the optimal value of \eqref{eq:W-DNN_re}, we write its Lagrangian as 
\begin{align*}
L(\bw_n)=\big(\sum_{i=1}^n \frac{\alpha_i w_{ni}}{n^{2/d}}\big)^2 + \delta \sum_{i=1}^n w_{ni}^2 + \nu(\sum_{i=1}^{n}w_{ni}-1), 
\end{align*}
where $\delta=(B_1)/(sB_2)$.

Similar to Section \ref{sec:pf_thm:opt_M-DNN}, replacing $l^{*}$ by $l^{\dag}$ in the optimization, we have 
\begin{equation}
w_{ni}^{\dag} = \frac{1}{l^{\dag}} + \frac{(l^{\dag})^{4/d}-(l^{\dag})^{2/d}\alpha_i}{\sum_{j=1}^{l^{\dag}} \alpha_j^2+\delta n^{4/d} -(l^{\dag})^{1+4/d}}. \label{optweight_W}
\end{equation}
For large $n,s$, we have
\begin{equation}
l^{\dag}= \Big\lceil \Big\{\frac{d(d+4)}{2(d+2)}\Big\}^{\frac{d}{d+4}}\delta^{\frac{d}{d+4}}n^{\frac{4}{d+4}}\Big\rceil=\Big\lceil \Big\{\frac{d(d+4)}{2(d+2)}\Big\}^{\frac{d}{d+4}}\Big(\frac{ B_1}{sB_2}\Big)^{\frac{d}{d+4}}n^{\frac{4}{d+4}}\Big\rceil.\nonumber
\end{equation}
Due to Assumption (w.1) in Section \ref{sec:defwnb}, we have $l^{\dag} \rightarrow \infty$ as $n \rightarrow \infty$. When $\gamma<  4/(d+4)$, plugging $l^{\dag}$ and $(\ref{sum_alpha})$ into $(\ref{optweight_W})$ yields the optimal weight and  \eqref{eq:opt_W-DNN}. 

When $\gamma\ge 4/(d+4)$, we have $s \succeq N^{4/(d+4)}$ and $n=O(N^{d/(d+4)})$. Similary to  \eqref{l_order}, we have, as $n,s\rightarrow \infty$, 
\begin{align*}
\big(\sum_{i=1}^n \frac{\alpha_i w_{ni}}{n^{2/d}}\big)^2  \asymp  \Big(\frac{{l^{\dag}}}{n}\Big)^{4/d} \succeq  \big(l^{\dag}\big)^{4/d}N^{-4/(d+4)} \succ N^{-4/(d+4)}. 
\end{align*}
The last inequality holds by $l^{\dag} \rightarrow \infty$ as $n \rightarrow \infty$. Therefore, along with \eqref{OWNN_re}, we have, as $n,s\rightarrow \infty$, 
\begin{align*}
\frac{{\rm Regret}(\widehat{\phi}_{n,s,w_n}^{W})}{{\rm Regret}(\widehat{\phi}_{N,\bw_N^*})} \asymp& \frac{ B_1\frac{\pi}{2s}  \sum_{i=1}^n w_{ni}^2 + B_2 \big(\sum_{i=1}^n \frac{\alpha_i w_{ni}}{n^{2/d}}\big)^2 }{N^{-4/(d+4)}} \\
\succeq &\frac{B_2 \big(\sum_{i=1}^n \frac{\alpha_i w_{ni}}{n^{2/d}}\big)^2}{N^{-4/(d+4)}} \rightarrow \infty. \nonumber  
\end{align*} 
This completes the proof of Corollary \ref{thm:opt_W-DNN}. \hfill $\blacksquare$

\subsection{Proof of Theorem \ref{thm:DNN_ci}}\label{sec:pf_thm:DNN_ci}
We will prove \eqref{eq:M-DNN_ci} and \eqref{eq:W-DNN_ci} in {\it Part 1} and {\it Part 2} of this section respectively. For the sake of simplicity, we omit $\bw_n$ in the subscript of such notations as $\widehat{\phi}_{n,s,\bw_n}^{M}$ and $S_{n,s,\bw_n}^{M}$.

{\it Part 1}: We use similar notations as those in Section \ref{sec:pf_thm:M-DNN_re}. Denote $\widehat{\phi}_{{\cal D}_1}^{M}$ and $\widehat{\phi}_{{\cal D}_2}^{M}$ as $\widehat{\phi}_{n,s}^{M}$ based on ${\cal D}_1$ and ${\cal D}_2$, which are i.i.d. copies of ${\cal D}$. Write $\bar{P}(x)=\pi_1 P_1 + (1-\pi_1) P_0$. We have 
\begin{align}
{\rm CIS}(\widehat{\phi}_{{\cal D}_2}^{M}) 
&= {\mathbb E}_{X}\big[{\mathbb P}_{{\cal D}_1,{\cal D}_2}\big(\widehat{\phi}_{{\cal D}_1}^{M}(X) \ne \widehat{\phi}_{{\cal D}_2}^{M}(X)\big|X \big)\big] \label{eq:CIS_proof_M} \\
&= {\mathbb E}_{X} \big[ {\mathbb P}_{{\cal D}_1,{\cal D}_2}\big(\widehat{\phi}_{{\cal D}_1}^{M}(X)=1, \widehat{\phi}_{{\cal D}_2}^{M}(X) =0 \big|X \big)\big] \nonumber\\
&\qquad\qquad 
 + {\mathbb E}_{X} \big[ {\mathbb P}_{{\cal D}_1,{\cal D}_2}\big(\widehat{\phi}_{{\cal D}_1}^{M}(X)=0, \widehat{\phi}_{{\cal D}_2}^{M}(X) =1 \big|X\big)\big] \nonumber\\
&= 2 {\mathbb E}_{X} \big[  {\mathbb P}_{{\cal D}_1}\big(\widehat{\phi}_{{\cal D}_1}^{M}(X)=0|X\big) \big(1- {\mathbb P}_{{\cal D}_1}\big(\widehat{\phi}_{{\cal D}_1}^{M}(X)=0|X\big)\big) \big]\nonumber \\
&= 2\int_{{\cal R}} {\mathbb P}(S_{n,s}^{M}(x) < 1/2 ) \big(1- {\mathbb P}(S_{n,s}^{M}(x) < 1/2)\big)d \bar{P}(x)\nonumber\\
&=  2\int_{{\cal R}} \big\{ {\mathbb P}(S_{n,s}^{M}(x) < 1/2) - \indi{\eta(x)<1/2}\big\} d\bar{P}(x) \nonumber\\
&\qquad
- \int_{{\cal R}} \big\{ {\mathbb P}^2(S_{n,s}^{M}(x) < 1/2) - \indi{\eta(x)<1/2} \big\} d \bar{P}(x). \nonumber
\end{align}

Next, we organize our proof in four steps similar to Section \ref{sec:pf_thm:M-DNN_re}. 
 
{\it Step 1}: Similar to {\it Step 1} of Section \ref{sec:pf_thm:M-DNN_re}, we have, uniformly for $\bw_n\in W_{n,\beta}$, 
\begin{align*}
&\int_{{\cal R}\cap{{\cal S}^{\epsilon_n}}} \big\{ {\mathbb P}(s_{n,s}^{M}(x) < 1/2) - \indi{\eta(x)<1/2} \big\} d\bar{P}(x)  \nonumber\\
=&\int_{{\cal S}} \int_{-\epsilon_n}^{\epsilon_n}\bar{f}(x_0^t) \big\{{\mathbb P}\big(S_{n,s}^{M}(x_0^t) < 1/2\big) - \indi{t<0}\big\} dt d\textrm{Vol}^{d-1}(x_0)\{1+o(1)\}, \nonumber\\
&\int_{{\cal R}\cap{{\cal S}^{\epsilon_n}}} \big\{ {\mathbb P}^2(s_{n,s}^{M}(x) < 1/2) - \indi{\eta(x)<1/2} \big\} d\bar{P}(x)  \nonumber\\
=&\int_{{\cal S}} \int_{-\epsilon_n}^{\epsilon_n}\bar{f}(x_0^t) \big\{{\mathbb P}^2\big(S_{n,s}^{M}(x_0^t) < 1/2\big) - \indi{t<0}\big\} dt d\textrm{Vol}^{d-1}(x_0)\{1+o(1)\}. \nonumber
\end{align*}

{\it Step 2}: Bound the contribution to CIS from ${\cal R}\backslash {\cal S}^{\epsilon_{n}}$. We have,
\begin{eqnarray}
&&\sup_{\bw_n\in W_{n,\beta}} \int_{{\cal R}\backslash {\cal S}^{\epsilon_n}} \big\{ {\mathbb P}\big(S_{n,s}^{M}(x) < 1/2\big)   \label{bound21_ci_M} \\
&&\qquad\qquad\qquad
- \indi{\eta(x)<1/2} \big\} d \bar{P}(x) = o(s_n^2/s+t_n^2), \nonumber\\
&&\sup_{\bw_n\in W_{n,\beta}} \int_{{\cal R}\backslash {\cal S}^{\epsilon_n}} \big\{ {\mathbb P}^2\big(S_{n,s}^{M}(x) < 1/2\big)  \label{bound22_ci_M}\\
&&\qquad\qquad\qquad
- \indi{\eta(x)<1/2} \big\} d \bar{P}(x) = o(s_n^2/s+t_n^2).  \nonumber
\end{eqnarray}
\eqref{bound21_ci_M} holds if we replace $P^\circ(x)$ by $\bar{P}(x)$ in {\it Step 2} of Section \ref{sec:pf_thm:M-DNN_re}. Furthermore, \eqref{bound22_ci_M} holds since 
\begin{eqnarray*}
&&\big|{\mathbb P}^2\big(S_{n,s}^{M}(x) < 1/2\big) - \indi{\eta(x)<1/2}\big| \\ 
&\le& 2\big|{\mathbb P}\big(S_{n,s}^{M}(x) < 1/2\big) - \indi{\eta(x)<1/2}\big|.
\end{eqnarray*}

{\it Step 3}: Bound the contribution to CIS from ${\cal S}^{\epsilon_{n}}\backslash{\cal S}^{\epsilon_{n,s}^{M}}$. We have,
\begin{align}
&\sup_{\bw_n\in W_{n,\beta}}\int_{{\cal S}} \int_{(-\epsilon_n,\epsilon_n)\backslash(-\epsilon_{n,s}^{M}, \epsilon_{n,s}^{M})}\bar{f}(x_0^t)  \big\{{\mathbb P}\big(S_{n,s}^{M}(x_0^t) < 1/2\big) 
\label{bound31_ci_M}\\
&\qquad\qquad\qquad\qquad
- \indi{t<0}\big\} dt d\textrm{Vol}^{d-1}(x_0)
= o(s_n^2/s+t_n^2), \nonumber \\
&\sup_{\bw_n\in W_{n,\beta}}\int_{{\cal S}} \int_{(-\epsilon_n,\epsilon_n)\backslash(-\epsilon_{n,s}^{M}, \epsilon_{n,s}^{M})}\bar{f}(x_0^t)  \big\{{\mathbb P}^2\big(S_{n,s}^{M}(x_0^t) < 1/2\big) 
\label{bound32_ci_M} \\
&\qquad\qquad\qquad\qquad
- \indi{t<0}\big\} dt d\textrm{Vol}^{d-1}(x_0)
= o(s_n^2/s+t_n^2). \nonumber
\end{align}

\eqref{bound31_ci_M} holds if we replace $\psi(x_0^t)$ by $\bar{f}(x_0^t)$ in {\it Step 3} of Section \ref{sec:pf_thm:M-DNN_re}. Furthermore, \eqref{bound32_ci_M} holds since 
\begin{eqnarray*}
\big|{\mathbb P}^2\big(S_{n,s}^{M}(x_0^t) < 1/2\big) - \indi{t<0}\big| 
\le 2\big|{\mathbb P}\big(S_{n,s}^{M}(x_0^t) < 1/2\big) - \indi{t<0}\big|.
\end{eqnarray*}

{\it Step 4}: In the end, we will show
\begin{align*}
&\int_{{\cal S}} \int_{-\epsilon_{n,s}^{M}}^{\epsilon_{n,s}^{M}} {\bar f}(x_0^t) \big\{{\mathbb P}\big(S_{n,s}^{M}(x_0^t) < 1/2\big) - \indi{t<0}\big\} dt d\textrm{Vol}^{d-1}(x_0) \nonumber\\
&-\int_{{\cal S}} \int_{-\epsilon_{n,s}^{M}}^{\epsilon_{n,s}^{M}} {\bar f}(x_0^t) \big\{{\mathbb P}^2\big(S_{n,s}^{M}(x_0^t) < 1/2\big) - \indi{t<0}\big\} dt d\textrm{Vol}^{d-1}(x_0) \nonumber\\
=& \frac{1}{2}B_3 \sqrt{\frac{\pi}{2s}}s_n +  o(\frac{s_n}{\sqrt{s}}+t_n). \nonumber
\end{align*}

Taylor expansion leads to
$$
\bar{f}(x_0^t) = \bar{f}(x_0) + \dot{\bar{f}}(x_0)^T\frac{\dot{\eta}(x_0)}{\|\dot{\eta}(x_0)\|} t + o(t).
$$
Hence, 
\begin{eqnarray}
&&\int_{{\cal S}} \int_{-\epsilon_{n,s}^{M}}^{\epsilon_{n,s}^{M}} {\bar f}(x_0^t) \big\{{\mathbb P}\big(S_{n,s}^{M}(x_0^t) < 1/2\big) - \indi{t<0}\big\} dt d\textrm{Vol}^{d-1}(x_0) \label{Taylor1_ci_M}\\
\quad &=& \big[ \int_{{\cal S}} \int_{-\epsilon_{n,s}^{M}}^{\epsilon_{n,s}^{M}} \bar{f}(x_0) \big\{{\mathbb P}\big(S_{n,s}^{M}(x_0^t) < 1/2\big) \nonumber - \indi{t<0}\big\} dt d\textrm{Vol}^{d-1}(x_0) \nonumber\\
&&+\int_{{\cal S}} \int_{-\epsilon_{n,s}^{M}}^{\epsilon_{n,s}^{M}} \frac{\dot{\bar{f}}(x_0)^T\dot{\eta}(x_0)t}{\|\dot{\eta}(x_0)\|} \big\{{\mathbb P}\big(S_{n,s}^{M}(x_0^t) < 1/2\big) \nonumber\\
&&\qquad\qquad\qquad\qquad\qquad
- \indi{t<0}\big\} dt d\textrm{Vol}^{d-1}(x_0) \big] \{1+o(1) \}. \nonumber
\end{eqnarray}


Next, we decompose
\begin{align}
&\int_{{\cal S}} \int_{-\epsilon_{n,s}^{M}}^{\epsilon_{n,s}^{M}}\bar{f}(x_0) \big\{ {\mathbb P}\big(S_{n,s}^{M}(x_0^t) < 1/2\big) - \indi{t<0} \big\} dt d\textrm{Vol}^{d-1}(x_0)\label{eq:decompose_510_ci_M} \\
=& \int_{{\cal S}} \int_{-\epsilon_{n,s}^{M}}^{\epsilon_{n,s}^{M}} \bar{f}(x_0) \big\{ \Phi\Big[\frac{\sqrt{s}\big(1/2-{\mathbb P}\big(S_{n}(x_0^t)\geq1/2\big)\big)}{\sqrt{{\mathbb P}\big(S_{n}(x_0^t)<1/2\big){\mathbb P}\big(S_{n}(x_0^t) \ge 1/2\big)}} \Big] \nonumber \\
& \qquad\qquad\qquad\qquad\qquad
- \indi{t<0}  \big\} dt d\textrm{Vol}^{d-1}(x_0) +R_{510}, \nonumber\\
&\int_{{\cal S}} \int_{-\epsilon_{n,s}^{M}}^{\epsilon_{n,s}^{M}}\frac{\dot{\bar{f}}(x_0)^T\dot{\eta}(x_0)t}{\|\dot{\eta}(x_0)\|}
 \big\{ {\mathbb P}\big(S_{n,s}^{M}(x_0^t) < 1/2\big)\label{eq:decompose_511_ci_M}\\
& \qquad\qquad\qquad\qquad\qquad
 - \indi{t<0} \big\} dt d\textrm{Vol}^{d-1}(x_0)  \nonumber \\
=& \int_{{\cal S}} \int_{-\epsilon_{n,s}^{M}}^{\epsilon_{n,s}^{M}} \frac{\dot{\bar{f}}(x_0)^T\dot{\eta}(x_0)t}{\|\dot{\eta}(x_0)\|} \big\{ \Phi\Big[\frac{\sqrt{s}\big(1/2-{\mathbb P}\big(S_{n}(x_0^t)\geq1/2\big)\big)}{\sqrt{{\mathbb P}\big(S_{n}(x_0^t)<1/2\big){\mathbb P}\big(S_{n}(x_0^t) \ge 1/2\big)}} \Big] \nonumber \\
& \qquad\qquad\qquad\qquad\qquad
- \indi{t<0}  \big\} dt  d\textrm{Vol}^{d-1}(x_0) +R_{511}. \nonumber
\end{align}

Applying \eqref{Edgeworth}, we have 
\begin{align*}
|R_{510}|\le&\int_{{\cal S}} \int_{-\epsilon_{n,s}^{M}}^{\epsilon_{n,s}^{M}}\bar{f}(x_0) \Big| {\mathbb P}\big(S_{n,s}^{M}(x_0^t) < 1/2\big)\\
&
-\Phi\Big[\frac{\sqrt{s}\big(1/2-{\mathbb P}\big(S_{n}(x_0^t)\geq1/2\big)\big)}{\sqrt{{\mathbb P}\big(S_{n}(x_0^t)<1/2\big){\mathbb P}\big(S_{n}(x_0^t) \ge 1/2\big)}} \Big] \Big|  dt d\textrm{Vol}^{d-1}(x_0) \nonumber\\
\leq &   O(\frac{1}{\sqrt{s}}) \int_{{\cal S}} \int_{-\epsilon_{n,s}^{M}}^{\epsilon_{n,s}^{M}}\bar{f}(x_0)  dt d\textrm{Vol}^{d-1}(x_0) = o(t_n+\frac{1}{\sqrt{s}}s_n),
\end{align*}
\begin{align*}
|R_{511}|\le&\int_{{\cal S}} \int_{-\epsilon_{n,s}^{M}}^{\epsilon_{n,s}^{M}}\frac{|\dot{\bar{f}}(x_0)^T\dot{\eta}(x_0)|}{\|\dot{\eta}(x_0)\|}|t| \Big| {\mathbb P}\big(S_{n,s}^{M}(x_0^t) < 1/2\big)\\
&
-\Phi\Big[\frac{\sqrt{s}\big(1/2-{\mathbb P}\big(S_{n}(x_0^t)\geq1/2\big)\big)}{\sqrt{{\mathbb P}\big(S_{n}(x_0^t)<1/2\big){\mathbb P}\big(S_{n}(x_0^t) \ge 1/2\big)}} \Big] \Big|  dt d\textrm{Vol}^{d-1}(x_0) \nonumber\\
\leq & O(\frac{1}{\sqrt{s}}) \int_{{\cal S}} \int_{-\epsilon_{n,s}^{M}}^{\epsilon_{n,s}^{M}}\frac{|\dot{\bar{f}}(x_0)^T\dot{\eta}(x_0)|}{\|\dot{\eta}(x_0)\|}|t|  dt d\textrm{Vol}^{d-1}(x_0) = o(t_n^2+\frac{1}{s}s_n^2).
\end{align*}

Next, we decompose
\begin{align}
&\int_{{\cal S}}\int_{-\epsilon_{n,s}^{M}}^{\epsilon_{n,s}^{M}}\bar{f}(x_0) \big\{ \Phi\big[\frac{\sqrt{s}\big(1/2-{\mathbb P}\big(S_{n}(x_0^t)\ge1/2\big)\big)}{\sqrt{{\mathbb P}\big(S_{n}(x_0^t)<1/2\big){\mathbb P}\big(S_{n}(x_0^t) \ge 1/2\big)}} \big] \label{eq:decompose_512_ci_M}\\
&\qquad\qquad\qquad\qquad
- \indi{t<0} \big\} dt d\textrm{Vol}^{d-1}(x_0)\nonumber\\
=& \int_{{\cal S}} \int_{-\epsilon_{n,s}^{M}}^{\epsilon_{n,s}^{M}}\bar{f}(x_0)\big \{\Phi\big[2\sqrt{s}\big({\mathbb P}\big(S_{n}(x_0^t)<1/2\big)-1/2\big) \big]\nonumber  \\
&\qquad\qquad\qquad\qquad
-\indi{t<0} \big\} dt d\textrm{Vol}^{d-1}(x_0) +R_{512}, \nonumber
\end{align}
\begin{align}
&\int_{{\cal S}}\int_{-\epsilon_{n,s}^{M}}^{\epsilon_{n,s}^{M}}\frac{\dot{\bar{f}}(x_0)^T\dot{\eta}(x_0)t}{\|\dot{\eta}(x_0)\|} \big\{ \Phi\big[\frac{\sqrt{s}\big(1/2-{\mathbb P}\big(S_{n}(x_0^t)\ge1/2\big)\big)}{\sqrt{{\mathbb P}\big(S_{n}(x_0^t)<1/2\big){\mathbb P}\big(S_{n}(x_0^t) \ge 1/2\big)}} \big] \label{eq:decompose_513_ci_M}\\
&\qquad\qquad\qquad\qquad
- \indi{t<0} \big\} dt d\textrm{Vol}^{d-1}(x_0)\nonumber\\
=& \int_{{\cal S}} \int_{-\epsilon_{n,s}^{M}}^{\epsilon_{n,s}^{M}}\frac{\dot{\bar{f}}(x_0)^T\dot{\eta}(x_0)t}{\|\dot{\eta}(x_0)\|}\big \{\Phi\big[2\sqrt{s}\big({\mathbb P}\big(S_{n}(x_0^t)<1/2\big)-1/2\big) \big]\nonumber  \\
&\qquad\qquad\qquad\qquad
-\indi{t<0} \big\} dt d\textrm{Vol}^{d-1}(x_0) +R_{513}. \nonumber
\end{align}

Applying \eqref{P_P}, we have
\begin{align*}
|R_{512}| \le & \int_{{\cal S}}\int_{-\epsilon_{n,s}^{M}}^{\epsilon_{n,s}^{M}}\bar{f}(x_0) \big| \Phi\Big[\frac{\sqrt{s}\big(1/2-{\mathbb P}\big(S_{n}(x_0^t) \geq 1/2\big)\big)}{\sqrt{{\mathbb P}\big(S_{n}(x_0^t)<1/2\big){\mathbb P}\big(S_{n}(x_0^t)\ge1/2\big)}} \Big] \nonumber\\
&\qquad 
-\Phi\big[2\sqrt{s}\big({\mathbb P}\big(S_{n}(x_0^t)<1/2\big)-1/2\big) \big] \big| dt d\textrm{Vol}^{d-1}(x_0)\nonumber \\
=& \int_{{\cal S}}\int_{-\epsilon_{n,s}^{M}}^{\epsilon_{n,s}^{M}}\bar{f}(x_0) \big| \Phi\Big[\frac{\sqrt{s}\big({\mathbb P}\big(S_{n}(x_0^t)<1/2\big)-1/2\big)}{\sqrt{{\mathbb P}\big(S_{n}(x_0^t)<1/2\big){\mathbb P}\big(S_{n}(x_0^t)\ge1/2\big)}} \Big] \nonumber\\
&\qquad 
-\Phi\big[2\sqrt{s}\big({\mathbb P}\big(S_{n}(x_0^t)<1/2\big)-1/2\big) \big] \big| dt d\textrm{Vol}^{d-1}(x_0) \nonumber\\
\leq&  O\Big(\frac{(\log(s))^2}{\sqrt{s}}\Big)\int_{{\cal S}}\int_{-\epsilon_{n,s}^{M}}^{\epsilon_{n,s}^{M}}\bar{f}(x_0) dt d\textrm{Vol}^{d-1}(x_0) =o(t_n+\frac{1}{\sqrt{s}}s_n). \nonumber
\end{align*}
\begin{align*}
|R_{513}| \le & \int_{{\cal S}}\int_{-\epsilon_{n,s}^{M}}^{\epsilon_{n,s}^{M}}\frac{|\dot{\bar{f}}(x_0)^T\dot{\eta}(x_0)||t|}{\|\dot{\eta}(x_0)\|} \big| \Phi\Big[\frac{\sqrt{s}\big(1/2-{\mathbb P}\big(S_{n}(x_0^t) \geq 1/2\big)\big)}{\sqrt{{\mathbb P}\big(S_{n}(x_0^t)<1/2\big){\mathbb P}\big(S_{n}(x_0^t)\ge1/2\big)}} \Big] \nonumber\\
&\qquad 
-\Phi\big[2\sqrt{s}\big({\mathbb P}\big(S_{n}(x_0^t)<1/2\big)-1/2\big) \big] \big| dt d\textrm{Vol}^{d-1}(x_0)\nonumber \\
=& \int_{{\cal S}}\int_{-\epsilon_{n,s}^{M}}^{\epsilon_{n,s}^{M}}\frac{|\dot{\bar{f}}(x_0)^T\dot{\eta}(x_0)||t|}{\|\dot{\eta}(x_0)\|} \big| \Phi\Big[\frac{\sqrt{s}\big({\mathbb P}\big(S_{n}(x_0^t)<1/2\big)-1/2\big)}{\sqrt{{\mathbb P}\big(S_{n}(x_0^t)<1/2\big){\mathbb P}\big(S_{n}(x_0^t)\ge1/2\big)}} \Big] \nonumber\\
&\qquad 
-\Phi\big[2\sqrt{s}\big({\mathbb P}\big(S_{n}(x_0^t)<1/2\big)-1/2\big) \big] \big| dt d\textrm{Vol}^{d-1}(x_0) \nonumber\\
\leq&  O\big(\frac{(\log(s))^2}{\sqrt{s}}\big)\int_{{\cal S}}\int_{-\epsilon_{n,s}^{M}}^{\epsilon_{n,s}^{M}}\frac{|\dot{\bar{f}}(x_0)^T\dot{\eta}(x_0)||t|}{\|\dot{\eta}(x_0)\|} dt d\textrm{Vol}^{d-1}(x_0) =o(t_n^2+\frac{s_n^2}{s}). \nonumber
\end{align*}

Next, we decompose
\begin{align}
&\int_{{\cal S}} \int_{-\epsilon_{n,s}^{M}}^{\epsilon_{n,s}^{M}}\bar{f}(x_0) \big\{ \Phi\big[2\sqrt{s}\big({\mathbb P}\big(S_{n}(x_0^t)<1/2\big)-1/2\big) \big] \label{eq:decompose_514_ci_M}\\
&\qquad\qquad\qquad\qquad
- \indi{t<0} \big\} dt d\textrm{Vol}^{d-1}(x_0) \nonumber\\
=& \int_{{\cal S}} \int_{-\epsilon_{n,s}^{M}}^{\epsilon_{n,s}^{M}}\bar{f}(x_0) \big\{ \Phi\big[2\sqrt{s}\big(\Phi\big(\frac{1/2-\mu_{n}(x_0^t)}{\sigma_{n}(x_0^t)}\big)-1/2\big) \big]  \nonumber\\
&\qquad\qquad\qquad\qquad
- \indi{t<0} \big\} dt d\textrm{Vol}^{d-1}(x_0) +R_{514}. \nonumber
\end{align}
\begin{align}
&\int_{{\cal S}} \int_{-\epsilon_{n,s}^{M}}^{\epsilon_{n,s}^{M}}\frac{\dot{\bar{f}}(x_0)^T\dot{\eta}(x_0)t}{\|\dot{\eta}(x_0)\|} \big\{ \Phi\big[2\sqrt{s}\big({\mathbb P}\big(S_{n}(x_0^t)<1/2\big)-1/2\big) \big] \label{eq:decompose_515_ci_M}\\
&\qquad\qquad\qquad\qquad
- \indi{t<0} \big\} dt d\textrm{Vol}^{d-1}(x_0) \nonumber\\
=& \int_{{\cal S}} \int_{-\epsilon_{n,s}^{M}}^{\epsilon_{n,s}^{M}}\frac{\dot{\bar{f}}(x_0)^T\dot{\eta}(x_0)t}{\|\dot{\eta}(x_0)\|} \big\{ \Phi\big[2\sqrt{s}\big(\Phi\big(\frac{1/2-\mu_{n}(x_0^t)}{\sigma_{n}(x_0^t)}\big)-1/2\big) \big]  \nonumber\\
&\qquad\qquad\qquad\qquad
- \indi{t<0} \big\} dt d\textrm{Vol}^{d-1}(x_0) +R_{515}. \nonumber
\end{align}

Applying \eqref{R13_temp}, we have
\begin{align*}
|R_{514}|
\le&  \int_{{\cal S}} \int_{-\epsilon_{n,s}^{M}}^{\epsilon_{n,s}^{M}}\bar{f}(x_0)\Big|\Phi\Big[2\sqrt{s}\big({\mathbb P}\big(S_{n}(x_0^t)<1/2\big)-1/2\big) \Big] \\ &\qquad
-\Phi\Big[2\sqrt{s}\big(\Phi\big(\frac{1/2-\mu_{n}(x_0^t)}{\sigma_{n}(x_0^t)}\big)-1/2\big) \Big] \Big|  dt d\textrm{Vol}^{d-1}(x_0)\\
=&o((\log(s))^{-2}) \int_{{\cal S}} \int_{-\epsilon_{n,s}^{M}}^{\epsilon_{n,s}^{M}} \bar{f}(x_0) dt d\textrm{Vol}^{d-1}(x_0)
=o(t_n+\frac{1}{\sqrt{s}}s_n).
\end{align*}
\begin{align*}
|R_{515}|
\le&  \int_{{\cal S}} \int_{-\epsilon_{n,s}^{M}}^{\epsilon_{n,s}^{M}}\frac{|\dot{\bar{f}}(x_0)^T\dot{\eta}(x_0)|}{\|\dot{\eta}(x_0)\|}|t|\Big|\Phi\Big[2\sqrt{s}\big({\mathbb P}\big(S_{n}(x_0^t)<1/2\big)-1/2\big) \Big] \\ &\qquad
-\Phi\Big[2\sqrt{s}\big(\Phi\big(\frac{1/2-\mu_{n}(x_0^t)}{\sigma_{n}(x_0^t)}\big)-1/2\big) \Big] \Big|  dt d\textrm{Vol}^{d-1}(x_0)\\
=&o((\log(s))^{-2}) \int_{{\cal S}} \int_{-\epsilon_{n,s}^{M}}^{\epsilon_{n,s}^{M}} \frac{|\dot{\bar{f}}(x_0)^T\dot{\eta}(x_0)|}{\|\dot{\eta}(x_0)\|}|t| dt d\textrm{Vol}^{d-1}(x_0)
=o(\frac{s_n^2}{s}+t_n^2).
\end{align*}

Next, we decompose
\begin{align}
&\int_{{\cal S}} \int_{-\epsilon_{n,s}^{M}}^{\epsilon_{n,s}^{M}}\bar{f}(x_0) \big\{ \Phi\big[2\sqrt{s}\big(\Phi\big(\frac{1/2-\mu_{n}(x_0^t)}{\sigma_{n}(x_0^t)}\big)-1/2\big) \big] \label{eq:decompose_516_ci_M}\\
&\qquad\qquad\qquad
- \indi{t<0}
 \big\} dt d\textrm{Vol}^{d-1}(x_0) \nonumber\\
=& \int_{{\cal S}} \int_{-\epsilon_{n,s}^{M}}^{\epsilon_{n,s}^{M}}\bar{f}(x_0)\big\{\Phi\big[\frac{1/2-\mu_{n}(x_0^t)}{\sqrt{\pi/(2s)}\sigma_{n}(x_0^t)} \big] \nonumber\\
&\qquad\qquad\qquad
-\indi{t<0}
 \big\} dt d\textrm{Vol}^{d-1}(x_0) +R_{516}. \nonumber
\end{align}
\begin{align}
&\int_{{\cal S}} \int_{-\epsilon_{n,s}^{M}}^{\epsilon_{n,s}^{M}}\frac{\dot{\bar{f}}(x_0)^T\dot{\eta}(x_0)t}{\|\dot{\eta}(x_0)\|} \big\{ \Phi\big[2\sqrt{s}\big(\Phi\big(\frac{1/2-\mu_{n}(x_0^t)}{\sigma_{n}(x_0^t)}\big)-1/2\big) \big] \label{eq:decompose_517_ci_M}\\
&\qquad\qquad\qquad
- \indi{t<0}
 \big\} dt d\textrm{Vol}^{d-1}(x_0) \nonumber\\
=& \int_{{\cal S}} \int_{-\epsilon_{n,s}^{M}}^{\epsilon_{n,s}^{M}}\frac{\dot{\bar{f}}(x_0)^T\dot{\eta}(x_0)t}{\|\dot{\eta}(x_0)\|}\big\{\Phi\big[\frac{1/2-\mu_{n}(x_0^t)}{\sqrt{\pi/(2s)}\sigma_{n}(x_0^t)} \big] \nonumber\\
&\qquad\qquad\qquad
-\indi{t<0}
 \big\} dt d\textrm{Vol}^{d-1}(x_0) +R_{517}. \nonumber
\end{align}

Applying \eqref{Normal_CDF_at_0}, we have
\begin{align*}
|R_{516}|
\le&\int_{{\cal S}} \int_{-\epsilon_{n,s}^{M}}^{\epsilon_{n,s}^{M}}\bar{f}(x_0)\Big|\Phi\Big[2\sqrt{s}\big(\Phi\big(\frac{1/2-\mu_{n}(x_0^t)}{\sigma_{n}(x_0^t)}\big)-1/2\big)\Big]\\
&\qquad\qquad
-\Phi\Big[\frac{1/2-\mu_{n}(x_0^t)}{\sqrt{\pi/(2s)}\sigma_{n}(x_0^t)}\Big] \Big| dt d\textrm{Vol}^{d-1}(x_0) \\
\leq& o\Big(\frac{1}{\sqrt{s}}\Big)\int_{{\cal S}} \int_{-\epsilon_{n,s}^{M}}^{\epsilon_{n,s}^{M}}\bar{f}(x_0) dt d\textrm{Vol}^{d-1}(x_0) =o(t_n+\frac{1}{\sqrt{s}}s_n). 
\end{align*} 
\begin{align*}
|R_{517}|
\le&\int_{{\cal S}} \int_{-\epsilon_{n,s}^{M}}^{\epsilon_{n,s}^{M}}\frac{|\dot{\bar{f}}(x_0)^T\dot{\eta}(x_0)|}{\|\dot{\eta}(x_0)\|}|t|\Big|\Phi\Big[2\sqrt{s}\big(\Phi\big(\frac{1/2-\mu_{n}(x_0^t)}{\sigma_{n}(x_0^t)}\Big)-1/2\big)\Big]\\
&\qquad\qquad
-\Phi\Big[\frac{1/2-\mu_{n}(x_0^t)}{\sqrt{\pi/(2s)}\sigma_{n}(x_0^t)}\Big] \Big| dt d\textrm{Vol}^{d-1}(x_0) \\
\leq&  o\Big(\frac{1}{\sqrt{s}}\Big)\int_{{\cal S}} \int_{-\epsilon_{n,s}^{M}}^{\epsilon_{n,s}^{M}}\frac{|\dot{\bar{f}}(x_0)^T\dot{\eta}(x_0)|}{\|\dot{\eta}(x_0)\|}|t| dt d\textrm{Vol}^{d-1}(x_0)  =o(t_n^2+\frac{1}{s}s_n^2). 
\end{align*} 

Next, we decompose
\begin{eqnarray}
&&\int_{{\cal S}} \int_{-\epsilon_{n,s}^{M}}^{\epsilon_{n,s}^{M}} \bar{f}(x_0) \big\{\Phi\Big(\frac{1/2 - \mu_{n}(x_0^t)}{\sqrt{\pi/(2s)}\sigma_{n}(x_0^t)} \Big) - \indi{t<0}\big\} dt d\textrm{Vol}^{d-1}(x_0) \label{eq:decompose_518_ci_M} \\
&=& \int_{{\cal S}} \int_{-\epsilon_{n,s}^{M}}^{\epsilon_{n,s}^{M}}\bar{f}(x_0) \big\{\Phi\Big(\frac{-2t\|\dot{\eta}(x_0)\|- 2a(x_0)t_n}{\sqrt{\pi/(2s)}s_{n}} \Big) \nonumber \\
&&\qquad\qquad\qquad\qquad\qquad
- \indi{t<0}\big\} dt d\textrm{Vol}^{d-1}(x_0) + R_{518}. \nonumber
\end{eqnarray}
\begin{eqnarray}
&&\int_{{\cal S}} \int_{-\epsilon_{n,s}^{M}}^{\epsilon_{n,s}^{M}}\frac{\dot{\bar{f}}(x_0)^T\dot{\eta}(x_0)t}{\|\dot{\eta}(x_0)\|} \big\{\Phi\Big(\frac{1/2 - \mu_{n}(x_0^t)}{\sqrt{\pi/(2s)}\sigma_{n}(x_0^t)} \Big) \label{eq:decompose_519_ci_M} \\
&&\qquad\qquad\qquad\qquad
- \indi{t<0}\big\} dt d\textrm{Vol}^{d-1}(x_0)  \nonumber\\
&=& \int_{{\cal S}} \int_{-\epsilon_{n,s}^{M}}^{\epsilon_{n,s}^{M}} \frac{\dot{\bar{f}}(x_0)^T\dot{\eta}(x_0)t}{\|\dot{\eta}(x_0)\|} \big\{\Phi\Big(\frac{-2t\|\dot{\eta}(x_0)\|- 2a(x_0)t_n}{\sqrt{\pi/(2s)}s_{n}} \Big) \nonumber \\
&&\qquad\qquad\qquad\qquad
- \indi{t<0}\big\} dt d\textrm{Vol}^{d-1}(x_0) + R_{519}. \nonumber
\end{eqnarray}
Similar to bounding $R_{15}$ in \eqref{eq:decompose_15_re_M}, we have \begin{align*}
&\int_{-\epsilon_{n,s}^{M}}^{\epsilon_{n,s}^{M}} \bar{f}(x_0) \Big| \Phi\Big(\frac{1/2 - \mu_{n}(x_0^t)}{\sqrt{\pi/(2s)}\sigma_{n}(x_0^t)}\Big) - \Phi\Big(\frac{-2t\|\dot{\eta}(x_0)\|- 2a(x_0)t_n}{\sqrt{\pi/(2s)}s_n} \Big)\Big| dt \nonumber\\
=& \bar{f}(x_0) (s_{n,s}^{M})^2 \int_{-\epsilon_{n,s}^{M}/s_{n,s}^{M}}^{\epsilon_{n,s}^{M}/s_{n,s}^{M}}  \Big| \Phi\Big(\frac{1/2-\mu_{n}(x_0^{r^{M}s_{n,s}^{M}})}{\sqrt{\pi/(2s)}\sigma_{n}(x_0^{r^{M}s_{n,s}^{M}})}\Big) \nonumber \\
&\qquad\qquad\qquad\qquad\qquad
-\Phi\big(-2\|\dot{\eta}(x_0)\|(r^{M}-r^{M}_{x_0})\big)  \Big| dr^{M} \nonumber \\
\le & \bar{f}(x_0) (s_{n,s}^{M})^2 \Big[ \int_{|r^{M}|\le \epsilon t_n/s_{n,s}^{M}} dr^{M} \nonumber\\
&
+ \epsilon^2 \int_{-\infty}^{\infty} (|r^{M}|+t_n/s_{n,s}^{M})\phi(\|\dot{\eta}(x_0)\||r^{M}-r^{M}_{x_0}|) dr^{M}\Big]=o(\frac{s_n}{\sqrt{s}}+t_n). \nonumber 
\end{align*}
The inequality above leads to $R_{518}=o(s_n^2/s+t_n^2)$. Similarly, we have
\begin{align*}
&\int_{-\epsilon_{n,s}^{M}}^{\epsilon_{n,s}^{M}} \frac{|\dot{\bar{f}}(x_0)^T\dot{\eta}(x_0)||t|}{\|\dot{\eta}(x_0)\|} \Big| \Phi\Big(\frac{1/2 - \mu_{n}(x_0^t)}{\sqrt{\pi/(2s)}\sigma_{n}(x_0^t)}\Big) \\
&\qquad\qquad\qquad\qquad\qquad
- \Phi\Big(\frac{-2t\|\dot{\eta}(x_0)\|- 2a(x_0)t_n}{\sqrt{\pi/(2s)}s_n} \Big)\Big| dt \\
=& \frac{|\dot{\bar{f}}(x_0)^T\dot{\eta}(x_0)|}{\|\dot{\eta}(x_0)\|} (s_{n,s}^{M})^2 \int_{-\epsilon_{n,s}^{M}/s_{n,s}^{M}}^{\epsilon_{n,s}^{M}/s_{n,s}^{M}} |r^{M}| \big| \Phi\Big(\frac{1/2-\mu_{n}(x_0^{r^{M}s_{n,s}^{M}})}{\sqrt{\pi/(2s)}\sigma_{n}(x_0^{r^{M}s_{n,s}^{M}})}\Big) \nonumber \\
&\qquad\qquad\qquad\qquad\qquad
-\Phi\big(-2\|\dot{\eta}(x_0)\|(r^{M}-r^{M}_{x_0})\big)  \big| dr^{M} \nonumber \\
\le & \frac{|\dot{\bar{f}}(x_0)^T\dot{\eta}(x_0)|}{\|\dot{\eta}(x_0)\|} (s_{n,s}^{M})^2 \Big[ \int_{|r^{M}|\le \epsilon t_n/s_{n,s}^{M}} |r^{M}| dr^{M} \nonumber\\
&
+ \epsilon^2 \int_{-\infty}^{\infty} |r^{M}|(|r^{M}|+t_n/s_{n,s}^{M})\phi(\|\dot{\eta}(x_0)\||r^{M}-r^{M}_{x_0}|) dr^{M}\Big]= o(\frac{s_n^2}{s}+t_n^2). \nonumber 
\end{align*}
The inequality above leads to $R_{519}=o(s_n^2/s+t_n^2)$.

Combining \eqref{Taylor1_ci_M} - \eqref{eq:decompose_519_ci_M}, we have
\begin{align}
&\int_{{\cal S}} \int_{-\epsilon_n}^{\epsilon_n} {\bar f}(x_0^t) \big\{{\mathbb P}\big(S_{n,s}^{M}(x_0^t) < 1/2\big) - \indi{t<0}\big\} dt d\textrm{Vol}^{d-1}(x_0) \label{exp1_ci_M}\\
=& \int_{{\cal S}} \int_{-\epsilon_n}^{\epsilon_n} \bar{f}(x_0) \big\{\Phi\big(\frac{-2t\|\dot{\eta}(x_0)\|- 2a(x_0)t_n}{\sqrt{\pi/(2s)}s_n} \big) \nonumber\\
&\qquad\qquad\qquad\qquad
- \indi{t<0}\big\} dt d\textrm{Vol}^{d-1}(x_0)\nonumber\\
&+\int_{{\cal S}} \int_{-\epsilon_n}^{\epsilon_n} \frac{\dot{\bar{f}}(x_0)^T\dot{\eta}(x_0)t}{\|\dot{\eta}(x_0)\|} \big\{\Phi\big(\frac{-2t\|\dot{\eta}(x_0)\|- 2a(x_0)t_n}{\sqrt{\pi/(2s)}s_n} \big) \nonumber\\
&\qquad\qquad\qquad\qquad
- \indi{t<0}\big\} dt d\textrm{Vol}^{d-1}(x_0) + o(s_n/\sqrt{s}+t_n). \nonumber
\end{align}
By similar arguments, we have
\begin{align}
&\int_{{\cal S}} \int_{-\epsilon_n}^{\epsilon_n} {\bar f}(x_0^t) \big\{{\mathbb P}^2\big(S_{n,s}^{M}(x_0^t) < 1/2\big) - \indi{t<0}\big\} dt d\textrm{Vol}^{d-1}(x_0) \label{exp2_ci_M}\\
=& \int_{{\cal S}} \int_{-\epsilon_n}^{\epsilon_n} \bar{f}(x_0) \big\{\Phi^2\big(\frac{-2t\|\dot{\eta}(x_0)\|- 2a(x_0)t_n}{\sqrt{\pi/(2s)}s_n} \big) \nonumber\\
&\qquad\qquad\qquad\qquad
- \indi{t<0}\big\} dt d\textrm{Vol}^{d-1}(x_0) \nonumber\\
&+\int_{{\cal S}} \int_{-\epsilon_n}^{\epsilon_n} \frac{\dot{\bar{f}}(x_0)^T\dot{\eta}(x_0)t}{\|\dot{\eta}(x_0)\|} \big\{\Phi^2\big(\frac{-2t\|\dot{\eta}(x_0)\|- 2a(x_0)t_n}{\sqrt{\pi/(2s)}s_n} \big) \nonumber\\
&\qquad\qquad\qquad\qquad
- \indi{t<0}\big\} dt d\textrm{Vol}^{d-1}(x_0)+ o(s_n/\sqrt{s}+t_n).  \nonumber
\end{align}

Finally, after substituting $t=\sqrt{\pi/(2s)}us_n/2$ in \eqref{exp1_ci_M} and \eqref{exp2_ci_M}, we have, up to $o(s_n/\sqrt{s}+t_n)$ difference,
\begin{align*}
\frac{{\rm CIS}(\widehat{\phi}_{n,s}^{M})}{2} =&\sqrt{\frac{\pi}{2s}}\frac{s_n}{2}\int_{{\cal S}} \int_{-\infty}^{\infty} \bar{f}(x_0) \big\{\Phi\big[-\|\dot{\eta}(x_0)\|u-\frac{2a(x_0)t_n}{\sqrt{\pi/(2s)}s_n} \big] \nonumber\\
&\qquad\qquad\qquad\qquad
- \indi{u<0} \big\} du d\textrm{Vol}^{d-1}(x_0) \nonumber\\
&+ \frac{\pi s_n^2}{8s}\int_{{\cal S}} \int_{-\infty}^{\infty} \frac{\dot{\bar{f}}(x_0)^T\dot{\eta}(x_0)}{\|\dot{\eta}(x_0)\|} u \big\{\Phi\big[-\|\dot{\eta}(x_0)\|u-\frac{2a(x_0)t_n}{\sqrt{\pi/(2s)}s_n}\big]  \nonumber\\
&\qquad\qquad\qquad\qquad
-\indi{u<0}\big\} du d\textrm{Vol}^{d-1}(x_0) \nonumber\\
&- \sqrt{\frac{\pi}{2s}}\frac{s_n}{2}\int_{{\cal S}} \int_{-\infty}^{\infty} \bar{f}(x_0) \big\{\Phi^2\big[-\|\dot{\eta}(x_0)\|u-\frac{2a(x_0)t_n}{\sqrt{\pi/(2s)}s_n} \big] \nonumber\\
&\qquad\qquad\qquad\qquad
- \indi{u<0}\big\} du d\textrm{Vol}^{d-1}(x_0) \nonumber\\
& - \frac{\pi s_n^2}{8s}\int_{{\cal S}} \int_{-\infty}^{\infty} \frac{\dot{\bar{f}}(x_0)^T\dot{\eta}(x_0)}{\|\dot{\eta}(x_0)\|} u \big\{\Phi^2\big[-\|\dot{\eta}(x_0)\|u-\frac{2a(x_0)t_n}{\sqrt{\pi/(2s)}s_n} \big] \nonumber\\
&\qquad\qquad\qquad\qquad
-\indi{u<0}\big\} du d\textrm{Vol}^{d-1}(x_0) \nonumber\\
=& I + II - III - IV.
\end{align*}

According to Lemma \ref{lemma:G}, we have
\begin{align*}
I - III =& \left [\int_{\cal S} \frac{{\bar{f}}(x_0)}{2\sqrt{\pi}\|\dot{\eta}(x_0)\|}  d\textrm{Vol}^{d-1}(x_0) \right] \sqrt{\frac{\pi}{2s}}s_n
= \frac{B_3}{2}\sqrt{\frac{\pi}{2s}}s_n,   \\
II - IV =& - \Big[\int_{\cal S} \frac{\dot{\bar{f}}(x_0)^T\dot{\eta}(x_0)a(x_0)}{2\sqrt{\pi}(\|\dot{\eta}(x_0)\|)^3}  d\textrm{Vol}^{d-1}(x_0)\Big] \sqrt{\frac{\pi}{2s}}s_nt_n = \frac{B_4}{2}\sqrt{\frac{\pi}{2s}}s_nt_n. 
\end{align*}
The desirable result is obtained by noting that $t_ns_n/\sqrt{s} = o(s_n/\sqrt{s}+t_n)$. This completes the proof of \eqref{eq:M-DNN_ci} in Theorem \ref{thm:DNN_ci}.

{\it Part 2}: Next, we will prove \eqref{eq:W-DNN_ci} in Theorem \ref{thm:DNN_ci}. We use similar notations as those in Section \ref{sec:pf_thm:W-DNN_re}. Similar to \eqref{eq:CIS_proof_M}, we have
\begin{align*}
{\rm CIS}(\widehat{\phi}_{n,s}^{W})/2
=&  \int_{{\cal R}} \big\{ {\mathbb P}(S_{n,s}^{W}(x) < 1/2 ) - \indi{\eta(x)<1/2}\big\} d\bar{P}(x) \nonumber\\
& -  \int_{{\cal R}} \big\{  {\mathbb P}^2(S_{n,s}^{W}(x) < 1/2) - \indi{\eta(x)<1/2} \big\}d \bar{P}(x).
\end{align*}

Next, we organize our proof in three steps similar to Section \ref{sec:pf_thm:M-DNN_re}. 

{\it Step 1}: Similar to {\it Step 1} in Section \ref{sec:pf_thm:W-DNN_re}, we have, uniformly for $\bw_n\in W_{n,\beta}$,
\begin{align*} 
&\int_{{\cal R}\cap{{\cal S}^{\epsilon_n}}} \big\{ {\mathbb P}(S_{n,s}^{W}(x) < 1/2) - \indi{\eta(x)<1/2} \big\} d \bar{P}(x)  \nonumber\\
=&\int_{{\cal S}} \int_{-\epsilon_n}^{\epsilon_n} {\bar f}(x_0^t) \big\{{\mathbb P}\big(S_{n,s}^{W}(x_0^t) < 1/2\big)
- \indi{t<0}\big\} dt d\textrm{Vol}^{d-1}(x_0)\{1+o(1)\},\nonumber\\
&\int_{{\cal R}\cap{{\cal S}^{\epsilon_n}}} \big\{ {\mathbb P}^2(S_{n,s}^{W}(x) < 1/2) - \indi{\eta(x)<1/2} \big\} d \bar{P}(x)  \nonumber\\
 =&\int_{{\cal S}} \int_{-\epsilon_n}^{\epsilon_n} {\bar f}(x_0^t) \big\{{\mathbb P}^2\big(S_{n,s}^{W}(x_0^t) < 1/2\big)  
- \indi{t<0}\big\} dt d\textrm{Vol}^{d-1}(x_0)\{1+o(1)\}.  \nonumber
\end{align*}

{\it Step 2}: Bound the contribution to CIS from ${\cal R}\backslash {\cal S}^{\epsilon_n}$. We have 
\begin{eqnarray}
&&\sup_{\bw_n\in W_{n,\beta}} \int_{{\cal R}\backslash {\cal S}^{\epsilon_n}} \big\{ {\mathbb P}\big(S_{n,s}^{W}(x) < 1/2\big) - \indi{\eta(x)<1/2} \big\} d \bar{P}(x) \label{bound1_ci_W}\\
&=& o(s_n^2/s+t_n^2), \nonumber\\
&&\sup_{\bw_n\in W_{n,\beta}} \int_{{\cal R}\backslash {\cal S}^{\epsilon_n}} \big\{ {\mathbb P}^2\big(S_{n,s}^{W}(x) < 1/2\big) - \indi{\eta(x)<1/2} \big\} d \bar{P}(x) \label{bound2_ci_W}\\
\quad &=& o(s_n^2/s+t_n^2). \nonumber
\end{eqnarray}
$(\ref{bound1_ci_W})$ holds by replacing $P^\circ(x)$ by $\bar{P}(x)$ in {\it Step 2} of Section \ref{sec:pf_thm:W-DNN_re}. Furthermore, $(\ref{bound2_ci_W})$ holds since
\begin{eqnarray*}
&&\big|{\mathbb P}^2\big(S_{n,s}^{W}(x) < 1/2\big) - \indi{\eta(x)<1/2}\big| \\ 
&\le& 2\big|{\mathbb P}\big(S_{n,s}^{W}(x) < 1/2\big) - \indi{\eta(x)<1/2}\big|.
\end{eqnarray*}

{\it Step 3}: In the end, we will show
\begin{eqnarray*}
&&\int_{{\cal S}} \int_{-\epsilon_n}^{\epsilon_n} {\bar f}(x_0^t) \big\{{\mathbb P}\big(S_{n,s}^{W}(x_0^t) < 1/2\big) - \indi{t<0}\big\} dt d\textrm{Vol}^{d-1}(x_0) \nonumber\\
&&-\int_{{\cal S}} \int_{-\epsilon_n}^{\epsilon_n} {\bar f}(x_0^t) \big\{{\mathbb P}^2\big(S_{n,s}^{W}(x_0^t) < 1/2\big) - \indi{t<0}\big\} dt d\textrm{Vol}^{d-1}(x_0) \nonumber\\
&=& \frac{1}{2}B_3 \frac{s_n}{\sqrt{s}} +  o(\frac{s_n}{\sqrt{s}}+t_n).
\end{eqnarray*}

Taylor expansion leads to
$$
\bar{f}(x_0^t) = \bar{f}(x_0) + \dot{\bar{f}}(x_0)^T\frac{\dot{\eta}(x_0)}{\|\dot{\eta}(x_0)\|} t + o(t).
$$
Hence, 
\begin{eqnarray}
&&\int_{{\cal S}} \int_{-\epsilon_n}^{\epsilon_n} {\bar f}(x_0^t) \big\{{\mathbb P}\big(S_{n,s}^{W}(x_0^t) < 1/2\big) - \indi{t<0}\big\} dt d\textrm{Vol}^{d-1}(x_0) \label{Taylor1_ci_W}\\
\quad &=& \big[ \int_{{\cal S}} \int_{-\epsilon_n}^{\epsilon_n} \bar{f}(x_0) \big\{{\mathbb P}\big(S_{n,s}^{W}(x_0^t) < 1/2\big) \nonumber - \indi{t<0}\big\} dt d\textrm{Vol}^{d-1}(x_0) \nonumber\\
&&+\int_{{\cal S}} \int_{-\epsilon_n}^{\epsilon_n} \frac{\dot{\bar{f}}(x_0)^T\dot{\eta}(x_0)t}{\|\dot{\eta}(x_0)\|} \big\{{\mathbb P}\big(S_{n,s}^{W}(x_0^t) < 1/2\big) \nonumber\\
&&\qquad\qquad\qquad\qquad\qquad
- \indi{t<0}\big\} dt d\textrm{Vol}^{d-1}(x_0) \big] \{1+o(1) \}. \nonumber
\end{eqnarray}

Next, we decompose
\begin{eqnarray}
&&\int_{{\cal S}} \int_{-\epsilon_n}^{\epsilon_n} {\bar f}(x_0) \big\{{\mathbb P}\big(S_{n,s}^{W}(x_0^t) < 1/2\big) - \indi{t<0}\big\} dt d\textrm{Vol}^{d-1}(x_0) \label{eq:decompose_521_ci_W}\\
&=& \int_{{\cal S}} \int_{-\epsilon_n}^{\epsilon_n} \bar{f}(x_0) \big\{\Phi\Big(\frac{1/2 - \mu_{n,s}(x_0^t)}{\sigma_{n,s}(x_0^t)} \Big) - \indi{t<0}\big\} dt d\textrm{Vol}^{d-1}(x_0) + R_{521}, \nonumber
\end{eqnarray}
\begin{align}
&\int_{{\cal S}} \int_{-\epsilon_n}^{\epsilon_n} \frac{\dot{\bar{f}}(x_0)^T\dot{\eta}(x_0)t}{\|\dot{\eta}(x_0)\|} \big\{{\mathbb P}\big(S_{n,s}^{W}(x_0^t) < 1/2\big) \label{eq:decompose_522_ci_W} \\
 &\qquad\qquad\qquad\qquad\qquad\qquad
-\indi{t<0}\big\} dt d\textrm{Vol}^{d-1}(x_0) \nonumber \\
=& \int_{{\cal S}} \int_{-\epsilon_n}^{\epsilon_n} \frac{\dot{\bar{f}}(x_0)^T\dot{\eta}(x_0)t}{\|\dot{\eta}(x_0)\|} \big\{\Phi\Big(\frac{1/2 - \mu_{n,s}(x_0^t)}{\sigma_{n,s}(x_0^t)} \Big)  \nonumber  \\
 &\qquad\qquad\qquad\qquad\qquad\qquad
- \indi{t<0}\big\} dt d\textrm{Vol}^{d-1}(x_0) + R_{522}. \nonumber
\end{align}

Applying \eqref{BE_thm_re_W} and \eqref{BE_greater_re_W}, we have
\begin{align*}
& \int_{-\epsilon_n}^{\epsilon_n} {\bar f}(x_0)   \Big| {\mathbb P}\big(S_{n,s}^{W}(x_0^t) < 1/2\big) - \Phi\Big(\frac{1/2 - \mu_{n,s}(x_0^t)}{\sigma_{n,s}(x_0^t)} \Big)\Big| dt \\
\leq & \int_{|t| \leq c_{33} t_n}   \frac{ c_{32}\sum_{l=1}^N w_{Nl}^3}{s_{n,s}^3} {\bar f}(x_0)  dt  \\
& + \int_{c_{33} t_n \leq |t| \leq \epsilon_n}   \frac{ c_{32}\sum_{l=1}^N w_{Nl}^3}{s_{n,s}^3} \frac{{\bar f}(x_0) }{1+ c_{34}^3 |t|^3/s_{n,s}^3}  dt  \\
\leq &  \frac{c_{32}\sum_{i=1}^n w_{ni}^3}{\sqrt{s}s_n^3} \int_{|t| \leq c_{33} t_n} {\bar f}(x_0)  dt \\
& + \frac{c_{32}\sum_{i=1}^n w_{ni}^3}{\sqrt{s}s_n^3} \int_{c_{33} t_n \leq |t| \leq \epsilon_n}  \frac{{\bar f}(x_0)}{ c_{34} \sqrt{s} |t|/s_n}  dt  =o(\frac{s_n}{\sqrt{s}}+t_n).
\end{align*}
The inequality above leads to $|R_{521}|=o(s_n/\sqrt{s}+t_n)$.  Similarly,
\begin{align*} 
&\int_{-\epsilon_n}^{\epsilon_n}   \frac{\big|\dot{\bar{f}}(x_0)^T\dot{\eta}(x_0)\big|}{\|\dot{\eta}(x_0)\|} |t| \Big| {\mathbb P}\big(S_{n,s}^{W}(x_0^t) < 1/2\big) - \Phi\Big(\frac{1/2 - \mu_{n,s}(x_0^t)}{\sigma_{n,s}(x_0^t)} \Big)\Big| dt  \\
\leq&  \frac{c_{32}\sum_{l=1}^N w_{Nl}^3}{s_{n,s}^3} \int_{|t| \leq c_{33} t_n}   \frac{\big|\dot{\bar{f}}(x_0)^T\dot{\eta}(x_0)\big|}{\|\dot{\eta}(x_0)\|} |t|  dt  \\
& + \frac{c_{32}\sum_{l=1}^N w_{Nl}^3}{s_{n,s}^3} \int_{c_{33} t_n \leq |t| \leq \epsilon_n} \frac{\big|\dot{\bar{f}}(x_0)^T\dot{\eta}(x_0)\big||t|}{\|\dot{\eta}(x_0)\|(1+ c_{34}^3 |t|^3/s_{n,s}^3)}  dt \\
\leq&  \frac{c_{32}\sum_{i=1}^n w_{ni}^3}{\sqrt{s}s_n^3} \int_{|t| \leq c_{33} t_n}  \frac{|\dot{\bar{f}}(x_0)^T\dot{\eta}(x_0)|}{\|\dot{\eta}(x_0)\|} |t|  dt \\
& + \frac{c_{32}\sum_{i=1}^n w_{ni}^3}{\sqrt{s}s_n^3} \int_{c_{33} t_n \leq |t| \leq \epsilon_n} \frac{|\dot{\bar{f}}(x_0)^T\dot{\eta}(x_0)||t|}{\|\dot{\eta}(x_0)\| c_{34}^2 |t|^2/s_{n,s}^2}  dt 
=o(\frac{s_n^2}{s}+t_n^2).
\end{align*}
The inequality above leads to $|R_{522}|=o(s_n^2/s+t_n^2)$.

Next, we decompose
\begin{eqnarray}
&&\int_{{\cal S}} \int_{-\epsilon_n}^{\epsilon_n} \bar{f}(x_0) \big\{\Phi\Big(\frac{1/2 - \mu_{n,s}(x_0^t)}{\sigma_{n,s}(x_0^t)} \Big) - \indi{t<0}\big\} dt d\textrm{Vol}^{d-1}(x_0) \label{eq:decompose_523_ci_W}\\
&=& \int_{{\cal S}} \int_{-\epsilon_n}^{\epsilon_n} \bar{f}(x_0) \big\{\Phi\Big(\frac{-2t\|\dot{\eta}(x_0)\|- 2a(x_0)t_n}{s_n/\sqrt{s}} \Big) \nonumber\\
&&\qquad\qquad\qquad\qquad\qquad
- \indi{t<0}\big\} dt d\textrm{Vol}^{d-1}(x_0) + R_{523}, \nonumber
\end{eqnarray}
\begin{eqnarray}
&&\int_{{\cal S}} \int_{-\epsilon_n}^{\epsilon_n} \frac{\dot{\bar{f}}(x_0)^T\dot{\eta}(x_0)t}{\|\dot{\eta}(x_0)\|} \big\{\Phi\Big(\frac{1/2 - \mu_{n,s}(x_0^t)}{\sigma_{n,s}(x_0^t)} \Big)\label{eq:decompose_524_ci_W}\\
&&\qquad\qquad\qquad\qquad\qquad
- \indi{t<0}\big\} dt d\textrm{Vol}^{d-1}(x_0) \nonumber \\
&&=\int_{{\cal S}} \int_{-\epsilon_n}^{\epsilon_n} \frac{\dot{\bar{f}}(x_0)^T\dot{\eta}(x_0)t}{\|\dot{\eta}(x_0)\|} \big\{\Phi\Big(\frac{-2t\|\dot{\eta}(x_0)\|- 2a(x_0)t_n}{s_n/\sqrt{s}} \Big)\nonumber\\
&&\qquad\qquad\qquad\qquad\qquad
 - \indi{t<0}\big\} dt d\textrm{Vol}^{d-1}(x_0) + R_{524}. \nonumber
\end{eqnarray}

Similar to bounding $R_{32}$ in \eqref{eq:decompose_32_re_W}, we have 
\begin{align*}
&\int_{-\epsilon_n}^{\epsilon_n} \bar{f}(x_0) \Big| \Phi\big(\frac{1/2 - \mu_{n,s}(x_0^t)}{\sigma_{n,s}(x_0^t)}\big)- \Phi\big(\frac{-2t\|\dot{\eta}(x_0)\|- 2a(x_0)t_n}{s_n/\sqrt{s}} \big)\Big| dt \nonumber\\
=& \bar{f}(x_0)s_{n,s} \int_{-\epsilon_n/s_{n,s}}^{\epsilon_n/s_{n,s}}  \big| \Phi\big(\frac{1/2-\mu_{n}(x_0^{r^{W}s_{n,s}})}{s^{-1/2}\sigma_{n}(x_0^{r^{W}s_{n,s}})}\big)   \\
&\qquad\qquad\qquad\qquad
-\Phi\big(-2\|\dot{\eta}(x_0)\|(r^{W}-r^{W}_{x_0})\big)  \big| dr^{W} \nonumber\\
\le & \bar{f}(x_0) s_{n,s}\Big[ \int_{|r^{W}|\le \epsilon t_n/s_{n,s}} dr^{W} \nonumber\\
&+ \epsilon^2 \int_{-\infty}^{\infty} (|r^{W}|+t_n/s_{n,s})\phi(\|\dot{\eta}(x_0)\||r^{W}-r^{W}_{x_0}|) dr^{W} \Big] = o(\frac{s_n}{\sqrt{s}} + t_n ).\nonumber
\end{align*}
The inequality above leads to $R_{523}=o(s_n/\sqrt{s}+t_n)$. Similarly,
\begin{align*}
& \int_{-\epsilon_n}^{\epsilon_n} \frac{|\dot{\bar{f}}(x_0)^T\dot{\eta}(x_0)||t|}{\|\dot{\eta}(x_0)\|} \Big| \Phi\big(\frac{1/2 - \mu_{n,s}(x_0^t)}{\sigma_{n,s}(x_0^t)}\big) -\Phi\big(\frac{-2t\|\dot{\eta}(x_0)\|- 2a(x_0)t_n}{s_n/\sqrt{s}} \big)\Big| dt \nonumber\\
=& \frac{|\dot{\bar{f}}(x_0)^T\dot{\eta}(x_0)|}{\|\dot{\eta}(x_0)\|}s_{n,s}^2 \int_{-\epsilon_n/s_{n,s}}^{\epsilon_n/s_{n,s}} |r^{W}| \big| \Phi\big(\frac{1/2-\mu_{n}(x_0^{r^{W}s_{n,s}})}{s^{-1/2}\sigma_{n}(x_0^{r^{W}s_{n,s}})}\big) \nonumber\\
&\qquad\qquad\qquad\qquad\qquad
-\Phi\big(-2\|\dot{\eta}(x_0)\|(r^{W}-r^{W}_{x_0})\big)  \big| dr^{W} \\
\le & \frac{|\dot{\bar{f}}(x_0)^T\dot{\eta}(x_0)|}{\|\dot{\eta}(x_0)\|} s_{n,s}^2\Big[   \int_{|r^{W}|\le \epsilon t_n/s_{n,s}} |r^{W}| dr^{W}  \nonumber\\
&+ \epsilon^2 \int_{-\infty}^{\infty}|r^{W}| (|r^{W}|
+t_n/s_{n,s})\phi(\|\dot{\eta}(x_0)\||r^{W}-r^{W}_{x_0}|) dr^{W} \Big] = o(s_n^2/s+t_n^2).\nonumber
\end{align*}
The inequality above leads to $R_{524}=o(s_n^2/s+t_n^2)$.

Combining \eqref{Taylor1_ci_W}--\eqref{eq:decompose_524_ci_W}, we have
\begin{align}
&\int_{{\cal S}} \int_{-\epsilon_n}^{\epsilon_n} {\bar f}(x_0^t) \big\{{\mathbb P}\big(S_{n,s}^{W}(x_0^t) < 1/2\big)- \indi{t<0}\big\} dt d\textrm{Vol}^{d-1}(x_0) \label{exp1_ci_W}\\
=& \int_{{\cal S}} \int_{-\epsilon_n}^{\epsilon_n} \bar{f}(x_0) \big\{\Phi\big(\frac{-2t\|\dot{\eta}(x_0)\|- 2a(x_0)t_n}{s_n/\sqrt{s}} \big)  \nonumber\\
&\qquad\qquad\qquad\qquad-  \indi{t<0}\big\} dt d\textrm{Vol}^{d-1}(x_0)\nonumber\\
&+\int_{{\cal S}} \int_{-\epsilon_n}^{\epsilon_n} \frac{\dot{\bar{f}}(x_0)^T\dot{\eta}(x_0)t}{\|\dot{\eta}(x_0)\|} \big\{\Phi\big(\frac{-2t\|\dot{\eta}(x_0)\|- 2a(x_0)t_n}{s_n/\sqrt{s}} \big) \nonumber\\
&\qquad\qquad\qquad\qquad
- \indi{t<0}\big\} dt d\textrm{Vol}^{d-1}(x_0) + o(s_n/\sqrt{s}+t_n). \nonumber
\end{align}

By similar arguments, we have
\begin{align}
&\int_{{\cal S}} \int_{-\epsilon_n}^{\epsilon_n} {\bar f}(x_0^t) \big\{{\mathbb P}^2\big(S_{n,s,\bw_n}(x_0^t) < 1/2\big)- \indi{t<0}\big\} dt d\textrm{Vol}^{d-1}(x_0) \label{exp2_ci_W}\\
=& \int_{{\cal S}} \int_{-\epsilon_n}^{\epsilon_n} \bar{f}(x_0) \big\{\Phi^2\big(\frac{-2t\|\dot{\eta}(x_0)\|- 2a(x_0)t_n}{s_n/\sqrt{s}} \big) \nonumber\\
&\qquad\qquad\qquad\qquad
- \indi{t<0}\big\} dt d\textrm{Vol}^{d-1}(x_0) \nonumber\\
&+\int_{{\cal S}} \int_{-\epsilon_n}^{\epsilon_n} \frac{\dot{\bar{f}}(x_0)^T\dot{\eta}(x_0)t}{\|\dot{\eta}(x_0)\|} \big\{\Phi^2\big(\frac{-2t\|\dot{\eta}(x_0)\|- 2a(x_0)t_n}{s_n/\sqrt{s}} \big) \nonumber\\
&\qquad\qquad\qquad\qquad
- \indi{t<0}\big\} dt d\textrm{Vol}^{d-1}(x_0)+ o(s_n/\sqrt{s}+t_n). \nonumber
\end{align}

Finally, after substituting $t=us_n/(2\sqrt{s})$ in (\ref{exp1_ci_W}) and  (\ref{exp2_ci_W}), we have, up to $o(s_n/\sqrt{s}+t_n)$ difference,
\begin{align*}
\frac{{\rm CIS}(\widehat{\phi}_{n,s}^{W})}{2} 
=&\frac{s_n}{2\sqrt{s}}\int_{{\cal S}} \int_{-\infty}^{\infty} \bar{f}(x_0) \big\{\Phi\big(-\|\dot{\eta}(x_0)\|u-\frac{2a(x_0)t_n}{s_n/\sqrt{s}} \big) \nonumber\\
&\qquad\qquad\qquad\qquad
- \indi{u<0} \big\} du d\textrm{Vol}^{d-1}(x_0) \nonumber\\
&+ \frac{s_n^2}{4s}\int_{{\cal S}} \int_{-\infty}^{\infty} \frac{\dot{\bar{f}}(x_0)^T\dot{\eta}(x_0)}{\|\dot{\eta}(x_0)\|} u \big\{\Phi\big(-\|\dot{\eta}(x_0)\|u-\frac{2a(x_0)t_n}{s_n/\sqrt{s}} \big)  \nonumber\\
&\qquad\qquad\qquad\qquad
-\indi{u<0}\big\} du d\textrm{Vol}^{d-1}(x_0) \nonumber\\
&- \frac{s_n}{2\sqrt{s}}\int_{{\cal S}} \int_{-\infty}^{\infty} \bar{f}(x_0) \big\{\Phi^2\big(-\|\dot{\eta}(x_0)\|u-\frac{2a(x_0)t_n}{s_n/\sqrt{s}} \big) \nonumber\\
&\qquad\qquad\qquad\qquad
- \indi{u<0}\big\} du d\textrm{Vol}^{d-1}(x_0) \nonumber\\
& - \frac{s_n^2}{4s}\int_{{\cal S}} \int_{-\infty}^{\infty} \frac{\dot{\bar{f}}(x_0)^T\dot{\eta}(x_0)}{\|\dot{\eta}(x_0)\|} u \big\{\Phi^2\big(-\|\dot{\eta}(x_0)\|u-\frac{2a(x_0)t_n}{s_n/\sqrt{s}} \big) \nonumber\\
&\qquad\qquad\qquad\qquad
-\indi{u<0}\big\} du d\textrm{Vol}^{d-1}(x_0) \nonumber\\
=& I + II - III - IV.
\end{align*}

According to Lemma \ref{lemma:G}, we have
\begin{align*}
I - III =& \left (\int_{\cal S} \frac{{\bar{f}}(x_0)}{2\sqrt{\pi}\|\dot{\eta}(x_0)\|}  d\textrm{Vol}^{d-1}(x_0) \right) \frac{s_n}{\sqrt{s}}
= \frac{B_3}{2}\frac{1}{\sqrt{s}} s_n,   \\
II - IV =& - \left (\int_{\cal S} \frac{\dot{\bar{f}}(x_0)^T\dot{\eta}(x_0)a(x_0)}{2\sqrt{\pi}(\|\dot{\eta}(x_0)\|)^3}  d\textrm{Vol}^{d-1}(x_0)\right) \frac{s_n}{\sqrt{s}}t_n = \frac{B_4}{2}\frac{1}{\sqrt{s}}s_nt_n. 
\end{align*}
The desirable result is obtained by noting that $t_ns_n/\sqrt{s} = o(s_n/\sqrt{s}+t_n)$. This concludes the proof of \eqref{eq:W-DNN_ci} in Theorem \ref{thm:DNN_ci}. \hfill $\blacksquare$

\subsection{Proof of Corollary \ref{thm:DNN_BNN}}\label{sec:pf_thm:DNN_BNN}
From \cite{S12}, we have 
\begin{eqnarray*}
\frac{{\rm Regret}(\widehat{\phi}_{N,\bw_N^*})}{{\rm Regret}(\widehat{\phi}_{N,k*})} \rightarrow 2^{4/(d+4)}\Big(\frac{d+2}{d+4}\Big)^{(2d+4)/(d+4)}, \\ 
\frac{{\rm Regret}(\widehat{\phi}_{N,q*})}{{\rm Regret}(\widehat{\phi}_{N,k*})} \rightarrow 2^{-4/(d+4)}\Gamma(2+2/d)^{2d/(d+4)}.  
\end{eqnarray*}
Therefore, 
\begin{equation}
\frac{{\rm Regret}(\widehat{\phi}_{N,q*})}{{\rm Regret}(\widehat{\phi}_{N,\bw_N^*})} \rightarrow Q^{''}=2^{-8/(d+4)}\Gamma(2+2/d)^{2d/(d+4)}\Big(\frac{d+4}{d+2}\Big)^{\frac{2d+4}{d+4}}.   \label{eq:opt_BNN_OWNN}
\end{equation}
By \eqref{eq:opt_BNN_OWNN}, \eqref{eq:opt_M-DNN} and \eqref{eq:opt_W-DNN}, we have
\begin{eqnarray*}
\frac{{\rm Regret}(\widehat{\phi}_{n,s,\bw_n^*}^{M})}{{\rm Regret}(\widehat{\phi}_{N,q*})} \rightarrow \frac{Q}{Q^{''}}>1  \;\;\;&\mbox{and}&\;\;\;
\frac{{\rm Regret}(\widehat{\phi}_{n,s,\bw_n^\dag}^{W})}{{\rm Regret}(\widehat{\phi}_{N,q*})} \rightarrow \frac{1}{Q^{''}}<1.
\end{eqnarray*}

From \cite{SQC16}, we have 
\begin{eqnarray*}
\frac{{\rm CIS}(\widehat{\phi}_{N,\bw_N^*})}{{\rm CIS}(\widehat{\phi}_{N,k*})} \rightarrow 2^{2/(d+4)}\Big(\frac{d+2}{d+4}\Big)^{(d+2)/(d+4)}, \\
\frac{{\rm CIS}(\widehat{\phi}_{N,q*})}{{\rm CIS}(\widehat{\phi}_{N,k*})} \rightarrow 2^{-2/(d+4)}\Gamma(2+2/d)^{d/(d+4)}.
\end{eqnarray*}
Therefore, 
\begin{eqnarray*}
\frac{{\rm CIS}(\widehat{\phi}_{N,q*})}{{\rm CIS}(\widehat{\phi}_{N,\bw_N^*})} \rightarrow \sqrt{Q^{''}}.
\end{eqnarray*}
Plugging $\bw_N^*$ in the general CIS formula given in \cite{SQC16} and plugging $\bw_n^*$ and $\bw_n^{\dag}$ in Theorem \ref{thm:DNN_ci}, we have
\begin{eqnarray*}
\frac{{\rm CIS}(\widehat{\phi}_{n,s,\bw_n^*}^{M})}{{\rm CIS}(\widehat{\phi}_{N,\bw_N^*})} \longrightarrow \sqrt{Q} \;\;\;&\mbox{and}&\;\;\; 
\frac{{\rm CIS}(\widehat{\phi}_{n,s,\bw_n^{\dag}}^{W})}{{\rm CIS}(\widehat{\phi}_{N,\bw_N^*})} \longrightarrow 1.
\end{eqnarray*}
Therefore, 
\begin{eqnarray*}
\frac{{\rm CIS}(\widehat{\phi}_{n,s,\bw_n^*}^{M})}{{\rm CIS}(\widehat{\phi}_{N,q*})} \rightarrow \sqrt{\frac{Q}{Q^{''}}}>1  \;\;\;&\mbox{and}&\;\;\;
\frac{{\rm CIS}(\widehat{\phi}_{n,s,\bw_n^{\dag}}^{W})}{{\rm CIS}(\widehat{\phi}_{N,q*})} \rightarrow \sqrt{\frac{1}{Q^{''}}}<1. \blacksquare
\end{eqnarray*}

\subsection{Lemmas}\label{sec:main_lemma}
In this section, we provide some lemmas. 
\begin{itemize}
    \item Lemma \ref{lemma:Phi}--Lemma \ref{lemma:G} are used for proving Theorem \ref{thm:M-DNN_re}.
    \item Lemma \ref{alpha} is used for proving Corollary \ref{thm:opt_M-DNN}.
\end{itemize}

\begin{lemma}
\label{lemma:Phi}
When $x$ is close to $0$ enough, we have \begin{equation*}
\Phi(x)-1/2=\frac{1}{\sqrt{2\pi}}x+O(x^3), 
\end{equation*} 
where $\Phi(x)$ is the standard normal distribution function.
\end{lemma}
\noindent {Proof of Lemma \ref{lemma:Phi}:} When $x$ is close to $0$ enough, by Taylor expansion of $\Phi(x)$ at 0, we have
\begin{eqnarray*}
\Phi(x) &=& \Phi(0) + \Phi'(0) x +\frac{1}{2} \Phi''(0) x^2 + O( \Phi'''(0) x^3) \\
&=& \frac{1}{2} +  \frac{1}{\sqrt{2\pi}}x + O( x^3 ).  \hfill \blacksquare
\end{eqnarray*}

\begin{lemma}
\label{lemma:mean_value} 
For constant $a>0$, we have
\begin{eqnarray}
|\Phi(ax_1)-\Phi(ax_2)| \le (a/2)|x_1 - x_2|, \label{eq:mean_value} 
\end{eqnarray}
where $\Phi(x)$ is the standard normal distribution function.
\end{lemma}
\noindent {Proof of Lemma \ref{lemma:mean_value}:}  
If $x_1=x_2$, \eqref{eq:mean_value} holds obviously. 

If $x_1<x_2$, by mean value theorem, there exists $x_0 \in (x_1,x_2)$ such that 
\begin{eqnarray*}
\Phi(ax_1)-\Phi(ax_2)  = \frac{1}{\sqrt{2\pi}}e^{-(ax_0)^2/2}a(x_1 - x_2). \end{eqnarray*}
Therefore,
\begin{eqnarray*}
|\Phi(ax_1)-\Phi(ax_2)| = \frac{1}{\sqrt{2\pi}}\exp\Big(-\frac{(ax_0)^2}{2}\Big)a|x_1 - x_2| \le (a/2)|x_1 - x_2|. \end{eqnarray*}

Similary, we can derive \eqref{eq:mean_value} when $x_1>x_2$.  \hfill $\blacksquare$

\begin{lemma}
\label{lemma:feller} 
\citep{feller1942} For all $x>0$, we have
\begin{eqnarray*}
1-\Phi(x) = \int_x^\infty \frac{1}{\sqrt{2\pi}}e^{-t^2/2}dt \le \frac{1}{x}\frac{e^{-x^2/2}}{\sqrt{2\pi}}.
\end{eqnarray*}
\end{lemma}
\noindent {Proof of Lemma \ref{lemma:feller}:}  
\begin{eqnarray*}
\int_x^\infty \frac{1}{\sqrt{2\pi}}e^{-t^2/2}dt \le \int_x^\infty \frac{t}{x}\frac{1}{\sqrt{2\pi}}e^{-t^2/2}dt = \frac{1}{x}\frac{e^{-x^2/2}}{\sqrt{2\pi}}.  \hfill \blacksquare
\end{eqnarray*}

\begin{lemma}
\label{lemma:dot} 
For $x_0 \in {\cal S} $, we have
\begin{eqnarray*}
2\bar{f}(x_0) \Vert\dot{\eta}(x_0)\Vert = \Vert\dot{\psi}(x_0)\Vert \;\;{\rm and}\;\; \dot{\psi}(x_0)^T\dot{\eta}(x_0) = \Vert \dot{\eta}(x_0) \Vert  \Vert\dot{\psi}(x_0)\Vert.
\end{eqnarray*}
\end{lemma}
\noindent {Proof of Lemma \ref{lemma:dot}:}  
By $\eta = {\mathbb P}(Y=1|X=x) = \frac{\pi_1 f_1}{ \pi_1 f_1 + (1-\pi_1) f_0}$, we have 
$$
\dot{\eta} = \frac{\pi_1(1-\pi_1) (\dot{f_1}f_0-f_1\dot{f_0})}{ (\pi_1 f_1 + (1-\pi_1) f_0)^2}.
$$
For $x_0 \in {\cal S}$, $\pi_1 f_1(x_0) = (1-\pi_1) f_0(x_0) = \frac{1}{2} \bar{f}(x_0)$, we have 
\begin{eqnarray*}
\dot{\eta}(x_0) &=& \frac{\pi_1(1-\pi_1) (\dot{f_1}(x_0)f_0(x_0)-f_1(x_0)\dot{f_0}(x_0))}{ [\pi_1 f_1(x_0) + (1-\pi_1) f_0(x_0)]^2} \\
 &=& \frac{ 1/2(\pi_1\dot{f_1}(x_0)-(1-\pi_1)\dot{f_0}(x_0))}{\bar{f}(x_0)} = \frac{ \dot{\psi}(x_0)}{2\bar{f}(x_0)}. 
 \end{eqnarray*}
Therefore, 
\begin{align*}
2\bar{f}(x_0) \Vert\dot{\eta}(x_0)\Vert =& \Vert\dot{\psi}(x_0)\Vert\;\;{\rm and}\;\; \\ \dot{\psi}(x_0)^T\dot{\eta}(x_0) =& 2\bar{f}(x_0)\dot{\eta}(x_0)^T\dot{\eta}(x_0) = \Vert \dot{\eta}(x_0) \Vert  \Vert\dot{\psi}(x_0)\Vert. \hfill \blacksquare
\end{align*}

\begin{lemma}
\label{lemma:G} 
\citep{SQC16} For any distribution function $G$, constant $a$, and constant $b>0$, we have
\begin{align*}
&\int_{-\infty}^{\infty} \big\{G(-bu- a ) - \indi{u<0}\big\} du = -\frac{1}{b}\big\{ a + \int_{-\infty}^{\infty} t dG(t) \big\}, \nonumber\\
&\int_{-\infty}^{\infty} u \big\{G(-bu- a ) - \indi{u<0}\big\} du \nonumber\\ 
&\qquad\qquad\qquad
= \frac{1}{b^2}\big\{ \frac{1}{2}a^2 + \frac{1}{2}\int_{-\infty}^{\infty} t^2dG(t) + a\int_{-\infty}^{\infty} t dG(t) \big\}. \nonumber \hfill \blacksquare
\end{align*}
\end{lemma}

\begin{lemma}
\label{alpha}\citep{SQC16}
Given $\alpha_i=i^{1+2/d}-(i-1)^{1+2/d}$, we have
\begin{eqnarray}
&&(1+\frac{2}{d})(i-1)^{\frac{2}{d}}\le \alpha_i \le (1+\frac{2}{d})i^{\frac{2}{d}}, \label{ineq}\\
&&\sum_{j=1}^k \alpha_j^2 = \frac{(d+2)^2}{d(d+4)}k^{1+4/d} \big\{1+O(\frac{1}{k})\big\}. \label{sum_alpha} \blacksquare
\end{eqnarray} 
\end{lemma}

\end{document}